\documentclass[12pt,reqno]{amsart}

\usepackage[sort,numbers]{natbib}

\usepackage[T1]{fontenc}
\usepackage{enumitem}
\setlist[enumerate,1]{wide, label=\textnormal{\textbf{(\roman*)}}, labelindent=0pt, itemsep=.1in}
\usepackage{hyperref}
\hypersetup{
  colorlinks   = true,
  citecolor    = blue,
  linkcolor    = blue
} 

\usepackage{amsmath,amsfonts,amsthm,color,tikz, comment, xcolor, xfrac,amssymb}
\usepackage[normalem]{ulem}
\usepackage{mathrsfs}
\usepackage{mathtools}
\usepackage{pdfsync}
\usepackage[font={scriptsize}]{caption}
\usepackage{environ}

\usepackage[left=1in, right=1in, top=1.1in,bottom=1.1in]{geometry}
\newcommand{\vk}{{\varkappa}}
\newcommand{\chik}{{\chi}}
\newcommand{\kAM}{{k}}
\newcommand{\aK}{{a}}
\newcommand{\specm}{{\hat{\Lambda}}}
\newcommand{\tfF}{{F}}
\newcommand{\tfPsi}{{\Psi}}
\newcommand{\nrho}{{\rho}}
\newcommand{\nrhop}{{\rho'}}

\newcommand{\R}{\mathbb R}

\newcommand{\cq}{\mathcal Q}

\newcommand{\ka}{\kappa}
\newcommand{\la}{\lambda}

\newcommand{\ze}{\zeta}
 
\newcommand{\lp}{\left(}
\newcommand{\rp}{\right)}

\newtheorem{theorem}{Theorem}[section]

\newtheorem{assumption}[theorem]{Assumption}

\newtheorem{definition}[theorem]{Definition}

\newtheorem{lemma}[theorem]{Lemma}

\newtheorem{proposition}[theorem]{Proposition}

\theoremstyle{remark}
\newtheorem{remark}[theorem]{Remark}

\theoremstyle{remark}
\newtheorem{example}[theorem]{Example}

 \numberwithin{equation}{section}
    
\newcommand{\bean}{\begin{eqnarray*}}
\newcommand{\eean}{\end{eqnarray*}}
\newcommand{\ben}{\begin{enumerate}}
\newcommand{\een}{\end{enumerate}}
\newcommand{\beq}{\begin{equation}}
\newcommand{\eeq}{\end{equation}}

\begin{document}

\title[Invariant measures in the whole space]{Invariant measures for stochastic heat equations\\ with dissipative drift in the whole space}

\author[E. Ashkarian]{Estepan Ashkarian}
\address{E. Ashkarian: Department of Mathematics, 
Purdue University, West Lafayette}
\email{eashkari@purdue.edu}

\author[S. Tindel]{Samy Tindel}
\address{S. Tindel: Department of Mathematics, 
Purdue University, West Lafayette}
\email{stindel@purdue.edu}

\begin{abstract}
We study the long-time behavior of a nonlinear stochastic heat equation on the whole space $\R^{d}$, driven by a Gaussian noise that is white in time and spatially homogeneous. The equation features a reaction term $f$ satisfying a one-sided Lipschitz condition and a bounded, globally Lipschitz diffusion coefficient $\sigma$. Solutions are constructed in a scale of weighted continuous function spaces following the framework of Salins~\cite{S}. Under an additional dissipativity condition on $f$~-- which covers reaction terms of the form $f(u)=u-u^{3}$~-- we prove that the transition semigroup admits an invariant measure on some weighted space $C_{\rho}$. This constitutes the first step toward a complete ergodic theory for dissipative stochastic heat equations on unbounded domains, a setting that has so far remained largely open in the literature.
\end{abstract}
\maketitle
\section{Introduction}

Ergodic theory for stochastic partial differential equations (SPDEs) has been the subject of a rich  literature over the past three decades. A foundational framework was initiated by Da~Prato and Zabczyk in their monograph~\cite{DPZ96}, which provides a systematic treatment of invariant measures and ergodicity for infinite-dimensional stochastic systems. Since then, this theory has been extended to a wide range of models, including stochastic Navier--Stokes equations with highly degenerate noise~\cite{HM}, stochastic heat equations with reflection~\cite{Z01}, and equations with multiplicative noise exhibiting phase transitions in their ergodic behavior~\cite{CE,GL}. We refer to~\cite{COTX} for a review of those results.
Remarkably, however, the arguably most basic case (a dissipative stochastic heat equation on the whole space $\R^{d}$ with a multiplicative noise) has not been given a definitive treatment. The aim of the present paper is to fill this gap.

\smallskip

Our goal is to study the long-time behavior of the nonlinear stochastic heat equation (NSHE) of the form
\begin{multline}\label{NSHE}
    \begin{cases}
        \partial_{t}u(t,x)=\frac{1}{2} \Delta u(t,x) +f(u(t,x))+\sigma(u(t,x))\dot{W}(t,x), \quad t\in[0,T], x\in \R^{d}\\
        u(0,x)=u_{0}(x)
    \end{cases},
\end{multline}
where $\Delta$ is the Laplacian on $\R^{d}$, $f$ and $\sigma$ are measurable real-valued functions representing the reaction and diffusion coefficients respectively, and $\dot{W}$ is a centered Gaussian noise which is white in time and spatially homogeneous. The precise assumptions on $f$, $\sigma$, the noise, and the initial condition $u_{0}$ are stated in Section~\ref{subsec: NoiseModel}--\ref{subsec: EU}.

\subsection{Prior results}\label{subsec:prior}
The existence of an invariant measure for a stochastic PDE is classically obtained by the Krylov--Bogoliubov argument. Namely, one shows that the time averages of the law of the solution are tight, and that the transition semigroup is Feller. Then every limit point of these averages is an invariant measure. The references below, which concern equation \eqref{NSHE} and its relatives, mostly differ by the way tightness is obtained. Those different approaches can depend on the spatial domain, on the drift and on the noise. We organize them accordingly.
\begin{enumerate}
\item \emph{Bounded domains.} In the framework of Da~Prato and Zabczyk \cite{DPZ96,DZ}, the equation is posed on a separable Hilbert space $H$, and tightness follows from the compactness of the semigroup $S(t)$ generated by the linear part (see \cite[Section~6.1]{DPZ96}). For an equation on a bounded domain $\mathcal{O}\subset\R^{d}$, this compactness is readily available. Within this framework, Cerrai~\cite{C} proved the existence of an invariant measure for reaction-diffusion systems with a multiplicative noise and a locally Lipschitz, monotone reaction term of polynomial growth. Her noise is white in space and time for $d=1$, and colored in space for $d>1$. See also Brze\'{z}niak and G\c{a}tarek~\cite{BG} for related results in Banach spaces. Let us also mention that \cite[Section~6.5]{DPZ96} treats dissipative systems in Banach spaces without any compactness, by a contraction argument which applies on unbounded domains as well. It requires an additive noise and a strictly dissipative drift, so that a reaction term like $f(u)=u-u^{3}$ is excluded.
\item \emph{The whole space, with the dissipativity of the heat kernel.} On $\R^{d}$ the heat semigroup is no longer compact, and the only dissipativity of the Laplacian is the pointwise decay $t^{-d/2}$ of the heat kernel, which is strong enough to produce uniform in time bounds only when $d\geq 3$. The first results in this setting are due to Tessitore and Zabczyk~\cite{TZ}, who considered equation \eqref{NSHE} with $f\equiv 0$ in weighted $L^{2}$ spaces. Misiats, Stanzhytskyi and Yip~\cite{MSY16,MSY20} also work in weighted $L^{2}$ spaces, for $d\geq 3$ and a nuclear $Q$-Wiener process, with a bounded $\sigma$ and a drift term which is bounded or only weakly dissipative in $u$. More recently, Gu and Li~\cite{GL}, Chen and Eisenberg~\cite{CE}, as well as~\cite{COTX}, established the existence of invariant measures, and phase transitions, for the equation without drift, under trace-class or spectral conditions on the noise correlation. These three references also require $d\geq 3$. Crucially, in each of them the existence of an invariant measure holds only in a \emph{small noise} regime. That is, the Lipschitz constant $L_{\sigma}$ of the diffusion coefficient $\sigma$ must be small enough to compensate for the strength of the noise correlation. When the noise intensity is too large, the moments of the solution blow up as $t\to\infty$ and no invariant measure exists.
\item \emph{The whole space, with the dissipativity of the drift.} Eckmann and Hairer~\cite{EH} handled the stochastic Ginzburg--Landau equation on $\R$, that is \eqref{NSHE} with the strongly dissipative drift $f(u)=u(1-|u|^{2})$ and an additive noise which is smooth (in fact analytic) in space and damped at infinity. Their solutions live in the space of bounded uniformly continuous functions, and the invariant measure is obtained for a weaker, weighted topology on that space. Indeed, the analyticity of the solution, inherited from the noise, provides the compactness that the heat semigroup lacks. Closer to our setting, Assing and Manthey~\cite{AM} proved the existence of invariant measures on $\R^{d}$, in weighted $L^{p}$ spaces, for a multiplicative noise and a possibly unbounded Lipschitz coefficient $\sigma$, under the dissipativity condition
\begin{equation}\label{eq: AM-dissip}
  u\,f(u)\leq C-\kAM u^{2}, \qquad u\in\R,
\end{equation}
where the constant $\kAM>0$ has to exceed a threshold depending on the Lipschitz constant of $\sigma$ and on the covariance of the noise. Their noise is either cylindrical, in dimension $d=1$ only, or a $Q$-Wiener process with a \emph{nuclear} covariance (that is $\operatorname{tr}(Q)<\infty$) for $d\geq 1$. This is a stronger spatial regularity than our Dalang condition~\eqref{eq:strongDalang}, which allows non-trace-class noises such as Riesz kernels. Compared with~\cite{AM}, our noise thus need not be trace class, in any dimension, at the price of a bounded diffusion coefficient $\sigma$.
\end{enumerate}

\subsection{Our contribution}\label{subsec:contribution}
In this paper we consider equation~\eqref{NSHE} in any dimension $d\geq 1$, with a drift $f$ which is dissipative at infinity, and we prove the existence of an invariant measure. Let us describe our assumptions informally, and refer to Sections~\ref{subsec: NoiseModel} and~\ref{subsec: EU} for the precise statements.

The noise $\dot{W}$ is a centered Gaussian noise, white in time and spatially homogeneous, whose spatial covariance $\Lambda$ has a spectral measure $\specm$ satisfying the strong Dalang condition
\begin{equation}\label{eq:intro-dalang}
\int_{\R^{d}}\frac{\specm(d\xi)}{1+|\xi|^{2(1-\eta)}}<\infty, \qquad\text{for some } \eta\in(0,1).
\end{equation}
This is the condition under which the solution of~\eqref{NSHE} is H\"older continuous in space and time, and it is satisfied for instance by the Riesz kernels $\Lambda(dx)=|x|^{-\alpha}dx$ with $0<\alpha<\min\{2,d\}$. Let us stress that no trace-class or nuclear condition is imposed on the noise, in contrast with the whole-space references of Section~\ref{subsec:prior}. As for the coefficients, the diffusion coefficient $\sigma$ is Lipschitz and bounded, and the reaction term $f$ is continuous, of polynomial growth, and satisfies the one-sided Lipschitz condition $f(u_{2})-f(u_{1})\leq\kappa(u_{2}-u_{1})$ for $u_{1}<u_{2}$, together with the dissipativity condition
\begin{equation}\label{eq:intro-dissip}
\big(f(u+h)-f(u)\big)\,\mathrm{sign}(h)\leq -c_{2}|h|^{\nu}+c_{1}\big(1+|u|^{\nu}\big), \qquad u,h\in\R,
\end{equation}
for some $\nu>1$ and $c_{1},c_{2}>0$. Roughly speaking, relation~\eqref{eq:intro-dissip} says that $f$ behaves like $-|u|^{\nu-1}u$ at infinity. The typical example we have in mind is the Allen--Cahn type nonlinearity $f(u)=u-u^{3}$, which satisfies~\eqref{eq:intro-dissip} with $\nu=3$. Note that $\kappa>0$ is allowed, so that $f$ is not dissipative near the origin and the deterministic equation may have several equilibria.

The solution of~\eqref{NSHE} is sought in the weighted spaces of Salins~\cite{S}. Namely, for $\theta>0$ we consider the space $C_{\theta}$ of continuous functions $\phi$ on $\R^{d}$ such that $\phi(x)/(1+|x|^{\theta})$ vanishes at infinity, equipped with the norm
\begin{equation}\label{eq:intro-Ctheta}
|\phi|_{C_{\theta}}:=\sup_{x\in\R^{d}}\frac{|\phi(x)|}{1+|x|^{\theta}} .
\end{equation}
This space contains unbounded functions, and the initial condition is any element of $C_{\theta}$. Our main result can be summarized as follows.
\begin{theorem}\label{theorem: main}
Under the assumptions above, let $\theta>0$ and $u_{0}\in C_{\theta}$, and let $u(t,\cdot\,;u_{0})$ be the solution of~\eqref{NSHE} started from $u_{0}$. Then for every $\rho>\nu\theta$, where $\nu$ is the exponent in \eqref{eq:intro-dissip}, the transition semigroup of~\eqref{NSHE} on $C_{\rho}$ admits an invariant probability measure. More precisely, the time averages of the law of $u(t,\cdot\,;u_{0})$ are tight in $C_{\rho}$, and every weak limit of these averages is an invariant measure.
\end{theorem}
The precise statement is Theorem~\ref{thm:existence-invariant} in Section~\ref{sec:invariance}. The strategy of its proof is described in Section~\ref{subsec:methodology} below.

\subsection{Methodology}\label{subsec:methodology}
Our global route to Theorem~\ref{theorem: main} is the classical one, due to Krylov and Bogoliubov. Namely, we prove that the time averages of the law of the solution are tight, and that the transition semigroup is Feller, so that every limit point of these averages is an invariant measure. However, the weights of the state space make both steps considerably harder than in a bounded domain, and the difficulties they raise are more than technicalities. Let us highlight the ingredients of our proof which we believe to be new.
\begin{enumerate}
\item \emph{A well-posedness theory in $C_{\theta}$.} The transition semigroup has to be defined on a separable Banach space of functions, and $C_{\theta}$ is the natural candidate. The well-posedness result of~\cite{S} does not start the equation from an arbitrary element of $C_{\theta}$, since it requires moments of the initial condition which are uniform over the centers $x_{0}$ of the weights. Theorem~\ref{thm:extend} removes this requirement when $\sigma$ is bounded. Its proof relies on the deterministic solution map $\mathcal{M}$ of~\cite{S} (see Section~\ref{sec:extend-M}), which handles the drift pathwise and whose Lipschitz property~\eqref{eq:MLip} in weighted norms is used throughout the paper. Our existence-uniqueness result also hinges on a Picard iteration in a translate $\mathcal{Y}_{\tau}$ of the space of~\cite{S} (see~\eqref{eq:Yspace}), in which the initial condition does not need any moment.
\item \emph{Continuity in the initial condition and a weighted-center norm.} The Feller property on $C_{\theta}$ requires the continuity of $u_{0}\mapsto u(t,\cdot\,;u_{0})$ for the norm of $C_{\theta}$, and the contraction argument of Theorem~\ref{thm:extend} does not give it. Indeed, the stochastic convolution estimate of~\cite{S} takes unweighted moments as input, and its weighted version (Proposition~\ref{prop:weighted-conv}) necessarily loses some exponent of the weight. Namely, with an input of exponent $\theta_{0}$, the output lives in $C_{x_{0},\theta}$ for $\theta>\theta_{0}$ only. We absorb this loss in a new norm $\|\cdot\|_{\mathcal{X}^{\theta_{0}}_{T}}$ (see Definition~\ref{def:X-weighted}), which penalizes the centers $x_{0}$ far from the origin, and we obtain in Theorem~\ref{thm:continuity} the Lipschitz continuity of $u_{0}\mapsto u(t,\cdot\,;u_{0})$ from $C_{\theta_{0}}$ to $C_{\theta}$, for every $\theta_{0}<\theta$.
\item \emph{A Feller property between two exponents.} As a consequence, the transition semigroup $\mathcal{P}_{t}$ is not shown to be Feller on one space, but to satisfy $\mathcal{P}_{t}(C_{b}(C_{\theta}))|_{C_{\theta_{0}}}\subset C_{b}(C_{\theta_{0}})$ for $\theta_{0}<\theta$ (see Proposition~\ref{prop:feller-two}). We then show that the Krylov--Bogoliubov argument is robust enough to accommodate this weaker property. Namely, it is run on $C_{\rho}$ with test functions which are continuous on the larger space $C_{\rho'}$, and these functions still generate the Borel $\sigma$-algebra of $C_{\rho}$ (see Lemma~\ref{lem:two-spaces}), which is all that the monotone class argument requires (see Theorem~\ref{thm:existence-invariant}).
\item \emph{Uniform bounds and tightness.} The tightness comes from uniform in time moment bounds on the solution in $C_{x_{0},\theta}$. They are proved in Section~\ref{sec:bounds} by decomposing $u$ into a stochastic convolution with a penalization $\beta$, whose moments decay as $\beta\to\infty$, and a random PDE which is handled pathwise thanks to the dissipativity of $f$. Those bounds are upgraded to a weighted Besov bound, and the compact embedding of the Besov space into $C_{\rho}$ for $\rho>\nu\theta$ gives the tightness in $C_{\rho}$ (see Section~\ref{sec:tightness}). This is where the loss $\rho>\nu\theta$ in Theorem~\ref{theorem: main} comes from: the factor $\nu$ is due to the polynomial growth of $f$ in the Besov bound, and the strict inequality is due to the compact embedding.
\end{enumerate}
Summarizing, although our proof follows the classical scheme, the amount of ingredients which are specific to the whole space and to the weighted spaces is not negligible, and we believe that they are of independent interest.

\subsection{General notation}\label{subsec:general-notation}
We close this introduction with some notation used throughout the paper. For a random variable $X$ and $p\geq 1$, we write $\|X\|_{L^{p}(\Omega)}:=\mathbf{E}[|X|^{p}]^{1/p}$, and $\mathrm{Law}(X)$ denotes the law of $X$. The indicator function of a set $A$ is $\mathbf{1}_{A}$, the Dirac mass at a point $\phi$ is $\delta_{\phi}$, and $\mu_{n}\Rightarrow\mu$ stands for the weak convergence of probability measures. The Fourier transform of a function $\phi$ on $\R^{d}$ is
\begin{equation}\label{eq:fourier-def}
\mathcal{F}\phi(\xi):=\int_{\R^{d}}e^{-i\xi\cdot x}\phi(x)\,dx, \qquad \xi\in\R^{d},
\end{equation}
with the same notation for the Fourier transform of a tempered measure, and $\mathcal{S}(\R^{d})$ is the Schwartz space of rapidly decreasing smooth functions. For a Banach space $E$, we write $C_{b}(E)$ and $B_{b}(E)$ for the spaces of bounded continuous and bounded Borel functions on $E$, both equipped with the supremum norm (see Section~\ref{subsec:functions-banach}). Finally, for two nonnegative quantities $a$ and $b$, the relation $a\lesssim b$ means that $a\leq Cb$ for a constant $C$ which does not depend on the variables at hand, and $a\lesssim_{\lambda,\theta}b$ specifies that $C$ depends on $\lambda$ and $\theta$.

\subsection{Outline of the paper}
The paper is organized as follows. In Section~\ref{sec:setting} we introduce the weighted function spaces, the noise and the standing assumptions on the coefficients. Then we prove the well-posedness of \eqref{NSHE} in the state space $C_{\theta}$ (Theorem~\ref{thm:extend}), together with the Markov property of the solution. Section~\ref{sec:feller} is devoted to the dependence of the solution on its initial condition. It contains the weighted estimate for the stochastic convolution, the continuity theorem in the weighted-center norm, the measurability and Feller properties of the transition semigroup between two exponents. In Section~\ref{sec:bounds} we establish the uniform in time bounds on the solution, first in the weighted spaces $C_{x_{0},\theta}$ and then in weighted Besov spaces. Section~\ref{sec:exist} gathers these ingredients. We prove the tightness of the time averages of the law of the solution, and then the existence of an invariant measure, which is Theorem~\ref{thm:existence-invariant}. The proofs of the auxiliary lemmas, namely the elementary facts on the weighted spaces, the moment estimates on stochastic convolutions and the measurability statements, are gathered in Appendices~\ref{app:weights} to~\ref{app:measurability}, so that the main text keeps to the arguments which are specific to the whole space setting.

\section{Setting and well-posedness}\label{sec:setting}
In this section we set up the framework in which equation $(\ref{NSHE})$ is solved. The function spaces are introduced in Section~\ref{subsec: funcspace} and the noise in Section~\ref{subsec: NoiseModel}. Section~\ref{subsec: EU} gathers the standing assumptions on the coefficients and on the initial condition. The well-posedness result on which the whole paper rests is Theorem~\ref{thm:extend} in Section~\ref{sec:extend}. It is an extension of the main theorem of \cite{S}, and its proof occupies the rest of that section. The Markov property of the solution is discussed in Section~\ref{subsec: markov}.

\subsection{Function Spaces}\label{subsec: funcspace}
For the solution of equation~\eqref{NSHE} on an infinite domain, the natural state space is a set of weighted continuous functions. This subsection introduces all the function spaces that will be used throughout the paper. We proceed in three steps. We first introduce weighted spaces of continuous functions, among which the state space $C_{\theta}$ for~\eqref{NSHE}. Then we define weighted H\"{o}lder spaces, which are the main tool in our tightness argument (Section~\ref{subsec: tight}). Finally we consider spaces of functions defined on the Banach space $C_{\theta}$, which are needed to formulate the Feller property of the associated transition semigroup.

\subsubsection{Spaces of continuous functions}
We start by introducing a notation for some standard unweighted function spaces: we call $C_{b}=C_{b}(\R^d)$ the space of bounded continuous functions $\phi \colon \R^d\to \R$, equipped with the norm $|\cdot|_{C_{b}}$, where
\begin{equation}\label{eq:normBound}
    |\phi|_{C_{b}}:=\sup_{x\in\R^d}|\phi(x)|.
\end{equation}
Later on, we will deal with functions which depend on both space and time. Therefore it might be useful to introduce, for $T>0$, the space of bounded continuous functions $\psi\colon [0,T]\times\R^d\to \R$, denoted by $C_{b}([0,T]\times \R^d)$ and equipped with the norm
\begin{equation}
    |\psi|_{C_{b}}:=\sup_{t\in [0,T]}\sup_{x\in \R^d}|\psi(t,x)|.
\end{equation}
We also consider the space of functions vanishing at infinity $C_{0}=C_{0}(\R^d)$, a closed subspace of $C_{b}(\R^d)$, where
\begin{equation}\label{eq: vanishinfty}
    \lim_{|x|\to \infty}|\phi(x)|=0, \quad \text{ for every } \phi\in C_{0}.
\end{equation}
Similarly, if we wish to consider functions depending on both space and time, we have a closed subspace of $C_{b}([0,T]\times\R^d)$, denoted by $C_{0}([0,T]\times \R^d)$, where 
\begin{equation}
    \lim_{|x|\to \infty}\sup_{t\in [0,T]}|\psi(t,x)|=0, \quad \text{ for every } \psi\in C_{0}([0,T]\times \R^d).
\end{equation}

We now introduce our weights for infinite domains.
\begin{definition}
Let $x_{0}\in \R^d$ and $\theta\geq0$ be fixed. We introduce the weights $w_{x_{0}}^\theta$ and $\hat{w}_{x_{0}}^{\theta}$ defined by
\begin{align}
 \label{eq: weight1}   w_{x_{0}}^{\theta}(x):=(1+|x-x_{0}|^{\theta})^{-1},
 \quad\text{and}\quad
\hat{w}_{x_{0}}^{\theta}(x):=(1+|x-x_{0}|^{2})^{-\theta/2}.
\end{align}
\end{definition}
\begin{remark}\label{rmk: weightEquiv} The reason we consider the two weights in (\ref{eq: weight1}) is because in \cite{S}, the weighted spaces are defined using the weights $w_{x_{0}}^{\theta}$. However, in the sequel we will also need to differentiate the weights in $x$. That is why we define the weight $\hat{w}_{x_{0}}^{\theta}$. This is no issue at all given that they are equivalent: there exists a universal constant $C_{1}>0$ such that
    \begin{equation}\label{eq: universalC1}
        \frac{1}{C_{1}(1+|x-x_{0}|^{\theta})}\leq (1+|x-x_{0}|^{2})^{-\theta/2}\leq \frac{C_{1}}{1+|x-x_{0}|^{\theta}}, \quad x\in \R^d,
    \end{equation}
    or equivalently $\frac{1}{C_{1}}w_{x_{0}}^{\theta}(x)\leq \hat{w}_{x_{0}}^{\theta}(x)\leq C_{1}w_{x_{0}}^{\theta}(x)$.
\end{remark}
The weights have at most polynomial growth, so that moving their center, or comparing two exponents, only costs a polynomial factor. The following elementary inequalities will be used throughout the paper.
\begin{lemma}\label{lem:weight-elem}
Let $\theta\geq 0$ and $x_{0}\in\R^{d}$, let $w_{x_{0}}^{\theta}$ be the weight of \eqref{eq: weight1}, and set $c_{\theta}:=\max\{1,2^{\theta-1}\}$. Then the following holds true.
\begin{enumerate}
\item\label{it:weight-triangle} For all $x,y\in\R^{d}$ we have
\begin{equation}\label{eq:weight-triangle}
1+|y-x_{0}|^{\theta}\leq c_{\theta}\bigl(1+|x-x_{0}|^{\theta}\bigr)\bigl(1+|x-y|^{\theta}\bigr),
\qquad\text{that is}\qquad
\frac{w_{x_{0}}^{\theta}(x)}{w_{x_{0}}^{\theta}(y)}\leq c_{\theta}\bigl(1+|x-y|^{\theta}\bigr).
\end{equation}
\item\label{it:weight-two-exp} For all $0\leq\theta_{0}\leq\theta$ and $x\in\R^{d}$ we have
\begin{equation}\label{eq:weight-two-exp}
w_{x_{0}}^{\theta}(x)\leq 2\,w_{x_{0}}^{\theta_{0}}(x).
\end{equation}
\end{enumerate}
\end{lemma}
\begin{proof}
The triangle inequality and the elementary bound $(a+b)^{\theta}\leq c_{\theta}(a^{\theta}+b^{\theta})$, valid for $a,b\geq 0$, give $1+|y-x_{0}|^{\theta}\leq 1+c_{\theta}(|x-x_{0}|^{\theta}+|x-y|^{\theta})$, and the right hand side is bounded by the product in \eqref{eq:weight-triangle}, as one sees by expanding this product. Dividing by $1+|x-x_{0}|^{\theta}$ gives the second form of \eqref{eq:weight-triangle}. As for item \ref{it:weight-two-exp}, for $|x-x_{0}|\geq 1$ we have $|x-x_{0}|^{\theta}\geq|x-x_{0}|^{\theta_{0}}$, so that $w_{x_{0}}^{\theta}(x)\leq w_{x_{0}}^{\theta_{0}}(x)$, and for $|x-x_{0}|<1$ we have $w_{x_{0}}^{\theta}(x)\leq 1\leq 2w_{x_{0}}^{\theta_{0}}(x)$, since $1+|x-x_{0}|^{\theta_{0}}<2$.
\end{proof}

\noindent
With the weights defined, we are now ready to introduce weighted spaces of functions defined on the whole space $\R^d$.
\begin{definition}\label{def: weightSpace} For a fixed $x_{0}\in \R^d$ and a parameter $\theta$, let $w_{x_{0}}^{\theta}$ be the weight defined\\ by $(\ref{eq: weight1})$. Then the weighted space $C_{x_{0},\theta}(\R^d)$ is defined as 
\begin{equation}\label{eq: weightSpace}
    C_{x_{0},\theta}(\R^d):=\left\{\phi \in C(\R^d): \lim_{|x|\to \infty}|w_{x_{0}}^{\theta}(x)\phi(x)|=0 \right\}.
\end{equation}
We  equip this space with the norm
\begin{equation}\label{eq: wsnorm}
    |\phi|_{C_{x_{0},\theta}}:=\sup_{x\in \R^d}\left|w_{x_{0}}^{\theta}(x)\phi(x)\right|.
\end{equation}
\end{definition}

\begin{remark}\label{rmk: pointwiseNorm}
    A direct consequence of $(\ref{eq: wsnorm})$ is the pointwise bound: for any $\phi\in C_{x_{0},\theta}(\R^d)$,
    \begin{equation}\label{eq: pointwiseNorm}
        |\phi(x_{0})|\leq |\phi|_{C_{x_{0},\theta}},
    \end{equation}
    obtained by setting $x=x_{0}$ in $(\ref{eq: wsnorm})$, so that $w_{x_{0}}^{\theta}(x_{0})=1$.
\end{remark}

\noindent
As mentioned above, in order to deal with systems like $(\ref{NSHE})$ we need to consider functions of both time and space. Below we extend Definition \ref{def: weightSpace}
 to this context.
 \begin{definition}\label{def: weightSpace2}
     As in Definition \ref{def: weightSpace}, we consider the weight $w_{x_{0}}^{\theta}$ defined by $(\ref{eq: weight1})$, for a fixed couple $(x_{0},\theta)\in \R^d \times \R_{+}$. For  a finite horizon $T$, we define $C_{x_{0},\theta}([0,T]\times \R^d)$ as
     \begin{equation}\label{eq: weightSpace2}
       C_{x_{0},\theta}([0,T]\times \R^d):=\left\{\psi \in C([0,T]\times\R^d): \lim_{|x|\to \infty}\sup_{t\in[0,T]}\left|w_{x_{0}}^{\theta}(x)\psi(t,x)\right|=0\right\}.
     \end{equation}
     This space is equipped with the norm
     \begin{equation}\label{eq: wsnorm2}
|\psi|_{C_{x_{0},\theta}([0,T]\times\R^d)}:=\sup_{t\in [0,T]}\sup_{x\in \R^d}\frac{|\psi(t,x)|}{1+|x-x_{0}|^{\theta}}=\sup_{t\in[0,T]}|w_{x_{0}}^{\theta}\psi(t,\cdot)|_{C_{0}}.
     \end{equation}
 \end{definition}
 \begin{remark}
     In both Definition \ref{def: weightSpace} and Definition \ref{def: weightSpace2}, $x_{0}$ is called the center of the weight.
 \end{remark}
\begin{remark}\label{rem: equiv-norms}
    The spaces $C_{0}$, $C_{0}([0,T]\times\R^{d})$, $C_{x_{0},\theta}(\R^{d})$ and $C_{x_{0},\theta}([0,T]\times\R^{d})$ are separable Banach spaces. Furthermore, for different centers $x_{0}\neq x_{1}\in \R^d$, the norms $|\cdot|_{C_{x_{0},\theta}}$ and $|\cdot|_{C_{x_{1},\theta}}$ are equivalent. However, in \cite{S}, it was required to have uniform estimates over the centers.
\end{remark}

Among the spaces $C_{x_{0},\theta}(\R^{d})$ of Definition~\ref{def: weightSpace}, the one with center $x_{0}=0$ will play a special role. Namely, we set
\begin{equation}\label{eq:Ctheta}
C_{\theta}:=C_{0,\theta}(\R^{d}),\qquad |\phi|_{C_{\theta}}:=|\phi|_{C_{0,\theta}}=\sup_{x\in\R^{d}}\frac{|\phi(x)|}{1+|x|^{\theta}} .
\end{equation}
The space $C_{\theta}$ contains unbounded functions, for instance $\phi(x)=(1+|x|^{2})^{\theta'/2}$ with $0<\theta'<\theta$. Let us explain why it is chosen as the state space.
\begin{remark}\label{rmk:stateSpace}
In order to show the existence of an invariant measure, one has to define a transition semigroup $\mathcal{P}_{t}\colon B_{b}(E)\to B_{b}(E)$ on a separable Banach space $E$ of functions defined pathwise, and to show that it is Feller, that is $\mathcal{P}_{t}(C_{b}(E))\subset C_{b}(E)$. The framework of \cite{S} requires moments of the initial condition which are uniform over the centers $x_{0}$. The two spaces it suggests are $C_{b}(\R^{d})$, which is not separable, and the space of random fields with moments uniform over the centers used in \cite{S} (see Remark~\ref{rmk:initial-fixed}), which is not a space of functions. This is why we work in the fixed-center space $C_{\theta}$ of \eqref{eq:Ctheta}, which is a separable Banach space by Remark~\ref{rem: equiv-norms}. The price to pay is that the existence and uniqueness result of \cite{S} has to be extended to initial conditions in $C_{\theta}$. This is the object of Theorem~\ref{thm:extend} below.
\end{remark}

\subsubsection{Spaces of Besov and H\"older functions}
The weighted continuous spaces introduced above provide the correct state space for the solution. To establish tightness of the family of time-marginal laws, however, we need spaces with additional regularity. The following definitions introduce weighted H\"{o}lder and Besov-type spaces, and collect the compact embedding result that will be one of the main analytical ingredient in Section~\ref{subsec: tight}.

\begin{definition}\label{def: compactSpace}
     For a given $x_{0}\in \R^d$ and $\theta>0$, let $C_{x_{0},\theta}(\R^d)$ be the space introduced in \\Definition \ref{def: weightSpace}. Then, for $s\in(0,1)$, we set
    \begin{equation}\label{eq: compactSpace}
    \mathscr{C}_{x_{0},\theta}^{s}(\R^d):=\{\psi \in C(\R^d): |\psi|_{C_{x_{0},\theta}}+[\psi]_{x_{0},\theta,s}<\infty\},
    \end{equation}
    where the weighted semi-norm $[\psi]_{s;x_{0},\theta}$ is defined by
    \begin{equation}\label{eq: csseminorm}
        [\psi]_{s;x_{0},\theta}:=\sup_{0<|h|<1}\sup_{x\in \R^d}\frac{|\psi(x+h)-\psi(x)|\hat{w}^{\theta}_{x_{0}}(x)}{|h|^{s}}.
    \end{equation}
\end{definition}
\begin{remark}
    In the language of weighted Besov space (see \cite[Chapter 4]{ET}), one can notice that the spaces $C_{x_{0},\theta}(\R^d)$ and $\mathscr{C}^{s}_{x_{0},\theta}(\R^d)$ from Definition \ref{def: weightSpace2} and Definition \ref{def: compactSpace} overlap with the Besov spaces $\mathcal{B}^{0}_{\infty,\infty}(w_{x_{0}}^{\theta})$ and $\mathcal{B}_{\infty,\infty}^{s}(w_{x_{0}}^{\theta})$.
\end{remark}
Next, we define a weighted fractional Sobolev type space which will provide us a useful continuous embedding.
\begin{definition}\label{def: Besov}
    Let $x_{0}\in \R^d$, $p\geq1$ and $\theta>0$. We define the space $\mathcal{B}_{p,p}^{s}(\hat{w}_{x_{0}}^{\theta})$ as
    \begin{equation}
        \mathcal{B}_{p,p}^{s}(\hat{w}_{x_{0}}^{\theta}):=\left\{ f\in \mathcal{S}'(\R^d):|f|_{\mathcal{B}_{p,p}^{s}(\hat{w}_{x_{0}}^{\theta})}<\infty\right\},
    \end{equation}
    where 
    \begin{equation}\label{eq: Besov-ppnorm}
        |f|_{\mathcal{B}_{p,p}^{s}(\hat{w}_{x_{0}}^{\theta})}:=\|f\|_{L^{p}(\hat{w}_{x_{0}}^{\theta})}
        +\left(\int_{|h|\leq 1} dh 
        \int_{\R^d}\frac{|f(x+h)-f(x)|^{p}\left(w_{x_{0}}^{\theta}(x)\right)^p}{|h|^{d+sp}} \, dx\right)^\frac{1}{p}
    \end{equation}
\end{definition}
Now we state a continuous and compact embedding result which we will implement to establish the existence of an invariant measure in Section \ref{subsec: tight}. Its proof can be found in \cite[Chapter 4]{ET}.
\begin{lemma}\label{lem: compact}
    Let $x_{0}\in \R^d$, $0<\widetilde{\theta}\leq\hat{\theta}<\theta$, $s>0$ and $\widetilde{s}>s+\frac{d}{p}$. Then the following series of embedding holds true:
    \begin{equation}
    \mathcal{B}_{p,p}^{\widetilde{s}}(w_{x_{0}}^{\widetilde{\theta}})\xhookrightarrow[\text{cont.}]{} \mathscr{C}^{s}_{x_{0},{\hat\theta}}(\R^d)\xhookrightarrow[\text{compact}]{} C_{x_{0},\theta}(\R^d) .
    \end{equation}
\end{lemma}

\subsubsection{Functions on Banach spaces}\label{subsec:functions-banach}
So far we have described spaces of functions defined on $\R^{d}$. The proof of the Feller property for the transition semigroup requires us to work with functions defined on the infinite-dimensional state space $C_{\theta}$ of \eqref{eq:Ctheta} itself. We therefore introduce, for a general Banach space $E$, the following spaces of bounded and differentiable functions on $E$. Namely we denote by $C_{b}(E)$ the space of bounded continuous functions $\varphi \colon E \to \R$, equipped with the sup norm
\begin{equation}\label{eq:normCbE}|\varphi|_{C_{b}(E)}:=\sup_{h\in E}|\varphi(h)|.\end{equation} 
Furthermore, as in \cite[Section~5]{C}, we denote by $C_{b}^{1}(E)$ the space of Fr\'{e}chet differentiable functions $\varphi \colon E \to \R$ having a continuous and bounded derivative, equipped with the norm
\begin{equation}\label{eq:normCb1E}
|\varphi|_{C_{b}^{1}(E)}:=|\varphi|_{C_{b}(E)}+\sup_{h\in E}|D\varphi(h)|_{E^{*}},
\end{equation}
where $E^{*}$ is the dual space of $E$. 

\subsection{Noise Model}\label{subsec: NoiseModel}
This subsection examines the noise which is considered in $(\ref{NSHE})$. We will work with the multiparametric setting developed by Walsh \cite{W} and Dalang \cite{D}. The following assumption on the noise $\dot{W}$ is considered.
\begin{assumption}\label{assump:noise}
    In $(\ref{NSHE})$, the noise $\dot{W}$ will be a centered Gaussian spatially homogeneous noise which is white in time, defined on a complete probability space $(\Omega,\mathcal{F},\mathcal{F}_{t},\mathbf{P})$. The filtration $\mathcal{F}_{t}$ satisfies the usual conditions, the noise is adapted to it, and for every $s\geq 0$ the increments of $\dot{W}$ after time $s$ are independent of $\mathcal{F}_{s}$. For instance, $\mathcal{F}_{t}$ can be the usual augmentation of the filtration generated by $\dot{W}$ and by a $\sigma$-field $\mathcal{G}$ independent of $\dot{W}$, which will carry the initial condition (see Assumption~\ref{assump: initial-cond-fixed}). In other words, there exists a positive and positive definite  tempered measure $\Lambda$ on $\R^d$, such that one formally has
    \begin{equation}\label{noise}
        \mathbf{E}[\dot{W}(s,x)\dot{W}(t,y)]=\delta(t-s)\Lambda(x-y).
    \end{equation}
    In $(\ref{noise})$, $\delta$ refers to the Dirac measure. Furthermore, the Fourier transform of $\Lambda$ is referred to as the spectral measure $\specm=\mathcal{F}\Lambda$, and we assume that it satisfies the so-called strong Dalang condition: there exists $\eta\in (0,1)$ such that
    \begin{equation}\label{eq:strongDalang}
    \int_{\R^d}\frac{1}{1+|\xi|^{2(1-\eta)}}\specm(d\xi)<+\infty.
        \end{equation}
        Moreover, we assume that the measure $\Lambda$ has a corresponding density $\ell$, or otherwise stated $\Lambda(dx)=\ell(x)\,dx$. 
        \end{assumption}
\begin{remark}
One can also express~\eqref{noise} as a relation for covariances of Wiener integrals. Specifically, for continuous and rapidly decreasing functions $\psi$ and $\phi$ on $\R_{+}\times\R^{d}$, the covariance of the corresponding Wiener integrals is given by
        \begin{equation}\label{eq:covWiener}
            \mathbf{E}\!\left[W(\psi)\,W(\phi)\right]
            = \int_{0}^{\infty}\int_{\R^{d}}\int_{\R^{d}} \psi(s,x)\, \ell(x-y) \, \phi(s,y)\,dy\,dx\,ds.
        \end{equation}
   In addition, the relationship between $\specm$ and $\ell$ is the following:
        \begin{equation}
\int_{\R^d}\ell(x)\varphi(x)\,dx=\int_{\R^d}\mathcal{F}\varphi(\xi)\specm(d\xi), \qquad \text{for all } \varphi\in\mathcal{S}(\R^d).
        \end{equation}
        It is also an easy exercise to establish that for $\varphi,\psi\in \mathcal{S}(\R^d)$ we have
        \begin{equation}\label{eq: relation}
            \int_{\R^d}\int_{\R^d}\varphi(x)\ell(x-y)\psi(y)\,dy\,dx=\int_{\R^d}\mathcal{F}\varphi(\xi)\overline{\mathcal{F}\psi(\xi)}\specm(d\xi),
        \end{equation}
        which gives an alternative expression for~\eqref{eq:covWiener}.
                \end{remark}
Assumption \ref{assump:noise} encompasses many examples of singular noise of interest. We refer to e.g.~ \cite{COTX} for a review. In the current contribution we will single out Riesz Kernels, which provide a canonical case of heavy tailed and singular at $0$ correlation, as our main example of application.
\begin{example} Set $\Lambda(dx)=|x|^{-\alpha}dx$ for $\alpha \in(0,d)$. We let the reader verify, as an easy exercise, that the corresponding spectral measure defined in Assumption \ref{assump:noise} is given by $\specm(d\xi)=C|\xi|^{\alpha-d} d\xi$ for a constant $C=C_{d,\alpha}$. The strong Dalang condition \eqref{eq:strongDalang} is then satisfied if and only if $\alpha<2(1-\eta)$, that is for every $\alpha\in(0,\min\{2,d\})$ with $\eta$ small enough.
\end{example}

\subsection{Mild solutions and standing assumptions}\label{subsec: EU}
In this subsection we give the mild formulation of $(\ref{NSHE})$ and we list the assumptions on the coefficients $f$, $\sigma$ and on the initial condition. The assumptions on $f$ and $\sigma$ are those of \cite{S}, up to the boundedness of $\sigma$ and to the dissipativity of $f$ (Assumption~\ref{assump: dissip}), which is the additional hypothesis of this paper. The assumption on the initial condition is weaker than in \cite{S}, since it only involves the fixed-center space $C_{\theta}$ of \eqref{eq:Ctheta}.
\begin{definition}\label{def: heatkernel}
The \emph{heat kernel} on $\R^d$ is
\begin{equation}\label{eq: heatkernel}
    G(t,x) := (2\pi t)^{-d/2}\exp\!\left(-\frac{|x|^{2}}{2t}\right), \qquad t>0,\; x\in\R^d.
\end{equation}
We shall use the following properties of $G$, all of which follow from the explicit expression \eqref{eq: heatkernel}.
\begin{enumerate}
\item\label{it:G-moment} \emph{Mass and moments.} For $t>0$, $G(t,\cdot)$ is a probability density, and by the change of variables $x=t^{1/2}y$ we have, for every $\theta\geq 0$,
\begin{equation}\label{eq:G-moment}
\int_{\R^{d}}G(t,x)\,|x|^{\theta}\,dx=m_{\theta}\,t^{\theta/2},
\qquad\text{where}\qquad
m_{\theta}:=\int_{\R^{d}}G(1,y)\,|y|^{\theta}\,dy<\infty .
\end{equation}
In particular $G(t,\cdot)$ integrates every polynomial.
\item\label{it:G-semigroup} \emph{Semigroup property.} For $s,t>0$ and $x\in\R^{d}$,
\begin{equation}\label{eq:G-semigroup}
G(t+s,x)=\int_{\R^{d}}G(t,x-y)\,G(s,y)\,dy .
\end{equation}
\item\label{it:G-fourier} \emph{Fourier transform.} With $\mathcal{F}$ defined by \eqref{eq:fourier-def}, for $t>0$ and $\xi\in\R^{d}$,
\begin{equation}\label{eq:G-fourier}
\mathcal{F}\big(G(t,\cdot)\big)(\xi)=\exp\Big(-\frac{t|\xi|^{2}}{2}\Big).
\end{equation}
\item\label{it:G-deriv} \emph{Derivative bounds.} A standard computation gives the gradient bound
\begin{equation}\label{eq: gradG}
    |\nabla_{x} G(t,x)| \leq C_{d}\, t^{-(d+1)/2}, \qquad t>0,\; x\in\R^d,
\end{equation}
where $C_{d}$ depends only on $d$. This bound will be used via the mean value theorem in $(\ref{eq: MVTforH})$. Similarly, the time derivative satisfies
\begin{equation}\label{eq: dtG}
    |\partial_{t} G(t,x)| \leq C_{d}\, t^{-(d+2)/2}, \qquad t>0,\; x\in\R^d,
\end{equation}
and the Laplacian satisfies $\partial_{t}G=\frac{1}{2}\Delta G$ as well as
\begin{equation}\label{eq:laplaceG}
\Delta_{x}G(t,x)=\Big(\frac{|x|^{2}}{t^{2}}-\frac{d}{t}\Big)G(t,x),\qquad\text{so that}\qquad
|\Delta_{x}G(t,x)|\leq\frac{1}{t}\Big(d+\frac{|x|^{2}}{t}\Big)G(t,x).
\end{equation}
\end{enumerate}
\end{definition}

\begin{definition}\label{def: mild}
    A mild solution to $(\ref{NSHE})$ is defined to be an adapted random field solution to the integral equation
    \begin{multline}\label{eq: mild}
  u(t,x)=\int_{\R^d}G(t,x-y)u_{0}(y)dy+\int_{0}^{t}\int_{\R^d}G(t-s,x-y)f(u(s,y))dyds\\
+\int_{0}^{t}\int_{\R^d}G(t-s,x-y)\sigma(u(s,y))W(ds,dy),
       \end{multline}
       where $G$ is the heat kernel defined in $(\ref{eq: heatkernel})$. In~(\ref{eq: mild}), the stochastic integral with respect to the noise $W$ is interpreted in the Dalang--Walsh sense (see \cite{D,W}).
\end{definition}
Let us now recall the semi--Lipschitz assumption and the growth condition on the drift coefficient $f$, which is used in \cite{S} to prove existence and uniqueness of a solution. This hypothesis is compatible with the dissipativity Assumption~\ref{assump: dissip} below, on which the ergodic properties rest.
\begin{assumption}\label{assump:reac}
    In $(\ref{NSHE})$, the reaction term $f\colon \R \to \R$ is continuous, and has the property that there exists $\kappa \in \R$ such that for any $u_1<u_2$ 
    \begin{equation}\label{assump: f1}
        f(u_2)-f(u_1)\leq \ka (u_2-u_1).
    \end{equation}
 Additionally, we impose on $f$ a very mild growth condition. Namely we suppose  that there exists a constant $K>0$ and $\nu\geq1$ such that
        \begin{equation}\label{assump: growth}
            |f(u)|\leq K(1+|u|^{\nu}) , \quad\text{ for every } u\in \R.
        \end{equation}
\end{assumption}
\begin{remark}\label{rmk: growth}
 In \cite{S}, the growth condition is assumed to be exponential,
    \begin{equation}\label{eq: expGrowth}
        |f(u)|\leq K \exp(K|u|^{\nu}), \quad \text{ for every } u\in \R
    \end{equation}
   in contrast to the polynomial growth condition in $(\ref{assump: growth})$. In this work, we restrict ourselves to the polynomial setting, as it is more convenient for the tightness arguments that follow. It also allows any weight exponent $\theta>0$ in the state space $C_{\theta}$ of \eqref{eq:Ctheta}, instead of the range $\theta<2/\nu$ of \cite{S} (see item \ref{it:M-eq} of Proposition~\ref{prop:Mmap} below).
\end{remark}
\begin{remark}\label{rmk: fequiv}
    As mentioned in {\cite[Proposition~2.4]{S}}, relation $(\ref{assump: f1})$ is equivalent to
    \begin{equation}\label{eq: equivf1}
        (f(u_{1})-f(u_{2}))\mathrm{sign}(u_{1}-u_{2})\leq \kappa |u_{1}-u_{2}|,\quad \text{for all } u_{1},u_{2}\in \R,
    \end{equation}
    where we define the sign function by
    \begin{equation}
        \mathrm{sign}(u)=\begin{cases}
            -1\quad \,\text{if}\,\,\,\, u<0\\
            0 \quad \,\,\,\,\hspace{0.1cm}\text{if}\,\,\,\, u=0\\
            1 \quad \hspace{0.04cm}\,\,\,\,\, \text{if}\,\,\,\, u>0.
        \end{cases}
    \end{equation}
\end{remark}
Assumption~\ref{assump:reac} is all that is needed for well-posedness. The ergodic properties of Sections~\ref{sec:bounds} and~\ref{sec:exist} rely on an additional dissipativity property of $f$, which we state now. It quantifies the one-sided Lipschitz condition \eqref{eq: equivf1} by a negative term of order $|h|^{\nu}$, and it is satisfied by the reaction term $f(u)=u-u^{3}$ with $\nu=3$.
\begin{assumption}\label{assump: dissip}
    The reaction term $f\colon \R \to \R$ satisfies
    \begin{equation}\label{eq: dissip}
        (f(u+h)-f(u))\,\mathrm{sign}(h)\leq -c_{2}|h|^{\nu}+c_{1}(1+|u|^{\nu}), \qquad u,h\in\R,
    \end{equation}
    where $c_{1},c_{2}>0$ and $\nu>1$ is the exponent of the growth condition \eqref{assump: growth} on $f$.
\end{assumption}
Let us stress that Assumption~\ref{assump: dissip} is not used before Section~\ref{sec:bounds}. In particular, the well-posedness Theorem~\ref{thm:extend} and the continuity results of Section~\ref{sec:feller} hold under Assumption~\ref{assump:reac} alone. Note also that the exponent $\nu$ in \eqref{eq: dissip} is the one of the growth condition \eqref{assump: growth}. This is a restriction, since taking $u=0$ in \eqref{eq: dissip} shows that the dissipation exponent never exceeds the growth exponent, while the converse may fail. It means that $f$ behaves like $-|u|^{\nu-1}u$ at infinity in both directions, as an odd-degree polynomial with negative leading coefficient does. The same convention is adopted in \cite[Theorem~5.5]{C}, where the growth exponent $m>1$ and the dissipativity condition (5.17) share the exponent $m$.

\noindent
Next, for the diffusion coefficient $\sigma$ in \eqref{NSHE}, we assume Lipschitz continuity and boundedness. Note that the boundedness is not a standard requirement. We explain in Remark~\ref{rmk:extend-compare} why it is needed.
\begin{assumption}\label{assump:diff}
    In $(\ref{NSHE})$, the diffusion coefficient $\sigma\colon \R \to \R$ is assumed to be globally Lipschitz and bounded with $L_{\sigma}:=\operatorname{Lip}(\sigma)$ and $M_{\sigma}:=\sup_{x\in \R^d}|\sigma(x)|$.
\end{assumption}

We now turn to the initial condition. In \cite{S}, the initial condition is a random field with moments uniform over the centers of the weights, namely $\sup_{x_{0}\in\R^{d}}\mathbf{E}[|u_{0}|^{p}_{C_{x_{0},\theta}}]<\infty$ (see the assumptions of \cite[Theorem~2.6]{S}). As explained in Remark~\ref{rmk:stateSpace}, we wish to start the equation from any element of the fixed-center space $C_{\theta}$. We thus consider the following assumption.
\begin{assumption}\label{assump: initial-cond-fixed}
The initial condition $\xi\colon\Omega\to C_{\theta}$ is an $\mathcal{F}_{0}$-measurable random variable with values in the space $C_{\theta}$ of \eqref{eq:Ctheta}, independent of the noise $\dot{W}$. Together with the weight exponent $\theta$, we fix an integrability exponent $p$. It will be the exponent of the moment estimates on the stochastic convolution (Proposition~\ref{prop:StochConvS}) and on the solution (item~\ref{it:extend-moment} of Theorem~\ref{thm:extend}). The parameters $\theta$ and $p$ satisfy
\begin{equation}\label{assump: theta-p}
\theta>0 \quad \text{and}\quad p>\max \left\{\frac{(1+\theta)(d+1)}{\theta},\frac{2(d+1)}{\eta} \right\},
\end{equation}
where $\eta$ is from the strong Dalang condition $(\ref{eq:strongDalang})$, and $d$ is the spatial dimension. The lower bound on $p$ in \eqref{assump: theta-p} is the one required in \cite[Theorem~2.6]{S}.
\end{assumption}
\begin{remark}\label{rmk:initial-fixed}
Let us compare Assumption~\ref{assump: initial-cond-fixed} with the initial condition of \cite{S} recalled above, and comment on the parameters $\theta$ and $p$.
\begin{enumerate}
\item For a deterministic initial condition, Assumption~\ref{assump: initial-cond-fixed} only says that $\xi\in C_{\theta}$, so that unbounded functions are allowed. For instance, take $0<\theta'<\theta$ and $\xi(x)=(1+|x|^{2})^{\theta'/2}$. Then $\xi\in C_{\theta}$, while $\sup_{x_{0}}|\xi|_{C_{x_{0},\theta}}\geq\sup_{x_{0}}|\xi(x_{0})|=\infty$ by the pointwise bound \eqref{eq: pointwiseNorm}. In contrast, a deterministic initial condition satisfying the assumption of \cite{S} has to be bounded, by the same pointwise bound.
\item No moment condition is imposed on $\xi$ in Assumption~\ref{assump: initial-cond-fixed}. Moments of $\xi$ will only be needed for moment bounds on the solution, see item \ref{it:extend-moment} of Theorem~\ref{thm:extend}.
\item As for the exponent $p$, the lower bound in \eqref{assump: theta-p} is what makes the stochastic convolution estimate of Proposition~\ref{prop:StochConvS} (quoted from \cite[Theorem~4.3]{S}) available, and that estimate is the main ingredient in the proof of Theorem~\ref{thm:extend}.
\item As for the exponent $\theta$, let us mention that \cite{S} also requires $\theta<2/\nu$, where $\nu$ is the exponent of the growth condition \eqref{assump: growth}. This restriction is due to the exponential growth \eqref{eq: expGrowth} assumed there, and item \ref{it:M-eq} of Proposition~\ref{prop:Mmap} below shows that it is not needed under the polynomial growth condition \eqref{assump: growth}.
\end{enumerate}
\end{remark}

\begin{remark}[Convention on the exponent $p$]\label{rmk:pstar}
The solution given by Theorem~\ref{thm:extend} below does not depend on $p$. Indeed, the uniqueness in item~\ref{it:extend-exist} of that theorem holds among the adapted mild solutions with paths in $C([0,T];C_{\theta})$, a class which does not involve $p$. Only the moment bound \eqref{eq:solution_bound} depends on $p$, and it holds for every $p$ satisfying \eqref{assump: theta-p} and \eqref{eq:initial_bound}. Moreover, since $\sigma$ is bounded, the stochastic convolutions of this paper have moments of all orders, so that $p$ can be taken as large as we wish. Accordingly, we say that a statement holds \emph{for $p\geq p^{*}$} when it holds for every $p$ larger than a finite exponent $p^{*}$, which depends only on $d$, $\eta$, $\nu$ and the weight exponents (such as $\theta$, $\theta_{0}$ or $\rho$) of the statement, and neither on the initial condition, nor on time, nor on the centers $x_{0}$. The value of $p^{*}$ changes from one statement to another, and is at least the lower bound in \eqref{assump: theta-p}, so that Theorem~\ref{thm:extend} applies. Within a proof, finitely many conditions of the form $p>c$, with $c$ depending on the same quantities, are invoked through the words ``since $p\geq p^{*}$'' without being made explicit. The constants in the estimates are allowed to depend on $p$.
\end{remark}

\subsection{Existence and uniqueness}\label{sec:extend}
This subsection is devoted to the well-posedness of $(\ref{NSHE})$ under the assumptions of Section~\ref{subsec: EU}. The main result of \cite{S}, namely \cite[Theorem~2.6]{S}, does not cover initial conditions in $C_{\theta}$, since it requires the moment condition $\sup_{x_{0}\in\R^{d}}\mathbf{E}[|u_{0}|^{p}_{C_{x_{0},\theta}}]<\infty$ of Remark~\ref{rmk:initial-fixed}, which is uniform over all centers $x_{0}\in\R^{d}$. We show below that, when $\sigma$ is bounded (Assumption~\ref{assump:diff}), equation \eqref{eq: mild} is well posed for every $\mathcal{F}_{0}$-measurable initial condition with values in $C_{\theta}$. The strategy is similar to the one of \cite[Section~6]{S}. Namely we perform a Picard iteration on the stochastic convolution, combined with the deterministic solution map of \cite[Section~5]{S}. The boundedness of $\sigma$ is what allows us to run this iteration without any moment uniform over the centers (see Remark~\ref{rmk:extend-compare}).

We proceed in five steps. In Section~\ref{sec:extend-tools}, we gather two elementary tools, a conjugation identity for the weight and the subdifferential of the sup norm, which are also used in Section~\ref{sec:bounds}. In Section~\ref{sec:extend-M}, we recall the deterministic mapping $\mathcal{M}$ and its properties from \cite[Section~5]{S}, and we extend its uniqueness property to the fixed-center space. In Section~\ref{sec:extend-fixed}, we record two elementary properties of the space $C_{\theta}$. In Section~\ref{sec:extend-conv}, we recall the stochastic convolution estimate from \cite[Theorem~4.3]{S}. Finally, Section~\ref{sec:extend-proof} contains the proof of the main result of this subsection, which we now state.
\begin{theorem}\label{thm:extend}
Let Assumptions~\ref{assump:noise}, \ref{assump:reac}, \ref{assump:diff} and~\ref{assump: initial-cond-fixed} hold, and let $C_{\theta}$ be the space defined in \eqref{eq:Ctheta}. Then the following holds true.
\begin{enumerate}
\item\label{it:extend-exist} For every $T>0$ there exists a unique adapted random field $\{u(t,x;\xi)\}_{(t,x)\in[0,T]\times\R^{d}}$ with paths in $C([0,T];C_{\theta})$ which is a mild solution of \eqref{eq: mild} with initial condition $u_{0}=\xi$, in the sense of Definition~\ref{def: mild}.
\item\label{it:extend-moment} Assume moreover that
\begin{equation}\label{eq:initial_bound}
\mathbf{E}\big[|\xi|^{p}_{C_{\theta}}\big]<\infty .
\end{equation}
Then there exists a constant $C>0$, depending on $T,p,\theta,d,\eta,\specm,\kappa,|f(0)|,L_{\sigma}$ and $M_{\sigma}$ but not on $\xi$, such that
\begin{equation}\label{eq:solution_bound}
\mathbf{E}\Big[\sup_{t\in[0,T]}|u(t,\cdot\,;\xi)|^{p}_{C_{\theta}}\Big]\leq C\Big(1+\mathbf{E}\big[|\xi|^{p}_{C_{\theta}}\big]\Big)<\infty .
\end{equation}
\end{enumerate}
\end{theorem}

\begin{remark}\label{rmk:extend-compare}
Let us compare Theorem~\ref{thm:extend} with the result of \cite{S}, namely \cite[Theorem~2.6]{S}.
\begin{enumerate}
\item\label{it:compare-unbdd} The initial condition is only assumed to lie in $C_{\theta}$, with no uniformity over the centers (see Remark~\ref{rmk:initial-fixed}). Accordingly, the moment bound \eqref{eq:solution_bound} is only stated for the center $x_{0}=0$, while the bound of \cite[Theorem~2.6]{S} is uniform over the centers.
\item\label{it:compare-sigma} The price to pay is the boundedness of $\sigma$ in Assumption~\ref{assump:diff}, i.e. $M_{\sigma}<\infty$. It is used in Steps~1, 2 and 5 of the proof below. It guarantees that the stochastic convolutions involved satisfy bounds uniform over the centers, whatever the initial condition is (see Remark~\ref{rmk:conv-bounded}). Note that \cite[Corollary~2.7]{S} also assumes $\sigma$ bounded in order to get uniqueness.
\item\label{it:compare-moment} Item \ref{it:extend-exist} of this remark holds without any moment assumption on $\xi$. In particular, equation \eqref{eq: mild} can be started from every deterministic $\phi\in C_{\theta}$, which is what the construction of a transition semigroup on $C_{\theta}$ requires (see Section~\ref{subsec: markov}).
\end{enumerate}
\end{remark}

\subsubsection{Preliminaries: a conjugation identity and subdifferentials}\label{sec:extend-tools}
The proof of Theorem~\ref{thm:extend} treats the drift term of \eqref{eq: mild} pathwise, through the deterministic mapping $\mathcal{M}$ of Section~\ref{sec:extend-M}. The uniqueness statement which we need for this mapping (Lemma~\ref{lem:Muniq} below) is obtained, as in \cite[Theorem~5.4]{S}, by differentiating in time the sup norm of a weighted difference of solutions. This rests on two elementary tools, which we gather here since they are used again in Section~\ref{sec: UniformBoundsU}: a conjugation identity between the Laplacian and the weight \eqref{eq: weight1}, and the subdifferential of the norm of $C_{0}$. The proof of the first one is a direct computation on the weight, and we defer it to Appendix~\ref{app:weights}.

\begin{lemma}\label{lem:conjugation}
Let $x_{0}\in\R^{d}$, $\theta>0$, and let $\hat{w}_{x_{0}}^{\theta}$ be the weight defined in \eqref{eq: weight1}. We define a differential operator $\mathcal{L}_{\theta}$ acting on functions $\psi\in C^{2}(\R^{d})$ by
\begin{equation}\label{eq:Ltheta-def}
\mathcal{L}_{\theta}\psi(x):=\frac{1}{2}\Delta\psi(x)+\frac{\theta\,(x-x_{0})\cdot\nabla\psi(x)}{1+|x-x_{0}|^{2}}+c_{\theta}(x)\,\psi(x),
\end{equation}
where the coefficient $c_{\theta}$ is given by
\begin{equation}\label{eq:Ltheta-coef}
c_{\theta}(x):=\frac{\theta\lp\frac{\theta+d}{2}-1\rp|x-x_{0}|^{2}+\frac{\theta}{2}\,d}{(1+|x-x_{0}|^{2})^{2}} .
\end{equation}
Then the following holds true.
\begin{enumerate}
\item\label{it:conj-identity} For every $\phi\in C^{2}(\R^{d})$ we have the conjugation identity
\begin{equation}\label{eq:conjugation}
\frac{1}{2} \hat{w}_{x_{0}}^{\theta}(x)\,\Delta\phi(x)=\mathcal{L}_{\theta}\big(\hat{w}_{x_{0}}^{\theta}\phi\big)(x),\qquad x\in\R^{d}.
\end{equation}
\item\label{it:conj-max} The coefficient $c_{\theta}$ is bounded on $\R^{d}$. Moreover, let $\psi\in C^{2}(\R^{d})\cap C_{0}(\R^{d})$ be nonzero and let $x^{*}\in\R^{d}$ be a point where $|\psi|$ attains its maximum, that is $|\psi(x^{*})|=|\psi|_{C_{0}}$. Then
\begin{equation}\label{eq:Ltheta-max}
\mathrm{sign}(\psi(x^{*}))\,\mathcal{L}_{\theta}\psi(x^{*})\leq |c_{\theta}|_{C_{b}}\,|\psi|_{C_{0}} .
\end{equation}
\end{enumerate}
\end{lemma}

We now turn to the second tool, namely the subdifferential of the sup norm, which is what allows one to differentiate the norm $|\psi(t,\cdot)|_{C_{0}}$ of a time dependent function $\psi$.

Recall that the space $C_{0}$ is defined by $(\ref{eq: vanishinfty})$.
We denote $C_{0}^{*}$ the dual space of $C_{0}=C_{0}(\R^d)$. It is a well-known fact (see e.g.~ \cite[Theorem C.18]{Con}) that $C_{0}^{*}$ is isomorphic to the space of signed regular finite Borel measures under the duality relation
\begin{equation}
    \left<\phi,\mu \right>_{C_{0},C_{0}^{*}}=\int_{\R^d}\phi(x)\mu(dx).
\end{equation}
The dual space $C_{0}^{*}$ is equipped with the total variation norm
\begin{equation}
    |\mu|_{C_{0}^{*}}:=\sup_{|\phi|_{C_{0}}\leq 1}\left<\phi,\mu\right>_{C_{0},C_{0}^{*}},
\end{equation}
where we recall that the norm in $C_{0}$ is a closed subspace of $C_{b}$ with norm $(\ref{eq:normBound})$. We are now ready to define the notion of subdifferential for a function in $C_{0}$.
\begin{definition}For any $\phi\in C_{0}(\R^d)$, the subdifferential is a subset of the dual space $C_{0}^{*}$ defined by
\begin{equation}\label{eq: subdifferential}
    \partial |\phi|_{C_{0}}:=\left\{\mu \in C_{0}^{*}: \left<\phi,\mu \right>_{C_{0}, C_{0}^{*}}=|\phi|_{C_{0}} \quad \text{and}\quad  |\mu|_{C_{0}^{*}}=1 \right\}.
\end{equation}
\end{definition}
The notion of subdifferential is developed e.g.~in \cite[Appendix D]{DZ}. The following two results, both taken from \cite[Propositions~3.3 and~3.5]{S}, describe the structure of subdifferential elements and their main application to bounding derivatives of sup norms.

\begin{proposition}\label{prop: subdiffEstimate}
Let $v \in C_{0}(\R^{d})$. Then the following holds true:
\begin{enumerate}
\item\label{item: subdiffChar} \emph{Characterization of subdifferential elements:}  
If $\mu \in \partial|v|_{C_{0}}$ as defined in~$(\ref{eq: subdifferential})$, then $\mu$ is a signed Borel measure of the form
\begin{equation}\label{eq: subdiffChar}
    \mu(dx) = \mathrm{sign}(v(x))\,\widetilde{\mu}(dx),
\end{equation}
where $\widetilde{\mu}$ is a positive Borel measure supported on the set $\{x \in \R^{d} : |v(x)| = |v|_{C_{0}}\}$ and satisfying $\widetilde{\mu}(\R^{d}) = 1$. In particular, for any $x_{0} \in \R^{d}$ with $|v(x_{0})| = |v|_{C_{0}}$, it is immediate from Definition~$(\ref{eq: subdifferential})$ that the measure $\delta_{x_{0}}\,\mathrm{sign}(v(x_{0}))$ is an element of  $\partial|v|_{C_{0}}$.

\item\label{item: subdiffEst} \emph{Upper-left-derivative estimate:} 
Let $v\colon [0,T]\times \R^{d} \to \R$ be a function that has continuous first derivatives and satisfies $v(t,\cdot)\in C_{0}$ for any $t\in [0,T]$. For any $t\in (0,T]$ and $\mu \in \partial|v(t,\cdot)|_{C_{0}}$ defined in~$(\ref{eq: subdifferential})$, the upper-left-derivative of the mapping $t\mapsto |v(t,\cdot)|_{C_{0}}$ satisfies
\begin{equation}\label{eq: derivEst}
    \frac{d^{-}}{dt}|v(t,\cdot)|_{C_{0}}\leq \left<\partial_{t}v(t,\cdot),\mu \right>_{C_{0},C_{0}^{*}},
\end{equation}
where the upper-left-derivative of a function $g\colon \R \to \R$ is defined as
\begin{equation}\label{eq: upperLeftDer}
    \frac{d^{-}}{dt}g(t):=\limsup_{h\downarrow 0}\frac{g(t)-g(t-h)}{h}.
\end{equation}
\end{enumerate}
\end{proposition}

\subsubsection{The mapping $\mathcal{M}$}\label{sec:extend-M}
We start by recalling the deterministic mapping introduced in \cite[Section~5]{S}. It will allow us to treat the drift term of \eqref{eq: mild} pathwise.

\begin{definition}\label{def:Mmap}
Let $T>0$, $z\in C_b([0,T]\times\mathbb{R}^d)$, and let $f\colon\mathbb{R}\to\mathbb{R}$ be the drift term of \eqref{eq: mild} satisfying Assumption~\ref{assump:reac}. We let $\mathcal{M}(z)$ be the unique function $v\in C_b([0,T]\times\mathbb{R}^d)$ solving the integral equation
\begin{equation}\label{eq:Mmap}
v(t,x)=\int_0^t\int_{\mathbb{R}^d}G(t-s,x-y)\,f\big(v(s,y)\big)\,dy\,ds+z(t,x),\qquad(t,x)\in[0,T]\times\mathbb{R}^d,
\end{equation}
where $G$ is the heat kernel of Definition~\ref{def: heatkernel}. The existence of such a $v$ is established in \cite[Theorem~5.2]{S}, and its uniqueness in $C_b([0,T]\times\mathbb{R}^d)$ in \cite[Theorem~5.4]{S}.
\end{definition}

The following proposition gathers the properties of $\mathcal{M}$ proved in \cite[Section~5]{S} which will be used below. Recall that the spaces $C_{x_0,\theta}([0,T]\times\mathbb{R}^d)$ were introduced in Definition~\ref{def: weightSpace2}, with norm \eqref{eq: wsnorm2}.

\begin{proposition}\label{prop:Mmap}
Let Assumption~\ref{assump:reac} hold, let $T>0$, let $\theta>0$, and let $\kappa$ be as in \eqref{assump: f1}. Then the following holds true.
\begin{enumerate}[label=\textnormal{(\roman*)}]
\item\label{it:M-ext} \emph{Extension.} The map $\mathcal{M}$ of Definition~\ref{def:Mmap} extends uniquely to a continuous map from $C_{0,\theta}([0,T]\times\mathbb{R}^d)$ into itself. For every $x_0\in\mathbb{R}^d$, this extension is also continuous from $C_{x_0,\theta}([0,T]\times\mathbb{R}^d)$ into itself. We still denote it by $\mathcal{M}$.
\item\label{it:M-eq} \emph{Integral equation.} For every $z\in C_{0,\theta}([0,T]\times\mathbb{R}^d)$, the function $v=\mathcal{M}(z)$ satisfies the integral equation \eqref{eq:Mmap}.
\item\label{it:M-lip} \emph{Lipschitz continuity uniformly in the center.} There exists a constant $C_{\theta,\kappa}>0$, depending only on $\theta$ and $\kappa$ and in particular independent of $x_0$ and $T$, such that for all $x_0\in\mathbb{R}^d$ and $z_1,z_2\in C_{{x_{0}},\theta}([0,T]\times\mathbb{R}^d)$:
\begin{equation}\label{eq:MLip}
\sup_{t\in[0,T]}\big|\mathcal{M}(z_1)(t,\cdot)-\mathcal{M}(z_2)(t,\cdot)\big|_{C_{x_0,\theta}}\le C_{\theta,\kappa}\,e^{C_{\theta,\kappa}T}\sup_{t\in[0,T]}\big|z_1(t,\cdot)-z_2(t,\cdot)\big|_{C_{x_0,\theta}} .
\end{equation}
\item\label{it:M-lin} \emph{Linear growth.} There exists a constant $C_{\theta,\kappa,f}>0$, depending only on $\theta$, $\kappa$ and $|f(0)|$, such that for all $x_0\in\mathbb{R}^d$ and $z\in C_{0,\theta}([0,T]\times\mathbb{R}^d)$:
\begin{equation}\label{eq:Mbound}
\sup_{t\in[0,T]}\big|\mathcal{M}(z)(t,\cdot)\big|_{C_{x_0,\theta}}\le C_{\theta,\kappa,f}\,e^{C_{\theta,\kappa,f}T}\Big(1+\sup_{t\in[0,T]}|z(t,\cdot)|_{C_{x_0,\theta}}\Big).
\end{equation}
\item\label{it:M-adapt} \emph{Adaptedness.} Let $z$ be a random field which is adapted with respect to an increasing right-continuous filtration and such that $\mathbf{P}\big(z\in C_{0,\theta}([0,T]\times\mathbb{R}^d)\big)=1$. Then $\mathcal{M}(z)$ is a well defined random field adapted to the same filtration.
\end{enumerate}
\end{proposition}

\begin{proof}
All the items are quoted from \cite[Section~5]{S}, except item \ref{it:M-eq}, which we deduce from items \ref{it:M-ext} and \ref{it:M-lin}. The growth condition on $f$ required there is the exponential bound \cite[(2.4)]{S}, which is implied by our polynomial bound \eqref{assump: growth} (see Remark~\ref{rmk: growth}).

Item \ref{it:M-ext} is \cite[Theorem~5.6]{S}, which gives the extension of $\mathcal{M}$ to $C_{x_0,\theta}([0,T]\times\mathbb{R}^d)$ for each center $x_0$. By Remark~\ref{rem: equiv-norms}, these spaces coincide as sets and have equivalent norms for all centers $x_0$. Moreover, $C_b([0,T]\times\mathbb{R}^d)$ is dense in each of them (see the proof of \cite[Theorem~5.6]{S}). Hence the extensions obtained for different centers coincide.

Items \ref{it:M-lip} and \ref{it:M-lin} are \cite[Theorem~5.4, relation (5.18)]{S} and \cite[Corollary~5.5, relation (5.29)]{S} for $z,z_1,z_2\in C_b([0,T]\times\mathbb{R}^d)$, and they are extended to $C_{x_0,\theta}([0,T]\times\mathbb{R}^d)$ in \cite[Theorem~5.6]{S}. The dependence of the constant on $|f(0)|$ in \eqref{eq:Mbound} comes from \cite[(5.30)]{S}.

Item \ref{it:M-eq} is the second assertion of \cite[Theorem~5.6]{S} when $\theta<2/\nu$, where $\nu$ is the exponent of the growth condition. This restriction ensures, under the exponential growth \eqref{eq: expGrowth} of \cite{S}, that $y\mapsto G(t-s,x-y)f(v(s,y))$ is integrable for $v\in C_{0,\theta}([0,T]\times\mathbb{R}^d)$. Under our polynomial growth condition \eqref{assump: growth}, the item holds for every $\theta>0$, by an easy elaboration of the proof of \cite[Theorem~5.6]{S} which we now sketch. One approximates $z$ by bounded forcing terms $z_{n}$, for which $\mathcal{M}(z_{n})$ solves \eqref{eq:Mmap} by Definition~\ref{def:Mmap}, and passes to the limit by dominated convergence. The only change with respect to \cite{S} is the dominating function, which is $C(1+|y|^{\theta\nu})G(t-s,x-y)$ instead of $C\exp(C|y|^{\theta\nu})G(t-s,x-y)$, owing to \eqref{assump: growth} and to the linear bound \eqref{eq:Mbound} on $\mathcal{M}(z_{n})$, which is available since items \ref{it:M-ext} and \ref{it:M-lin} have been established first. By item \ref{it:G-moment} of Definition~\ref{def: heatkernel}, this dominating function is integrable for every $\theta>0$. In particular, for every $(t,x)\in[0,T]\times\R^{d}$, the function $(s,y)\mapsto G(t-s,x-y)f(v(s,y))$ is integrable on $[0,t]\times\R^{d}$.

Item \ref{it:M-adapt} is \cite[Proposition~5.7]{S}, which is stated there for $\theta<2/\nu$ because its justification is the integral equation \eqref{eq:Mmap}. Once item \ref{it:M-eq} is available for every $\theta>0$, so is item \ref{it:M-adapt}.
\end{proof}

Proposition~\ref{prop:Mmap} does not settle the uniqueness of solutions to the integral equation \eqref{eq:Mmap} in the class we are working with. Indeed, item~\ref{it:M-ext} only says that the continuous extension of $\mathcal{M}$ is unique, and item~\ref{it:M-eq} only says that $\mathcal{M}(z)$ is a solution of \eqref{eq:Mmap} for $z\in C_{0,\theta}([0,T]\times\R^{d})$. Uniqueness of solutions is proved in \cite[Theorem~5.4]{S}, but in the class $C_{b}([0,T]\times\R^{d})$ and for a bounded forcing term $z$ only. We need it in the larger class $C_{0,\theta}([0,T]\times\R^{d})$, in which the extension of item~\ref{it:M-ext} takes its values. This is the object of Lemma~\ref{lem:Muniq} below, which is used twice in the proof of Theorem~\ref{thm:extend}. The argument, like the one of \cite[Theorem~5.4]{S}, bounds the sup norm of a weighted difference of solutions by differentiating it in time, with the two tools of Section~\ref{sec:extend-tools}: the conjugation identity of Lemma~\ref{lem:conjugation} gives the equation satisfied by the weighted difference, and the subdifferential estimate of Proposition~\ref{prop: subdiffEstimate} turns it into a differential inequality for its sup norm.
\begin{lemma}\label{lem:Muniq}
Let Assumption~\ref{assump:reac} hold, let $T>0$ and $\theta>0$. Consider the integral equation \eqref{eq:Mmap} of Definition~\ref{def:Mmap}, where now the forcing term $z$ is an element of $C_{0,\theta}([0,T]\times\R^{d})$ and the solution $v$ is only required to belong to $C_{0,\theta}([0,T]\times\R^{d})$, instead of $C_{b}([0,T]\times\R^{d})$ for both of them. Then any such solution $v$ satisfies
\begin{equation}\label{eq:Muniq}
v=\mathcal{M}(z),
\end{equation}
where $\mathcal{M}(z)$ is the extension given by item~\ref{it:M-ext} of Proposition~\ref{prop:Mmap}.
\end{lemma}

\begin{proof}
Throughout the proof we fix the center $x_{0}=0$ of the space $C_{0,\theta}([0,T]\times\R^{d})$, and we write $\hat{w}^{\theta}$ for the weight $\hat{w}_{0}^{\theta}$ of \eqref{eq: weight1}. We also set $u:=\mathcal{M}(z)$. By item~\ref{it:M-eq} of Proposition~\ref{prop:Mmap}, $u$ is an element of $C_{0,\theta}([0,T]\times\R^{d})$ which solves the integral equation \eqref{eq:Mmap}. We thus have two solutions $v$ and $u$ of \eqref{eq:Mmap} in $C_{0,\theta}([0,T]\times\R^{d})$, with the same forcing term $z$, and we have to prove that they coincide. We proceed in four steps. We first write the equation satisfied by the difference $v-u$, then the equation satisfied by its weighted version. We then bound the time derivative of the sup norm of the weighted difference, and we conclude by a comparison argument.

\smallskip

\noindent
{\it Step 1: The difference of the two solutions.} We set
\begin{equation}\label{eq:Muniq-diff}
w:=v-u,
\qquad\text{and}\qquad
g(s,y):=f(v(s,y))-f(u(s,y)),\qquad (s,y)\in[0,T]\times\R^{d}.
\end{equation}
Since $v$ and $u$ both satisfy \eqref{eq:Mmap} with the same $z$, the forcing term cancels in the difference and we get
\begin{equation}\label{eq:Muniq-potential}
w(t,x)=\int_{0}^{t}\int_{\R^{d}}G(t-s,x-y)\,g(s,y)\,dy\,ds,\qquad (t,x)\in[0,T]\times\R^{d}.
\end{equation}
Let us check that the integral in \eqref{eq:Muniq-potential} converges absolutely. Owing to the growth condition \eqref{assump: growth} and to the pointwise bound $|v(s,y)|\leq (1+|y|^{\theta})|v|_{C_{0,\theta}([0,T]\times\R^{d})}$ stemming from the definition \eqref{eq: wsnorm2} of the norm (and similarly for $u$), there exists a constant $C>0$ depending on $K,\nu,\theta$ and the norms of $v,u$ such that
\begin{equation}\label{eq:Muniq-g-growth}
|g(s,y)|\leq C\lp 1+|y|^{\theta\nu}\rp ,\qquad (s,y)\in[0,T]\times\R^{d}.
\end{equation}
The heat kernel $G(t-s,x-\cdot)$ integrates any polynomial (item \ref{it:G-moment} of Definition~\ref{def: heatkernel}), so that \eqref{eq:Muniq-potential} converges absolutely, for every value of $\theta>0$. Relation \eqref{eq:Muniq-potential} expresses $w$ as the heat potential of $g$. Therefore, as in \cite[(5.3)]{S}, $w$ is differentiable and solves
\begin{equation}\label{eq:Muniq-heat}
\partial_{t}w(t,x)=\frac{1}{2}\Delta w(t,x)+g(t,x),
\qquad w(0,x)=0 .
\end{equation}

\smallskip

\noindent
{\it Step 2: The weighted difference.} We define the weighted difference
\begin{equation}\label{eq:Muniq-tilde-w}
\widetilde{w}(t,x):=\hat{w}^{\theta}(x)\,w(t,x)=\frac{w(t,x)}{(1+|x|^{2})^{\theta/2}} .
\end{equation}
Since $v$ and $u$ belong to $C_{0,\theta}([0,T]\times\R^{d})$, the definition \eqref{eq: weightSpace2} of this space and the equivalence~\eqref{eq: universalC1} of the weights yield
\begin{equation}\label{eq:Muniq-vanish}
\lim_{|x|\to\infty}\sup_{t\in[0,T]}|\widetilde{w}(t,x)|=0 ,
\end{equation}
that is $\widetilde{w}\in C_{0}([0,T]\times\R^{d})$. In particular $t\mapsto\widetilde{w}(t,\cdot)$ is a continuous path with values in $C_{0}(\R^{d})$, and the function $t\mapsto|\widetilde{w}(t,\cdot)|_{C_{0}}$ is continuous on $[0,T]$ and vanishes at $t=0$. Multiplying \eqref{eq:Muniq-heat} by $\hat{w}^{\theta}$ and invoking the conjugation identity \eqref{eq:conjugation} of Lemma~\ref{lem:conjugation}, we obtain that $\widetilde{w}$ solves
\begin{equation}\label{eq:Muniq-weighted-eq}
\partial_{t}\widetilde{w}(t,x)=\mathcal{L}_{\theta}\widetilde{w}(t,x)+\hat{w}^{\theta}(x)\,g(t,x),
\end{equation}
where $\mathcal{L}_{\theta}$ is the operator defined in \eqref{eq:Ltheta-def}, with center $x_{0}=0$.

\smallskip

\noindent
{\it Step 3: Differential inequality for the sup norm.} Fix $t\in(0,T]$ such that $\widetilde{w}(t,\cdot)$ is nonzero (otherwise there is nothing to prove at time $t$). By \eqref{eq:Muniq-vanish} there exists a point $x_{t}\in\R^{d}$ such that $|\widetilde{w}(t,x_{t})|=|\widetilde{w}(t,\cdot)|_{C_{0}}$. According to item~\ref{item: subdiffChar} of Proposition~\ref{prop: subdiffEstimate} (quoted from \cite[Proposition~3.3]{S}), the measure $\mu_{t}:=\delta_{x_{t}}\,\mathrm{sign}(\widetilde{w}(t,x_{t}))$ belongs to the subdifferential $\partial|\widetilde{w}(t,\cdot)|_{C_{0}}$ of \eqref{eq: subdifferential}. Hence the upper-left-derivative estimate \eqref{eq: derivEst} of item~\ref{item: subdiffEst} of the same proposition, applied with this $\mu_{t}$ and the equation \eqref{eq:Muniq-weighted-eq}, reads
\begin{equation}\label{eq:Muniq-dini}
\frac{d^{-}}{dt}|\widetilde{w}(t,\cdot)|_{C_{0}}
\leq
\mathrm{sign}(\widetilde{w}(t,x_{t}))\,\mathcal{L}_{\theta}\widetilde{w}(t,x_{t})
+\mathrm{sign}(\widetilde{w}(t,x_{t}))\,\hat{w}^{\theta}(x_{t})\,g(t,x_{t}) .
\end{equation}
We bound the two terms in the right hand side of \eqref{eq:Muniq-dini} separately. The first one is handled by Lemma~\ref{lem:conjugation}, item~\ref{it:conj-max}, applied to $\psi=\widetilde{w}(t,\cdot)$ and $x^{*}=x_{t}$. Namely, setting $c_{\theta}^{*}:=|c_{\theta}|_{C_{b}}$, relation \eqref{eq:Ltheta-max} gives
\begin{equation}\label{eq:Muniq-term1}
\mathrm{sign}(\widetilde{w}(t,x_{t}))\,\mathcal{L}_{\theta}\widetilde{w}(t,x_{t})\leq c_{\theta}^{*}\,|\widetilde{w}(t,\cdot)|_{C_{0}} .
\end{equation}
As for the second term, the weight is positive, so that $\mathrm{sign}(\widetilde{w}(t,x_{t}))=\mathrm{sign}(w(t,x_{t}))=\mathrm{sign}(v(t,x_{t})-u(t,x_{t}))$. Thanks to the one-sided Lipschitz bound \eqref{eq: equivf1} on $f$, written at the points $u_{1}=v(t,x_{t})$ and $u_{2}=u(t,x_{t})$, we thus get
\begin{equation}\label{eq:Muniq-term2}
\mathrm{sign}(\widetilde{w}(t,x_{t}))\,\hat{w}^{\theta}(x_{t})\,g(t,x_{t})
\leq \kappa\,\hat{w}^{\theta}(x_{t})\,|v(t,x_{t})-u(t,x_{t})|
=\kappa\,|\widetilde{w}(t,x_{t})|
=\kappa\,|\widetilde{w}(t,\cdot)|_{C_{0}} .
\end{equation}
Here we have used the fact that \eqref{assump: f1} with the constant $\kappa$ implies the same relation with $\max(\kappa,0)$, so that we may and do assume $\kappa\geq 0$. Plugging \eqref{eq:Muniq-term1} and \eqref{eq:Muniq-term2} into \eqref{eq:Muniq-dini}, we end up with
\begin{equation}\label{eq:Muniq-gronwall}
\frac{d^{-}}{dt}|\widetilde{w}(t,\cdot)|_{C_{0}}\leq \lp c_{\theta}^{*}+\kappa\rp|\widetilde{w}(t,\cdot)|_{C_{0}} ,
\qquad t\in(0,T],
\end{equation}
which trivially also holds at the times $t$ where $\widetilde{w}(t,\cdot)$ vanishes.

\smallskip

\noindent
{\it Step 4: Conclusion.} The function $\varphi(t):=|\widetilde{w}(t,\cdot)|_{C_{0}}$ is continuous on $[0,T]$, satisfies $\varphi(0)=0$ and the linear differential inequality \eqref{eq:Muniq-gronwall}. Set $c:=c_{\theta}^{*}+\kappa$. By the comparison theorem for differential inequalities involving the upper-left derivative \cite[Theorem~1.3]{T}, applied to \eqref{eq:Muniq-gronwall} and to the linear equation $\dot{z}=cz$ with $z(0)=\varphi(0)=0$, we get $\varphi(t)\leq z(t)=0$ for every $t\in[0,T]$. Hence $\widetilde{w}\equiv 0$ on $[0,T]\times\R^{d}$, that is $v=u=\mathcal{M}(z)$, which is our claim \eqref{eq:Muniq}.
\end{proof}

\subsubsection{Two properties of the fixed-center space}\label{sec:extend-fixed}
We now record two elementary facts about the space $C_\theta$ of \eqref{eq:Ctheta}. The first one identifies the space-time space of Definition~\ref{def: weightSpace2} with a space of continuous paths. The second one controls on the space $C_{\theta}$, the convolution between the initial condition and the heat operator. In what follows we set, for $t\ge0$ and $\phi\in C_\theta$,
\begin{equation}\label{eq:heatsemigroup}
H_t\phi(x):=\int_{\mathbb{R}^d}G(t,x-y)\,\phi(y)\,dy,\qquad H_0\phi:=\phi,
\end{equation}
and we recall that $c_{\theta}=\max\{1,2^{\theta-1}\}$ is the constant of Lemma~\ref{lem:weight-elem} and $m_{\theta}$ the Gaussian moment of \eqref{eq:G-moment}.

Let us start with the identification of the space-time space $C_{0,\theta}([0,T]\times\R^{d})$ with a space of paths, whose proof is deferred to Appendix~\ref{app:weights}.
\begin{lemma}\label{lem:fc-paths}
Let $\theta>0$ and $T>0$, and let $C_{\theta}$ be the space defined in \eqref{eq:Ctheta}. Then, as sets, $C_{0,\theta}([0,T]\times\mathbb{R}^d)=C([0,T];C_\theta)$. Moreover, for every $\psi$ in this space,
\begin{equation*}
|\psi|_{C_{0,\theta}([0,T]\times\mathbb{R}^d)}=\sup_{t\in[0,T]}|\psi(t,\cdot)|_{C_\theta}.
\end{equation*}
\end{lemma}

Lemma~\ref{lem:fc-paths} describes the paths of the random fields we shall work with. We now turn to the effect of the heat semigroup $H_{t}$ on an element $\phi\in C_{\theta}$, which is the contribution of the initial condition to the mild equation \eqref{eq: mild}. The proof of the corresponding Lemma~\ref{lem:fc-heat} rests on the Gaussian moment bound \eqref{eq:G-moment} and on an approximation of $\phi$ by elements of $C_{0}$. It is deferred to Appendix~\ref{app:weights}.
\begin{lemma}\label{lem:fc-heat}
Let $\theta>0$ and $T>0$, let $C_{\theta}$ be the space defined in \eqref{eq:Ctheta}, and let $H_{t}$ be the heat semigroup \eqref{eq:heatsemigroup}. Then for every $\phi\in C_\theta$ the following holds true.
\begin{enumerate}
\item\label{it:fc-heat-bound} For all $t\ge0$, the weighted supremum of $H_{t}\phi$ defined by \eqref{eq:Ctheta} is finite. Namely, with the constant $c_{\theta}$ of Lemma~\ref{lem:weight-elem} and the Gaussian moment $m_{\theta}$ of \eqref{eq:G-moment},
\begin{equation}\label{eq:U0Bound}
|H_t\phi|_{C_\theta}\le c_\theta\big(1+m_\theta\,t^{\theta/2}\big)\,|\phi|_{C_\theta},\qquad t\ge0 .
\end{equation}
\item\label{it:fc-heat-cont} For all $t\ge0$ we have $H_{t}\phi\in C_{\theta}$, and the map $t\mapsto H_t\phi$ belongs to $C([0,T];C_\theta)$.
\end{enumerate}
\end{lemma}

Note that we will apply Lemma~\ref{lem:fc-heat} to the convolution between the heat kernel and the initial condition. Namely, for $\xi$ as in Assumption~\ref{assump: initial-cond-fixed} we set
\begin{equation}\label{eq:U0}
U_0(t,x):=H_t\xi(x)=\int_{\mathbb{R}^d}G(t,x-y)\,\xi(y)\,dy .
\end{equation}
The map $\phi\mapsto(t\mapsto H_t\phi)$ is linear, and by Lemma~\ref{lem:fc-heat} it is bounded from $C_{\theta}$ to $C([0,T];C_{\theta})$. Therefore $U_0$ is an $\mathcal{F}_0$-measurable random variable with values in $C([0,T];C_\theta)=C_{0,\theta}([0,T]\times\mathbb{R}^d)$ (Lemma~\ref{lem:fc-paths}), and by \eqref{eq:U0Bound}
\begin{equation}\label{eq:U0Bound-sup}
\sup_{t\in[0,T]}|U_0(t,\cdot)|_{C_\theta}\le c_\theta\big(1+m_\theta T^{\theta/2}\big)\,|\xi|_{C_\theta}\qquad\mathbf{P}\text{-a.s.}
\end{equation}

\subsubsection{The stochastic convolution}\label{sec:extend-conv}
The second ingredient is the following estimate on stochastic convolutions, established in \cite[Theorem~4.3]{S}. For an adapted random field $h$ we write
\begin{equation}\label{eq:StochConv}
\mathcal{Z}^h(t,x):=\int_0^t\int_{\mathbb{R}^d}G(t-s,x-y)\,h(s,y)\,W(ds,dy),\qquad(t,x)\in[0,T]\times\mathbb{R}^d,
\end{equation}
where the stochastic integral is understood in the Dalang--Walsh sense, as in Definition~\ref{def: mild}.

\begin{proposition}[{\cite[Theorem~4.3]{S}}]\label{prop:StochConvS}
Let Assumption~\ref{assump:noise} hold, let $T>0$, and let $\theta$ and $p$ satisfy \eqref{assump: theta-p}. There exists a constant $C_{T,p,\theta}>0$ such that for every adapted random field $h$ satisfying $\sup_{t\in[0,T]}\sup_{x\in\mathbb{R}^d}\mathbf{E}|h(t,x)|^p<\infty$, the stochastic convolution \eqref{eq:StochConv} satisfies
\begin{equation}\label{eq:StochConvBound}
\sup_{x_0\in\mathbb{R}^d}\mathbf{E}\Big[\sup_{t\in[0,T]}|\mathcal{Z}^h(t,\cdot)|^p_{C_{x_0,\theta}}\Big]\le C_{T,p,\theta}\sup_{t\in[0,T]}\sup_{x_0\in\mathbb{R}^d}\mathbf{E}|h(t,x_0)|^p ,
\end{equation}
and $\mathbf{P}\big(\mathcal{Z}^h\in C_{x_0,\theta}([0,T]\times\mathbb{R}^d)\big)=1$ for all $x_0\in\mathbb{R}^d$. Furthermore, for fixed $p$ and $\theta$,
\begin{equation}\label{eq:CTpTheta}
\lim_{T\to0}C_{T,p,\theta}=0 .
\end{equation}
\end{proposition}

Note that the condition $p>(1+\theta)(d+1)/\theta$ in \eqref{assump: theta-p} is equivalent to $\theta>(d+1)/(p-(d+1))$, which is the condition required in \cite[Theorem~4.3]{S}. Also, the constant $C_{T,p,\theta}$ depends on the spectral measure $\specm$ and on the exponent $\eta$ of \eqref{eq:strongDalang}. We do not indicate this dependence.

\begin{remark}\label{rmk:conv-bounded}
Two features of Proposition~\ref{prop:StochConvS} will be used repeatedly.
\begin{enumerate}[label=\textnormal{(\roman*)}]
\item\label{it:conv-bdd} \emph{Bounded integrands.} If $h=\sigma(v)$ for an adapted random field $v$, then Assumption~\ref{assump:diff} gives $\mathbf{E}|h(t,x)|^p\le M_\sigma^p$. Hence \eqref{eq:StochConvBound} yields
\begin{equation}\label{eq:StochConvSigma}
\sup_{x_0\in\mathbb{R}^d}\mathbf{E}\Big[\sup_{t\in[0,T]}|\mathcal{Z}^{\sigma(v)}(t,\cdot)|^p_{C_{x_0,\theta}}\Big]\le C_{T,p,\theta}\,M_\sigma^p ,
\end{equation}
\emph{whatever the field $v$ is}. In particular, no information on the growth of $v$, or of the initial condition, is needed. This is the place where the boundedness of $\sigma$ enters.
\item\label{it:conv-lin} \emph{Linearity.} For two such integrands, $\mathcal{Z}^{h_1}-\mathcal{Z}^{h_2}=\mathcal{Z}^{h_1-h_2}$ almost surely, by linearity of the stochastic integral.
\end{enumerate}
\end{remark}

\subsubsection{Proof of Theorem~\ref{thm:extend}}\label{sec:extend-proof}
Before giving the proof, let us explain how the mapping $\mathcal{M}$ enters. Assume for a moment that $u$ is a solution of \eqref{eq: mild} as in Theorem~\ref{thm:extend}, and set $Z:=\mathcal{Z}^{\sigma(u)}$ as in \eqref{eq:StochConv}. With the notation \eqref{eq:U0}, relation \eqref{eq: mild} reads
\begin{equation*}
u(t,x)=\int_{0}^{t}\int_{\R^{d}}G(t-s,x-y)\,f(u(s,y))\,dy\,ds+\big(U_{0}+Z\big)(t,x).
\end{equation*}
In other words, $u$ solves \eqref{eq:Mmap} with $z=U_{0}+Z$. If $u$ has paths in $C([0,T];C_{\theta})$, then by Lemma~\ref{lem:fc-paths} and Lemma~\ref{lem:Muniq} this means $u=\mathcal{M}(U_{0}+Z)$. Hence $u$ solves the fixed point equation
\begin{equation}\label{eq:fixedpoint}
u=\mathcal{M}\big(U_{0}+\mathcal{Z}^{\sigma(u)}\big).
\end{equation}
Conversely, let $u$ be an adapted random field with paths in $C([0,T];C_{\theta})$ which satisfies \eqref{eq:fixedpoint}. Then item~\ref{it:M-eq} of Proposition~\ref{prop:Mmap} asserts that $u$ solves \eqref{eq:Mmap} with $z=U_{0}+\mathcal{Z}^{\sigma(u)}$, and this relation is exactly \eqref{eq: mild} with $u_{0}=\xi$. Our proof thus consists in solving \eqref{eq:fixedpoint} by a contraction argument. Let us mention that our strategy departs slightly from~\cite[Section~6]{S}.
Indeed, the latter reference solves the equation $Z=\mathcal{Z}^{\sigma(\mathcal{M}(U_{0}+Z))}$ for the stochastic convolution $Z$, instead of~\eqref{eq:fixedpoint}. The two formulations are equivalent through $u=\mathcal{M}(U_{0}+Z)$, and we find it more natural to iterate on $u$ directly. Our proof is divided in five steps. Steps 1 to 3 construct the solution on a small interval $[0,T_{0}]$, Step 4 proves the moment bound and extends the solution to $[0,T]$, and Step 5 proves uniqueness.

\begin{proof}
Throughout the proof we fix $T>0$ and $\theta,p$ satisfying \eqref{assump: theta-p}. The filtration $(\mathcal{F}_{t})_{t\geq 0}$ of Assumption~\ref{assump:noise} satisfies the usual conditions, as required in Proposition~\ref{prop:Mmap}\ref{it:M-adapt}. The value of the constant $C$ may change from line to line, but it never depends on $\xi$.

\smallskip

\noindent
{\it Step 1: The space and the Picard map.} For $\tau\in(0,T]$, let $\mathcal{X}_{\tau}$ be the set of adapted random fields $Y$ on $[0,\tau]\times\R^{d}$ such that $Y\in C_{0,\theta}([0,\tau]\times\R^{d})$ almost surely and
\begin{equation}\label{eq:Xnorm}
\|Y\|_{\mathcal{X}_{\tau}}:=\bigg(\sup_{x_{0}\in\R^{d}}\mathbf{E}\Big[\sup_{t\in[0,\tau]}|Y(t,\cdot)|^{p}_{C_{x_{0},\theta}}\Big]\bigg)^{1/p}<\infty .
\end{equation}
This is the metric used in \cite[Section~6]{S}. Note that the supremum over the centers sits \emph{outside} the expectation. It is easily seen that the space $(\mathcal{X}_{\tau},\|\cdot\|_{\mathcal{X}_{\tau}})$, with fields identified up to indistinguishability, is a Banach space.

The solution $u$ cannot be sought in $\mathcal{X}_{\tau}$ itself. Indeed, Assumption~\ref{assump: initial-cond-fixed} imposes no moment on $\xi$, and even for a deterministic $\xi\in C_{\theta}$ the supremum over the centers of $|\xi|_{C_{x_{0},\theta}}$ may be infinite (see Remark~\ref{rmk:initial-fixed}). We therefore work in a translate of $\mathcal{X}_{\tau}$. To this aim, observe that by \eqref{eq:U0Bound-sup} and Lemma~\ref{lem:fc-paths}, $U_{0}$ of \eqref{eq:U0} is $\mathcal{F}_{0}$-measurable, hence adapted, with paths in $C_{0,\theta}([0,\tau]\times\R^{d})$. Items~\ref{it:M-ext} and~\ref{it:M-adapt} of Proposition~\ref{prop:Mmap} thus ensure that $\mathcal{M}(U_{0})$ is a well defined adapted random field with paths in $C_{0,\theta}([0,\tau]\times\R^{d})$. We set
\begin{equation}\label{eq:Yspace}
\mathcal{Y}_{\tau}:=\big\{u=\mathcal{M}(U_{0})+Y\,:\,Y\in\mathcal{X}_{\tau}\big\},
\qquad
d_{\tau}(u_{1},u_{2}):=\|u_{1}-u_{2}\|_{\mathcal{X}_{\tau}} .
\end{equation}
Every $u\in\mathcal{Y}_{\tau}$ is an adapted random field with paths in $C_{0,\theta}([0,\tau]\times\R^{d})$ almost surely, and $(\mathcal{Y}_{\tau},d_{\tau})$ is a complete metric space, being a translate of the Banach space $\mathcal{X}_{\tau}$.

For $u\in\mathcal{Y}_{\tau}$ we define
\begin{equation}\label{eq:Gamma}
\Gamma(u):=\mathcal{M}\big(U_{0}+\mathcal{Z}^{\sigma(u)}\big)\qquad\text{on }[0,\tau]\times\R^{d} .
\end{equation}
Let us check that $\Gamma$ maps $\mathcal{Y}_{\tau}$ into itself. Let $u\in\mathcal{Y}_{\tau}$. Since $u$ is adapted and continuous, $\sigma(u)$ is adapted and continuous, hence predictable, and it is bounded by $M_{\sigma}$. By Proposition~\ref{prop:StochConvS} and \eqref{eq:StochConvSigma}, the stochastic convolution $Z:=\mathcal{Z}^{\sigma(u)}$ has a version in $C_{0,\theta}([0,\tau]\times\R^{d})$, which belongs to $\mathcal{X}_{\tau}$ and satisfies
\begin{equation}\label{eq:GammaBound}
\|Z\|^{p}_{\mathcal{X}_{\tau}}\leq C_{\tau,p,\theta}\,M_{\sigma}^{p} .
\end{equation}
Hence $U_{0}+Z$ is adapted with paths in $C_{0,\theta}([0,\tau]\times\R^{d})$, and items~\ref{it:M-ext} and~\ref{it:M-adapt} of Proposition~\ref{prop:Mmap} ensure that $\Gamma(u)=\mathcal{M}(U_{0}+Z)$ is a well defined adapted random field with paths in $C_{0,\theta}([0,\tau]\times\R^{d})$. To this aim it remains to see that $\Gamma(u)-\mathcal{M}(U_{0})$ belongs to $\mathcal{X}_{\tau}$. We invoke the Lipschitz bound \eqref{eq:MLip} of the map $\mathcal{M}$ with $z_{1}:=U_{0}+Z$ and $z_{2}:=U_{0}$, so that $z_{1}-z_{2}=Z$. Since the constant in \eqref{eq:MLip} does not depend on the center $x_{0}$, we obtain, almost surely and for every $x_{0}\in\R^{d}$,
\begin{equation}\label{eq:Gamma-affine}
\sup_{s\in[0,\tau]}\big|\Gamma(u)(s,\cdot)-\mathcal{M}(U_{0})(s,\cdot)\big|_{C_{x_{0},\theta}}
\leq C_{\theta,\kappa}\,e^{C_{\theta,\kappa}\tau}\sup_{s\in[0,\tau]}|Z(s,\cdot)|_{C_{x_{0},\theta}} .
\end{equation}
Raising \eqref{eq:Gamma-affine} to the power $p$, taking expectations and then the supremum over $x_{0}$, and invoking \eqref{eq:GammaBound}, we end up with $\|\Gamma(u)-\mathcal{M}(U_{0})\|_{\mathcal{X}_{\tau}}<\infty$. Hence $\Gamma(\mathcal{Y}_{\tau})\subset\mathcal{Y}_{\tau}$.

\smallskip

\noindent
{\it Step 2: Contraction.} Let $u_{1},u_{2}\in\mathcal{Y}_{\tau}$, and set $Z_{i}:=\mathcal{Z}^{\sigma(u_{i})}$ for $i=1,2$, so that $\Gamma(u_{i})=\mathcal{M}(U_{0}+Z_{i})$. We invoke the Lipschitz bound \eqref{eq:MLip} with $z_{i}:=U_{0}+Z_{i}$. The crucial observation is that the heat flow $U_{0}$ cancels in the difference, that is $z_{1}-z_{2}=Z_{1}-Z_{2}$. As in \eqref{eq:Gamma-affine}, we obtain, almost surely and for every $x_{0}\in\R^{d}$,
\begin{equation}\label{eq:contr1}
\sup_{s\in[0,\tau]}\big|\Gamma(u_{1})(s,\cdot)-\Gamma(u_{2})(s,\cdot)\big|_{C_{x_{0},\theta}}\leq C_{\theta,\kappa}\,e^{C_{\theta,\kappa}\tau}\sup_{s\in[0,\tau]}|Z_{1}(s,\cdot)-Z_{2}(s,\cdot)|_{C_{x_{0},\theta}} .
\end{equation}
Next, we bound the right hand side of \eqref{eq:contr1}. By Remark~\ref{rmk:conv-bounded}\ref{it:conv-lin}, $Z_{1}-Z_{2}=\mathcal{Z}^{\sigma(u_{1})-\sigma(u_{2})}$. The integrand $\sigma(u_{1})-\sigma(u_{2})$ is adapted and bounded by $2M_{\sigma}$, so that it satisfies the hypothesis of Proposition~\ref{prop:StochConvS}. Applying \eqref{eq:StochConvBound}, and then the Lipschitz property of $\sigma$ (Assumption~\ref{assump:diff}), we get
\begin{equation}\label{eq:contr2}
\|Z_{1}-Z_{2}\|^{p}_{\mathcal{X}_{\tau}}\leq C_{\tau,p,\theta}\,L_{\sigma}^{p}\sup_{s\in[0,\tau]}\sup_{x_{0}\in\R^{d}}\mathbf{E}\big|u_{1}(s,x_{0})-u_{2}(s,x_{0})\big|^{p} .
\end{equation}
Finally, we bound the right hand side of \eqref{eq:contr2} by the distance \eqref{eq:Yspace}. For each fixed $x_{0}$, the pointwise bound \eqref{eq: pointwiseNorm} of Remark~\ref{rmk: pointwiseNorm}, applied \emph{with center $x_{0}$}, gives $|u_{1}(s,x_{0})-u_{2}(s,x_{0})|\leq|u_{1}(s,\cdot)-u_{2}(s,\cdot)|_{C_{x_{0},\theta}}$. Hence
\begin{equation}\label{eq:contr3}
\sup_{s\in[0,\tau]}\sup_{x_{0}\in\R^{d}}\mathbf{E}\big|u_{1}(s,x_{0})-u_{2}(s,x_{0})\big|^{p}\leq\sup_{x_{0}\in\R^{d}}\mathbf{E}\Big[\sup_{s\in[0,\tau]}|u_{1}(s,\cdot)-u_{2}(s,\cdot)|^{p}_{C_{x_{0},\theta}}\Big]=d_{\tau}(u_{1},u_{2})^{p}.
\end{equation}
Raising \eqref{eq:contr1} to the power $p$, taking expectations and the supremum over $x_{0}$, and then plugging \eqref{eq:contr3} into \eqref{eq:contr2}, we have thus shown that
\begin{equation}\label{eq:contraction}
d_{\tau}\big(\Gamma(u_{1}),\Gamma(u_{2})\big)^{p}\leq q_{\tau}\,d_{\tau}(u_{1},u_{2})^{p},\qquad\text{where}\quad q_{\tau}:=C_{\tau,p,\theta}\,L_{\sigma}^{p}\,C_{\theta,\kappa}^{p}\,e^{pC_{\theta,\kappa}\tau}.
\end{equation}
Owing to \eqref{eq:CTpTheta}, we can choose $T_{0}\in(0,T]$ such that
\begin{equation}\label{eq:T0choice}
q_{T_{0}}<1 .
\end{equation}
We stress that $T_{0}$ depends only on $p,\theta,\kappa,L_{\sigma}$ and on the constants of Proposition~\ref{prop:StochConvS}. In particular it depends neither on the initial condition $\xi$ nor on its law.

\smallskip

\noindent
{\it Step 3: A solution on $[0,T_{0}]$.} By \eqref{eq:contraction}, \eqref{eq:T0choice} and the completeness of $\mathcal{Y}_{T_{0}}$ established in Step~1, Banach's fixed point theorem provides a unique $u\in\mathcal{Y}_{T_{0}}$ such that $\Gamma(u)=u$. That is, $u$ solves \eqref{eq:fixedpoint} on $[0,T_{0}]$. As an element of $\mathcal{Y}_{T_{0}}$, $u$ is adapted with paths in $C_{0,\theta}([0,T_{0}]\times\R^{d})=C([0,T_{0}];C_{\theta})$, where we have used Lemma~\ref{lem:fc-paths}. Hence, as explained before the proof and summarized in~\eqref{eq:fixedpoint}, $u$ is a mild solution of \eqref{eq: mild} on $[0,T_{0}]$. In addition, setting $Z:=\mathcal{Z}^{\sigma(u)}$, so that $u=\mathcal{M}(U_{0}+Z)$, the bound \eqref{eq:GammaBound} reads
\begin{equation}\label{eq:Zbound}
\sup_{x_{0}\in\R^{d}}\mathbf{E}\Big[\sup_{t\in[0,T_{0}]}|Z(t,\cdot)|^{p}_{C_{x_{0},\theta}}\Big]\leq C_{T_{0},p,\theta}\,M_{\sigma}^{p}<\infty .
\end{equation}

\smallskip

\noindent
{\it Step 4: Moment bound and extension to $[0,T]$.} We first assume that \eqref{eq:initial_bound} holds. We apply the linear growth bound \eqref{eq:Mbound} with center $x_{0}=0$ and $z=U_{0}+Z$. Combining it with \eqref{eq:U0Bound-sup}, we get, almost surely,
\begin{equation*}
\sup_{t\in[0,T_{0}]}|u(t,\cdot)|_{C_{\theta}}\leq C_{\theta,\kappa,f}\,e^{C_{\theta,\kappa,f}T_{0}}\Big(1+c_{\theta}\big(1+m_{\theta}T_{0}^{\theta/2}\big)|\xi|_{C_{\theta}}+\sup_{t\in[0,T_{0}]}|Z(t,\cdot)|_{C_{\theta}}\Big).
\end{equation*}
Raising this inequality to the power $p$, taking expectations, and invoking \eqref{eq:initial_bound} and \eqref{eq:Zbound} (with $x_{0}=0$), we end up with
\begin{equation}\label{eq:solution_bound_T0}
\mathbf{E}\Big[\sup_{t\in[0,T_{0}]}|u(t,\cdot)|^{p}_{C_{\theta}}\Big]\leq C\Big(1+\big(1+T_{0}^{\theta/2}\big)^{p}\,\mathbf{E}\big[|\xi|^{p}_{C_{\theta}}\big]+C_{T_{0},p,\theta}M_{\sigma}^{p}\Big)<\infty ,
\end{equation}
for a constant $C=C_{p,\theta,\kappa,f,\theta,T_{0}}$. This is \eqref{eq:solution_bound} on the interval $[0,T_{0}]$.

Since $T_{0}$ is independent of the initial condition, a concatenation argument extends the local solution to $[0,T]$. Namely, we successively restart the equation at the endpoints of intervals of length at most $T_{0}$, taking the terminal value of the preceding solution as the new initial condition. This terminal value is measurable with respect to the current $\sigma$-field and takes its values in $C_{\theta}$, and Lemma~\ref{lem:restart} below ensures that the concatenated process satisfies \eqref{eq: mild} on $[0,T]$. Furthermore, the terminal value satisfies \eqref{eq:initial_bound} by \eqref{eq:solution_bound_T0}, and the iteration of \eqref{eq:solution_bound_T0} over the finitely many intervals yields \eqref{eq:solution_bound}.

Let us stress that the moment condition \eqref{eq:initial_bound} has only been used to derive \eqref{eq:solution_bound_T0}, hence \eqref{eq:solution_bound}. Steps 1 to 3 and the concatenation argument hold for every $\xi$ satisfying Assumption~\ref{assump: initial-cond-fixed}. This proves the existence part of item \ref{it:extend-exist}, as well as item \ref{it:extend-moment}.

\smallskip

\noindent
{\it Step 5: Uniqueness.} Let $\tilde{u}$ be another adapted random field with paths in $C([0,T];C_{\theta})$ which solves \eqref{eq: mild} with $u_{0}=\xi$, and set $\tilde{Z}:=\mathcal{Z}^{\sigma(\tilde{u})}$. Since $\sigma(\tilde{u})$ is adapted, continuous and bounded by $M_{\sigma}$, Proposition~\ref{prop:StochConvS} and \eqref{eq:StochConvSigma} show that $\tilde{Z}$ has a version in $C_{0,\theta}([0,T]\times\R^{d})$ which belongs to $\mathcal{X}_{\tau}$ for every $\tau\in(0,T]$. Next, relation \eqref{eq: mild} for $\tilde{u}$ states that $\tilde{u}$ solves \eqref{eq:Mmap} with $z=U_{0}+\tilde{Z}$, and both $\tilde{u}$ and $U_{0}+\tilde{Z}$ belong to $C_{0,\theta}([0,T]\times\R^{d})$ almost surely by Lemma~\ref{lem:fc-paths}. Therefore Lemma~\ref{lem:Muniq} yields
\begin{equation*}
\tilde{u}=\mathcal{M}(U_{0}+\tilde{Z})=\mathcal{M}\big(U_{0}+\mathcal{Z}^{\sigma(\tilde{u})}\big)=\Gamma(\tilde{u}).
\end{equation*}
Moreover, the same argument as in \eqref{eq:Gamma-affine}, with $Z$ replaced by $\tilde{Z}$, gives $\|\tilde{u}-\mathcal{M}(U_{0})\|_{\mathcal{X}_{\tau}}<\infty$, so that $\tilde{u}\in\mathcal{Y}_{\tau}$ for every $\tau\in(0,T]$. That is, $\tilde{u}$ is a fixed point of $\Gamma$ in $\mathcal{Y}_{T_{0}}$. By the uniqueness part of Banach's fixed point theorem, $\tilde{u}=u$ on $[0,T_{0}]$. In particular $\tilde{u}(T_{0},\cdot)=u(T_{0},\cdot)$. We can then extend this identity to the whole interval $[0,T]$ by the same patching technique as above.
This completes the proof of item \ref{it:extend-exist}.
\end{proof}

\subsection{Markov property and transition semigroup}\label{subsec: markov}
We begin with the restart property of equation \eqref{eq: mild}, which is behind the Markov property of the solution, and which we shall also use in Section~\ref{sec:feller}. For $s\geq 0$, the shifted noise $W^{s}$ is defined by
\begin{equation}\label{eq:shifted-noise}
\int_{0}^{t}\int_{\R^{d}}h(r,y)\,W^{s}(dr,dy):=\int_{s}^{s+t}\int_{\R^{d}}h(r-s,y)\,W(dr,dy),\qquad t\geq 0,
\end{equation}
for every random field $h$ which is adapted to the shifted filtration $(\mathcal{F}_{s+t})_{t\geq 0}$. Since $\dot{W}$ is white in time, $W^{s}$ is a noise satisfying Assumption~\ref{assump:noise} with respect to $(\mathcal{F}_{s+t})_{t\geq 0}$, with the same covariance as $W$, and it is independent of $\mathcal{F}_{s}$.

\begin{lemma}\label{lem:restart}
Let Assumptions~\ref{assump:noise}, \ref{assump:reac} and \ref{assump:diff} hold, let $\xi$ satisfy Assumption~\ref{assump: initial-cond-fixed}, and let $T>0$ and $s\in[0,T)$. Let $u$ be an adapted random field with paths in $C([0,T];C_{\theta})$, and set $v(t,x):=u(s+t,x)$ for $(t,x)\in[0,T-s]\times\R^{d}$. Then $u$ is a mild solution of \eqref{eq: mild} on $[0,T]$ with initial condition $\xi$ if and only if the following two conditions hold.
\begin{enumerate}
\item\label{it:restart-before} The restriction of $u$ to $[0,s]\times\R^{d}$ is a mild solution of \eqref{eq: mild} on $[0,s]$ with initial condition $\xi$.
\item\label{it:restart-after} The field $v$ is a mild solution of \eqref{eq: mild} on $[0,T-s]$ with initial condition $u(s,\cdot)$, driven by the shifted noise $W^{s}$ of \eqref{eq:shifted-noise} with the filtration $(\mathcal{F}_{s+t})_{t\geq 0}$.
\end{enumerate}
In particular, let $u=u(\cdot\,;\xi)$ be the solution given by Theorem~\ref{thm:extend}. Then $u(s,\cdot)$ is an $\mathcal{F}_{s}$-measurable random variable with values in $C_{\theta}$, independent of $W^{s}$, and the field $v$ is the solution given by Theorem~\ref{thm:extend} on $[0,T-s]$ for the noise $W^{s}$ and the initial condition $u(s,\cdot)$.
\end{lemma}

\begin{proof}
Only the equivalence between the two items and the mild equation on $[0,T]$ is used in the proof of Theorem~\ref{thm:extend} (Step 4 there); the last assertion of the lemma is a consequence of that theorem, which is why it is proved after it. The equivalence is a standard fact from stochastic calculus. It is obtained by splitting the time integrals of \eqref{eq: mild} at the time $s$, and by invoking the semigroup property \eqref{eq:G-semigroup} of the heat kernel $G$ together with the stochastic Fubini theorem. As for the last assertion, $u(s,\cdot)$ is $\mathcal{F}_{s}$-measurable with values in $C_{\theta}$ since $u$ is adapted with paths in $C([0,T];C_{\theta})$, and it is independent of $W^{s}$ because $W^{s}$ is independent of $\mathcal{F}_{s}$. Hence $u(s,\cdot)$ satisfies Assumption~\ref{assump: initial-cond-fixed} with respect to the filtration $(\mathcal{F}_{s+t})_{t\geq 0}$ and the noise $W^{s}$. By item \ref{it:restart-after} and the uniqueness in Theorem~\ref{thm:extend}, $v$ is the solution given by that theorem for these data.
\end{proof}

Item \ref{it:extend-exist} of Theorem~\ref{thm:extend} allows us to start equation \eqref{eq: mild} from every deterministic $\phi\in C_{\theta}$. We denote by $u(t,\cdot\,;\phi)$ the corresponding solution. This is the basis for the definition of the transition semigroup of \eqref{NSHE} on the state space $C_{\theta}$.
\begin{definition}\label{def: transition}
Under the assumptions of Theorem~\ref{thm:extend}, let $B_{b}(C_{\theta})$ be the space of bounded Borel functions $F\colon C_{\theta}\to\R$. For $t\geq 0$, $F\in B_{b}(C_{\theta})$ and $\phi\in C_{\theta}$, we set
\begin{equation}\label{eq: transition}
\mathcal{P}_{t}F(\phi):=\mathbf{E}\big[F(u(t,\cdot\,;\phi))\big].
\end{equation}
For a probability measure $\mu$ on $C_{\theta}$, we also define the measure $\mathcal{P}_{t}^{*}\mu$ by the duality relation
\begin{equation}\label{eq: Pt}
\int_{C_{\theta}}\mathcal{P}_{t}F(\phi)\,\mu(d\phi)=\int_{C_{\theta}}F(\phi)\,\mathcal{P}_{t}^{*}\mu(d\phi)=\mathbf{E}\big[F(u(t,\cdot\,;\xi))\big],\qquad F\in B_{b}(C_{\theta}),
\end{equation}
where $\xi$ is any initial condition satisfying Assumption~\ref{assump: initial-cond-fixed} with $\mathrm{Law}(\xi)=\mu$. In other words, $\mathcal{P}_{t}^{*}\mu=\mathrm{Law}(u(t,\cdot\,;\xi))$. The measurability of $\phi\mapsto\mathcal{P}_{t}F(\phi)$, which makes the first integral in \eqref{eq: Pt} meaningful, and the last identity in \eqref{eq: Pt}, which identifies the measure $\mathcal{P}_{t}^{*}\mu$ with the law of the solution started from a random initial condition, are established in Section~\ref{subsec:feller-two} (Proposition~\ref{prop:feller-two} and Lemma~\ref{lem:plug-in}).
\end{definition}
\begin{remark}\label{rmk:markov-property}
Let us comment on Definition~\ref{def: transition}.
\begin{enumerate}
\item The solution $u$ to equation \eqref{NSHE} is a Markov process under the assumptions of Theorem~\ref{thm:extend}. Specifically, if $\varphi\in C_{b}(C_{\theta})$, then for all $t,h\geq 0$ we have
\begin{equation}\label{eq: markov}
\mathbf{E}\big[\varphi(u(t+h,\cdot\,;\xi))\,|\,\mathcal{F}_{t}\big]=\mathcal{P}_{h}\varphi(u(t,\cdot\,;\xi)) .
\end{equation}
This property is standard and can be proved as in \cite{DPZ96}, thanks to the uniqueness in item \ref{it:extend-exist} of Theorem~\ref{thm:extend} and to the restart property of Lemma~\ref{lem:restart}. Relation \eqref{eq: markov} entails the semigroup property $\mathcal{P}_{t+h}=\mathcal{P}_{t}\mathcal{P}_{h}$, and by duality $\mathcal{P}_{t+h}^{*}=\mathcal{P}_{t}^{*}\mathcal{P}_{h}^{*}$.
\item Relation \eqref{eq: transition} requires the map $\phi\mapsto u(t,\cdot\,;\phi)$ to be measurable from $C_{\theta}$ to $C_{\theta}$, and the Feller property $\mathcal{P}_{t}(C_{b}(C_{\theta}))\subset C_{b}(C_{\theta})$ requires its continuity in probability. The contraction argument of Section~\ref{sec:extend-proof} gives neither of them. Indeed, it closes in the norm \eqref{eq:Xnorm}, whose supremum over the centers $x_{0}$ turns the initial condition term into the supremum norm $|\phi_{1}-\phi_{2}|_{C_{b}}$. These questions are taken up in Section~\ref{sec:feller}. There we shall prove a version of the estimate on the stochastic convolution whose input is weighted (see Proposition~\ref{prop:weighted-conv}), and we shall derive from it the continuity of $\phi\mapsto u(t,\cdot\,;\phi)$ from a space $C_{\theta_{0}}$ to the space $C_{\theta}$, for exponents $\theta_{0}<\theta$. This yields a Feller property between the two exponents $\theta_{0}$ and $\theta$, which is what the existence of an invariant measure requires in Section~\ref{sec:invariance}.
\item In the framework of \cite{S}, the semigroup $\mathcal{P}_{t}$ could not be defined pointwise on the state space, since the initial condition had to satisfy the uniform moment condition of Remark~\ref{rmk:initial-fixed}. Only the dual operator $\mathcal{P}_{t}^{*}$ on the probability measures with a $p$-th moment uniform over the centers was available, together with the duality relation \eqref{eq: Pt} for $\mu$-almost every $\phi$. This is the main gain of Theorem~\ref{thm:extend}.
\end{enumerate}
\end{remark}

\section{Continuity with respect to the initial condition and Feller property}\label{sec:feller}

In this section we start the study of the dependence of the solution $u(t,\cdot\,;\phi)$ of Theorem~\ref{thm:extend} on its initial condition $\phi$, with the Feller property of the transition semigroup $\mathcal{P}_{t}$ of Definition~\ref{def: transition} in view. As explained in item (ii) of Remark~\ref{rmk:markov-property}, the Feller property requires the continuity of $\phi\mapsto u(t,\cdot\,;\phi)$ for the norm of the state space $C_{\theta}$. Unfortunately, the contraction argument of Section~\ref{sec:extend-proof} does not give this continuity. Indeed it is based on Proposition~\ref{prop:StochConvS} for the stochastic convolution, which takes a supremum of unweighted moments as input. The first step is therefore a version of Proposition~\ref{prop:StochConvS} whose input is a supremum of \emph{weighted} moments, centered at the same point $x_{0}$ as the weight in the output. This is the object of Sections~\ref{subsec:weighted-conv} to \ref{subsec:weighted-conv-proof}. We shall see that the exponent of the weight is still lost in the Kolmogorov step, and by the same amount as in \cite[Theorem~4.3]{S}. The second step, carried out in Sections~\ref{subsec:X-weighted} and \ref{subsec:continuity}, absorbs this loss in a norm which penalizes the centers far from the origin. It yields the continuity of $\phi\mapsto u(t,\cdot\,;\phi)$ from a space $C_{\theta_{0}}$ to the space $C_{\theta}$, for exponents $\theta_{0}<\theta$. Finally, Section~\ref{subsec:feller-two} draws the consequences for the transition semigroup $\mathcal{P}_{t}$, namely its measurability on $C_{\theta_{0}}$ and a Feller property between the two exponents.
\subsection{A weighted estimate for the stochastic convolution}\label{subsec:weighted-conv}

As mentioned above, Proposition~\ref{prop:StochConvS} has respective input and output defined by
\begin{equation}\label{eq:StochConvS-io}
\sup_{t\in[0,T]}\sup_{x\in\R^{d}}\mathbf{E}\big[|h(t,x)|^{p}\big]
\qquad\text{and}\qquad
\mathbf{E}\Big[\sup_{t\in[0,T]}\big|\mathcal{Z}^{h}(t,\cdot)\big|^{p}_{C_{x_{0},\theta}}\Big].
\end{equation}
The exponent $\theta$ has thus to be interpreted as a loss, paid in the passage from pointwise moments to a supremum over space. In this section we show that when the input moments are themselves weighted, with a weight $w_{x_{0}}^{\theta_{0}}$, the output is a supremum weighted by $w_{x_{0}}^{\theta_{0}+\theta_{1}}$, where $\theta_{1}$ is subject to the same lower bound as $\theta$ in Proposition~\ref{prop:StochConvS}. Let us insist on the fact that the center $x_{0}$ is the same for the input and the output. This is what makes the estimate useful for the fixed-center space $C_{\theta}$ of \eqref{eq:Ctheta}.

Let us define the weighted moment functional we shall work with. For $x_{0}\in\R^{d}$, $\theta_{0}\geq 0$, $T>0$ and an adapted random field $h=\{h(s,y)\,;\,(s,y)\in[0,T]\times\R^{d}\}$, we set
\begin{equation}\label{eq:Nfunctional}
\mathcal{N}^{x_{0}}_{p,\theta_{0},T}(h):=\sup_{s\in[0,T]}\sup_{y\in\R^{d}}w_{x_{0}}^{\theta_{0}}(y)\,\|h(s,y)\|_{L^{p}(\Omega)},
\end{equation}
where $w_{x_{0}}^{\theta_{0}}$ is the weight of \eqref{eq: weight1}. For $\theta_{0}=0$, the $p$-th power of \eqref{eq:Nfunctional} is the input of Proposition~\ref{prop:StochConvS} in \eqref{eq:StochConvS-io}. Observe that the definition of \eqref{eq:Nfunctional} does not require any continuity of $h$. The main result of this section is the following.

\begin{proposition}\label{prop:weighted-conv}
Let Assumption~\ref{assump:noise} hold, let $T>0$, let $\theta_{0}\geq 0$ and $\theta_{1}>0$, and set $\theta:=\theta_{0}+\theta_{1}$. Let $p\geq p^{*}$ in the sense of Remark~\ref{rmk:pstar}. Then there exists a constant $C_{T}=C_{T,p,\theta_{0},\theta_{1}}$, which does not depend on $x_{0}$ and satisfies
\begin{equation}\label{eq:CT-to-zero}
\lim_{T\to 0}C_{T}=0,
\end{equation}
such that the following holds. For every $x_{0}\in\R^{d}$ and every adapted random field $h$ on $[0,T]\times\R^{d}$ with $\mathcal{N}^{x_{0}}_{p,\theta_{0},T}(h)<\infty$, the stochastic convolution $\mathcal{Z}^{h}$ of \eqref{eq:StochConv} admits a modification with paths in $C_{x_{0},\theta}([0,T]\times\R^{d})$, and this modification satisfies
\begin{equation}\label{eq:weighted-conv}
\mathbf{E}\Big[\sup_{t\in[0,T]}\big|\mathcal{Z}^{h}(t,\cdot)\big|^{p}_{C_{x_{0},\theta}}\Big]\leq C_{T}\,\mathcal{N}^{x_{0}}_{p,\theta_{0},T}(h)^{p}.
\end{equation}
\end{proposition}

For $\theta_{0}=0$, Proposition~\ref{prop:weighted-conv} is Proposition~\ref{prop:StochConvS}, with $\theta=\theta_{1}$. For $\theta_{0}>0$ it says that the loss of exponent is still $\theta_{1}$, whatever the exponent $\theta_{0}$ of the input is. Let us also mention that in our setting the integrands are of the form $\sigma(v)$ for a random field $v$, and are thus bounded by $M_{\sigma}$ thanks to Assumption~\ref{assump:diff}. Hence all moments are available, which is what allows us to take $p\geq p^{*}$ as in Remark~\ref{rmk:pstar}, and the loss $\theta_{1}$ can be made as small as we wish by increasing $p$. We shall see in the sequel that a vanishing loss is not needed for the Feller property. The proof is an elaboration of the two steps of \cite[Theorem~4.3]{S}. The first one consists in moment bounds on the pointwise values and on the increments of $\mathcal{Z}^{h}$. The second one is the Kolmogorov criterion of \cite[Corollary~4.2]{S}, applied to the weighted field $\hat{w}_{x_{0}}^{\theta_{0}}\mathcal{Z}^{h}$. The proof is given in Section~\ref{subsec:weighted-conv-proof}. We first gather, in Section~\ref{subsec:conv-tools}, the four tools it needs.

\subsection{Four tools}\label{subsec:conv-tools}

The first tool is the spectral bound which is behind all the moment estimates of this paper. It concerns the Fourier transform of powers of the heat kernel, integrated against the spectral measure $\specm$ of Assumption~\ref{assump:noise}. Its proof is a computation on the explicit expression \eqref{eq: heatkernel} of the heat kernel, and we defer it to Appendix~\ref{app:conv}.

\begin{lemma}\label{lem: spectral}
Let $\specm$ be the spectral measure satisfying the strong Dalang condition~\eqref{eq:strongDalang} with exponent $\eta\in(0,1)$. For any $\hat{\alpha}\in(0,1]$ and $\tau>0$, the Fourier transform of $G(\tau,\cdot)^{\hat{\alpha}}$ is
\begin{equation}\label{eq: FourierPow}
    \mathcal{F}\bigl(G(\tau,\cdot)^{\hat\alpha}\bigr)(\xi) = C_{\hat{\alpha},d}\,\tau^{d(1-\hat\alpha)/2}\exp\!\left(-\frac{|\xi|^{2}\tau}{2\hat\alpha}\right),
    \quad\text{with}\quad C_{\hat\alpha,d}\in(0,\infty).
\end{equation}
Moreover, for any $\vk\in(1-\eta,1)$ there exists a constant $C_{\hat\alpha,d,\specm,\eta,\vk}<\infty$ such that for all $\tau>0$,
\begin{equation}\label{eq: spectralBound}
    \int_{\R^{d}}\bigl|\mathcal{F}\bigl(G(\tau,\cdot)^{\hat\alpha}\bigr)(\xi)\bigr|^{2}\specm(d\xi)\leq C_{\hat\alpha,d,\specm,\eta,\vk}\,\tau^{d(1-\hat\alpha)}\bigl(1+\tau^{-\vk}\bigr).
\end{equation}
\end{lemma}

In the sequel we shall mostly use this lemma with $\hat{\alpha}=1$. In that case $C_{1,d}=1$ in \eqref{eq: FourierPow}, so that \eqref{eq: spectralBound} reads
\begin{equation}\label{eq:spectral-one}
\int_{\R^{d}}e^{-\tau|\xi|^{2}}\,\specm(d\xi)\leq C_{\specm,\eta,\vk}\bigl(1+\tau^{-\vk}\bigr),\qquad\tau>0,\quad\vk\in(1-\eta,1).
\end{equation}

The second tool is the Kolmogorov continuity criterion on $\R^{d}$ of \cite{S}, which produces a weighted supremum out of pointwise moment bounds. We quote it as follows.

\begin{proposition}[{\cite[Corollary~4.2]{S}}]\label{prop: Salins}
Let $X\colon[0,T]\times\R^{d}\times\Omega\to\R$ be a jointly measurable random field, let $x_{0}\in\R^{d}$, and assume that there exist $A_{T}>0$, $p>d+1$ and $\lambda,\gamma\in(0,1)$ such that for all $x,y\in\R^{d}$ and all $s,t\in[0,T]$,
\begin{align}
\label{eq: propSalins1}     \mathbf{E}\big[|X(0,x_{0})|^{p}\big]&\leq A_{T},\\
\label{eq: propSalins2}  \mathbf{E}\big[|X(t,x)-X(t,y)|^{p}\big]&\leq A_{T}\,|x-y|^{\lambda p},\\
\label{eq: propSalins3}  \mathbf{E}\big[|X(t,x)-X(s,x)|^{p}\big]&\leq A_{T}\,|t-s|^{\gamma p}.
\end{align}
Assume moreover that
\begin{equation}\label{eq: propSalins4}
p\lambda>d+1,\qquad p\gamma>d+1,\qquad\text{and}\qquad\theta>\frac{p\lambda}{p-(d+1)}.
\end{equation}
Then $X$ admits a modification with paths in $C_{x_{0},\theta}([0,T]\times\R^{d})$, and there exists a constant $C_{p,d,\lambda,\gamma,\theta,T}>0$, independent of $x_{0}$ and $A_{T}$ and nondecreasing in $T$, such that this modification satisfies
\begin{equation}\label{eq: propSalins5}
\mathbf{E}\Big[\sup_{t\in[0,T]}|X(t,\cdot)|^{p}_{C_{x_{0},\theta}}\Big]\leq C_{p,d,\lambda,\gamma,\theta,T}\,A_{T}.
\end{equation}
\end{proposition}

Let us stress that the lower bound on $\theta$ in \eqref{eq: propSalins4} involves the H\"older exponent $\lambda$ of the space increments, and not the exponent $\gamma$ of the time increments. This is how the condition reads in \cite[Corollary~4.2]{S}, where the proof maps $\R^{d}$ onto the unit ball through $z\mapsto x_{0}+z(1-|z|)^{-1/\theta}$ and the exponent $\theta$ only interacts with the space variable. Let us also mention that the statement of \cite[Corollary~4.2]{S} does not display the dependence of its constant on $T$. This dependence cannot be avoided. 

We shall also use the version of this criterion for a random field in the space variable only, which is \cite[Theorem~4.1]{S}. We quote it in the same way.

\begin{proposition}[{\cite[Theorem~4.1]{S}}]\label{prop: Salins-space}
Let $X\colon\R^{d}\times\Omega\to\R$ be a measurable random field, let $x_{0}\in\R^{d}$, and assume that there exist $A>0$, $p>d$ and $\lambda\in(0,1)$ such that for all $x,y\in\R^{d}$,
\begin{equation}\label{eq: propSalins-space1}
\mathbf{E}\big[|X(x_{0})|^{p}\big]\leq A,
\qquad\text{and}\qquad
\mathbf{E}\big[|X(x)-X(y)|^{p}\big]\leq A\,|x-y|^{\lambda p}.
\end{equation}
Assume moreover that
\begin{equation}\label{eq: propSalins-space2}
p\lambda>d\qquad\text{and}\qquad\theta>\frac{p\lambda}{p-d}.
\end{equation}
Then $X$ admits a modification with paths in $C_{x_{0},\theta}(\R^{d})$, and there exists a constant $C_{p,d,\lambda,\theta}>0$, independent of $x_{0}$ and $A$, such that this modification satisfies
\begin{equation}\label{eq: propSalins-space3}
\mathbf{E}\Big[|X|^{p}_{C_{x_{0},\theta}}\Big]\leq C_{p,d,\lambda,\theta}\,A.
\end{equation}
\end{proposition}

The third tool is elementary. It allows to absorb a polynomial weight, and also the factor produced by the Laplacian bound \eqref{eq:laplaceG} of the heat kernel, into a heat kernel with a larger variance. Note that Lemma~\ref{lem:weight-absorb} below is proved in Appendix~\ref{app:weights}.

\begin{lemma}\label{lem:weight-absorb}
Let $\theta_{0}\geq 0$ and $a\in(0,1]$. There exists a constant $C=C_{d,\theta_{0},a}$ such that for all $\tau>0$ and $z\in\R^{d}$,
\begin{equation}\label{eq:absorb-power}
\bigl(1+|z|^{\theta_{0}}\bigr) \, G(\tau,z)^{a}\leq C\bigl(1+\tau^{\theta_{0}/2}\bigr)\,G(2\tau,z)^{a},
\end{equation}
and
\begin{equation}\label{eq:absorb-laplace}
\bigl(1+|z|^{\theta_{0}}\bigr)|\Delta_{z}G(\tau,z)|\leq\frac{C}{\tau}\bigl(1+\tau^{\theta_{0}/2}\bigr)\,G(2\tau,z).
\end{equation}
\end{lemma}

The fourth tool is also elementary. The weight inequality \eqref{eq:weight-triangle} of Lemma~\ref{lem:weight-elem} will be invoked repeatedly in the sequel, in order to transfer the weight from the integration variable $y$ to the point $x$ where the stochastic convolution is evaluated. We complement it with two bounds on the smooth weight $\hat{w}_{x_{0}}^{\theta_{0}}$ of \eqref{eq: weight1}, which we shall use in Section~\ref{subsec:weighted-conv-proof} because it can be differentiated. The proof of Lemma~\ref{lem:weights} below is deferred to Appendix~\ref{app:weights}.

\begin{lemma}\label{lem:weights}
Let $\theta_{0}\geq 0$ and $x_{0}\in\R^{d}$, and let $\hat{w}_{x_{0}}^{\theta_{0}}$ be the smooth weight of \eqref{eq: weight1}. Then $\hat{w}_{x_{0}}^{\theta_{0}}$ is differentiable, and for all $x,y\in\R^{d}$ we have
\begin{equation}\label{eq:what-grad}
|\nabla\hat{w}_{x_{0}}^{\theta_{0}}(x)|\leq\theta_{0}\,\hat{w}_{x_{0}}^{\theta_{0}}(x),
\qquad
\frac{\hat{w}_{x_{0}}^{\theta_{0}}(x)}{\hat{w}_{x_{0}}^{\theta_{0}}(y)}\leq\bigl(1+|x-y|\bigr)^{\theta_{0}}.
\end{equation}
\end{lemma}

\subsection{Proof of Proposition~\ref{prop:weighted-conv}}\label{subsec:weighted-conv-proof}

We now prove Proposition~\ref{prop:weighted-conv}, following the two steps mentioned in Section~\ref{subsec:weighted-conv} after its statement. The first step is Lemma~\ref{lem:weighted-moments} below, namely the moment bounds on the pointwise values and on the increments of $\mathcal{Z}^{h}$. The second step, which follows the statement of the lemma, is the Kolmogorov criterion of Proposition~\ref{prop: Salins} applied to a weighted field. In the three bounds below, the weight of the input appears as a factor $(1+|x-x_{0}|^{\theta_{0}})$ at the point of evaluation, and the powers of $T$ are the ones which will give the limit \eqref{eq:CT-to-zero}. The proof of Lemma~\ref{lem:weighted-moments} rests on a master inequality for stochastic convolutions with deterministic kernels, combined with the spectral bound of Lemma~\ref{lem: spectral} and the absorption inequalities of Lemma~\ref{lem:weight-absorb}. It is deferred to Appendix~\ref{app:conv}.

\begin{lemma}\label{lem:weighted-moments}
Let Assumption~\ref{assump:noise} hold, let $T>0$, $p\geq 2$, $x_{0}\in\R^{d}$ and $\theta_{0}\geq 0$. Let $h$ be an adapted random field on $[0,T]\times\R^{d}$ such that the quantity $N:=\mathcal{N}^{x_{0}}_{p,\theta_{0},T}(h)$ of \eqref{eq:Nfunctional} is finite, and let $Z:=\mathcal{Z}^{h}$ be the stochastic convolution \eqref{eq:StochConv}. Let $\chik\in(0,\eta)$ and $\lambda\in(0,\chik)$, and set $c_{T}:=(1+T)^{\theta_{0}+1}$. Then there exists a constant $C=C_{p,d,\theta_{0},\specm,\eta,\chik,\lambda}$, independent of $T$, $x_{0}$ and $h$, such that for all $t,t_{1},t_{2}\in[0,T]$ and $x,x_{1},x_{2}\in\R^{d}$ the following holds true.
\begin{enumerate}
\item\label{it:wm-point} \emph{Pointwise bound.} We have
\begin{equation}\label{eq:wm-point}
\|Z(t,x)\|_{L^{p}(\Omega)}\leq C\,c_{T}^{1/2}\,T^{\frac{\chik}{2}}\,N\,\bigl(1+|x-x_{0}|^{\theta_{0}}\bigr).
\end{equation}
\item\label{it:wm-space} \emph{Space increments.} We have
\begin{equation}\label{eq:wm-space}
\|Z(t,x_{1})-Z(t,x_{2})\|_{L^{p}(\Omega)}\leq C\,c_{T}^{1/2}\,T^{\frac{\chik-\lambda}{2}}\,N\,\bigl(1+|x_{1}-x_{0}|^{\theta_{0}}+|x_{2}-x_{0}|^{\theta_{0}}\bigr)\,|x_{1}-x_{2}|^{\lambda}.
\end{equation}
\item\label{it:wm-time} \emph{Time increments.} We have
\begin{equation}\label{eq:wm-time}
\|Z(t_{2},x)-Z(t_{1},x)\|_{L^{p}(\Omega)}\leq C\,c_{T}^{1/2}\,N\,\bigl(1+|x-x_{0}|^{\theta_{0}}\bigr)\,|t_{2}-t_{1}|^{\frac{\chik}{2}}.
\end{equation}
\end{enumerate}
\end{lemma}

Lemma~\ref{lem:weighted-moments} is the first step of the proof of Proposition~\ref{prop:weighted-conv}. We now turn to the second one, and apply the Kolmogorov criterion of Proposition~\ref{prop: Salins} to the weighted field
\begin{equation}\label{eq:Ztilde}
\tilde{Z}(t,x):=\hat{w}_{x_{0}}^{\theta_{0}}(x)\,Z(t,x),\qquad(t,x)\in[0,T]\times\R^{d},
\end{equation}
where $Z=\mathcal{Z}^{h}$ and $\hat{w}_{x_{0}}^{\theta_{0}}$ is the smooth weight of \eqref{eq: weight1}. Observe that the weight $\hat{w}_{x_{0}}^{\theta_{0}}$ is chosen, rather than $w_{x_{0}}^{\theta_{0}}$, because it can be differentiated, see Lemma~\ref{lem:weights}. We also recall from \eqref{eq: universalC1} that $\hat{w}_{x_{0}}^{\theta_{0}}\leq C_{1}w_{x_{0}}^{\theta_{0}}$, so that
\begin{equation}\label{eq:what-kills-weight}
\hat{w}_{x_{0}}^{\theta_{0}}(x)\bigl(1+|x-x_{0}|^{\theta_{0}}\bigr)\leq C_{1},\qquad x\in\R^{d}.
\end{equation}

\begin{proof}[Proof of Proposition~\ref{prop:weighted-conv}]
We write $N:=\mathcal{N}^{x_{0}}_{p,\theta_{0},T}(h)$ and $c_{T}:=(1+T)^{\theta_{0}+1}$ as in Lemma~\ref{lem:weighted-moments}, and we may assume $N>0$. The proof is divided in three steps. We first choose the exponents, then verify the hypotheses of Proposition~\ref{prop: Salins} for the field $\tilde{Z}$ of \eqref{eq:Ztilde}, and finally go back to $Z$.

\smallskip

\noindent
{\it Step 1: choice of the exponents.} We fix $\lambda\in(0,\min\{\eta,\theta_{1}\})$, $\chik\in(\lambda,\eta)$ and $\varepsilon\in(0,\chik)$, and we set $\gamma:=(\chik-\varepsilon)/2$. These exponents do not depend on $p$, and we have $\chik\in(0,\eta)$ and $\lambda<\chik$ as required in Lemma~\ref{lem:weighted-moments}. Since $p\geq p^{*}$, we may assume that
\begin{equation}\label{eq:choice-lambda}
p\lambda>d+1,\qquad p\gamma>d+1\qquad\text{and}\qquad\theta_{1}>\frac{p\lambda}{p-(d+1)} ,
\end{equation}
where the last condition holds for $p$ large since $p\lambda/(p-(d+1))\to\lambda<\theta_{1}$. Observe that relation \eqref{eq:choice-lambda} is condition \eqref{eq: propSalins4} with $\theta_{1}$ in place of $\theta$.

\smallskip

\noindent
{\it Step 2: hypotheses of Proposition~\ref{prop: Salins} for $\tilde{Z}$.} Let $C$ be the constant of Lemma~\ref{lem:weighted-moments} for the exponents $\chik,\lambda$ of Step~1, and set
\begin{equation}\label{eq:AT-def}
A_{T}:=C^{p}\,(2C_{1})^{p}\,2^{p\theta_{0}}(1+\theta_{0})^{p}\,c_{\theta_{0}}^{p}\,c_{T}^{p/2}\Bigl(T^{\frac{\chik}{2}}+T^{\frac{\chik-\lambda}{2}}+T^{\frac{\varepsilon}{2}}\Bigr)^{p}N^{p}.
\end{equation}
Every power of $T$ in \eqref{eq:AT-def} is positive, so that $A_{T}/N^{p}\to 0$ as $T\to 0$, and $A_{T}$ does not depend on $x_{0}$. We check \eqref{eq: propSalins1}--\eqref{eq: propSalins3} with this $A_{T}$ and the exponents $\lambda,\gamma$ of Step~1.

Relation \eqref{eq: propSalins1} is trivial, since $Z(0,\cdot)=0$. For the time increments \eqref{eq: propSalins3}, we multiply~\eqref{eq:wm-time} by $\hat{w}_{x_{0}}^{\theta_{0}}(x)$ and use \eqref{eq:what-kills-weight}, which gives
\begin{equation}\label{eq:Ztilde-time}
\|\tilde{Z}(t_{2},x)-\tilde{Z}(t_{1},x)\|_{L^{p}(\Omega)}\leq C\,C_{1}\,c_{T}^{1/2}\,N\,|t_{2}-t_{1}|^{\frac{\chik}{2}}
\leq C\,C_{1}\,c_{T}^{1/2}\,T^{\frac{\varepsilon}{2}}\,N\,|t_{2}-t_{1}|^{\gamma},
\end{equation}
where we have used $|t_{2}-t_{1}|\leq T$ in the second inequality. This implies \eqref{eq: propSalins3}. Similarly, multiplying~\eqref{eq:wm-point} by $\hat{w}_{x_{0}}^{\theta_{0}}(x)$ gives the uniform bound
\begin{equation}\label{eq:Ztilde-point}
\sup_{t\in[0,T]}\sup_{x\in\R^{d}}\|\tilde{Z}(t,x)\|_{L^{p}(\Omega)}\leq C\,C_{1}\,c_{T}^{1/2}\,T^{\frac{\chik}{2}}\,N,
\end{equation}
which will be employed later.

We now turn to the space increments \eqref{eq: propSalins2}. Let $t\in[0,T]$ and $x_{1},x_{2}\in\R^{d}$. If $|x_{1}-x_{2}|>1$, then \eqref{eq:Ztilde-point} and the triangle inequality give
\begin{equation}\label{eq:Ztilde-space-far}
\|\tilde{Z}(t,x_{1})-\tilde{Z}(t,x_{2})\|_{L^{p}(\Omega)}\leq 2C\,C_{1}\,c_{T}^{1/2}\,T^{\frac{\chik}{2}}\,N\,|x_{1}-x_{2}|^{\lambda}.
\end{equation}
If $|x_{1}-x_{2}|\leq 1$, we write
\begin{equation}\label{eq:Ztilde-product}
\tilde{Z}(t,x_{1})-\tilde{Z}(t,x_{2})=\hat{w}_{x_{0}}^{\theta_{0}}(x_{1})\bigl(Z(t,x_{1})-Z(t,x_{2})\bigr)+\bigl(\hat{w}_{x_{0}}^{\theta_{0}}(x_{1})-\hat{w}_{x_{0}}^{\theta_{0}}(x_{2})\bigr)Z(t,x_{2}),
\end{equation}
and we bound the two terms separately. For the first one, we apply \eqref{eq:weight-triangle} with $x=x_{1}$ and $y=x_{2}$, and then $|x_{1}-x_{2}|\leq 1$, to get $1+|x_{2}-x_{0}|^{\theta_{0}}\leq 2c_{\theta_{0}}(1+|x_{1}-x_{0}|^{\theta_{0}})$. Hence \eqref{eq:wm-space} and \eqref{eq:what-kills-weight} yield
\begin{equation}\label{eq:Ztilde-space-1}
\hat{w}_{x_{0}}^{\theta_{0}}(x_{1})\,\|Z(t,x_{1})-Z(t,x_{2})\|_{L^{p}(\Omega)}\leq 2c_{\theta_{0}}\,C\,C_{1}\,c_{T}^{1/2}\,T^{\frac{\chik-\lambda}{2}}\,N\,|x_{1}-x_{2}|^{\lambda}.
\end{equation}
For the second term in \eqref{eq:Ztilde-product}, the mean value theorem and the first relation in \eqref{eq:what-grad} give $|\hat{w}_{x_{0}}^{\theta_{0}}(x_{1})-\hat{w}_{x_{0}}^{\theta_{0}}(x_{2})|\leq\theta_{0}|x_{1}-x_{2}|\sup_{\zeta}\hat{w}_{x_{0}}^{\theta_{0}}(\zeta)$, the supremum being over the segment $[x_{1},x_{2}]$. On this segment, since we are assuming $|x_{1}-x_{2}|\leq 1$, we have $|\zeta-x_{2}|\leq 1$. Therefore the second relation in \eqref{eq:what-grad} gives $\hat{w}_{x_{0}}^{\theta_{0}}(\zeta)\leq 2^{\theta_{0}}\hat{w}_{x_{0}}^{\theta_{0}}(x_{2})$. Therefore, by \eqref{eq:Ztilde-point} and $|x_{1}-x_{2}|\leq|x_{1}-x_{2}|^{\lambda}$,
\begin{equation}\label{eq:Ztilde-space-2}
\bigl|\hat{w}_{x_{0}}^{\theta_{0}}(x_{1})-\hat{w}_{x_{0}}^{\theta_{0}}(x_{2})\bigr|\,\|Z(t,x_{2})\|_{L^{p}(\Omega)}\leq\theta_{0}2^{\theta_{0}}\,C\,C_{1}\,c_{T}^{1/2}\,T^{\frac{\chik}{2}}\,N\,|x_{1}-x_{2}|^{\lambda}.
\end{equation}
Gathering \eqref{eq:Ztilde-space-far}, \eqref{eq:Ztilde-product}, \eqref{eq:Ztilde-space-1} and \eqref{eq:Ztilde-space-2}, and raising to the power $p$, we obtain \eqref{eq: propSalins2} with the constant $A_{T}$ of \eqref{eq:AT-def}. The field $\tilde{Z}$ is jointly measurable, as $Z$ is.

Let us now explain briefly why the exponent $\theta_{1}$ enters at this point. The moment bounds \eqref{eq:Ztilde-time}--\eqref{eq:Ztilde-space-2} on $\tilde{Z}$ are uniform in $x$, since the weight $\hat{w}_{x_{0}}^{\theta_{0}}$ has absorbed the factor $(1+|x-x_{0}|^{\theta_{0}})$ of Lemma~\ref{lem:weighted-moments}. However, Proposition~\ref{prop: Salins} does not produce a bounded modification out of uniform moment bounds. It produces a modification in a weighted space $C_{x_{0},\theta}([0,T]\times\R^{d})$, with an exponent $\theta$ subject to the last condition in \eqref{eq: propSalins4}. This is where the loss $\theta_{1}$ comes from. Namely, the exponents satisfy \eqref{eq: propSalins4} with $\theta_{1}$ in place of $\theta$, by \eqref{eq:choice-lambda}. Proposition~\ref{prop: Salins} thus provides a modification of $\tilde{Z}$ with paths in $C_{x_{0},\theta_{1}}([0,T]\times\R^{d})$, such that
\begin{equation}\label{eq:Ztilde-conclusion}
\mathbf{E}\Big[\sup_{t\in[0,T]}\sup_{x\in\R^{d}}\bigl(w_{x_{0}}^{\theta_{1}}(x)\,|\tilde{Z}(t,x)|\bigr)^{p}\Big]\leq C_{p,d,\lambda,\gamma,\theta_{1},T}\,A_{T}.
\end{equation}

\smallskip

\noindent
{\it Step 3: back to $Z$.} Dividing the modification of $\tilde{Z}$ by the continuous positive function $\hat{w}_{x_{0}}^{\theta_{0}}$ gives a modification of $Z$, which we still denote by $Z$, with continuous paths on $[0,T]\times\R^{d}$. It is then easy to change the weights $\hat{w}$ for $w$. Namely applying~\eqref{eq:Ztilde-conclusion} and the elementary bound~\eqref{eq: universalC1}, one can easily change the weight $\hat{w}_{x_{0}}^{\theta_{0}}$ defining $\tilde{Z}$ for $w_{x_{0}}^{\theta_{0}}$. We get that 
 $Z$ has paths in $C_{x_{0},\theta}([0,T]\times\R^{d})$, and \eqref{eq:Ztilde-conclusion} yields
\begin{equation}\label{eq:Z-conclusion}
\mathbf{E}\Big[\sup_{t\in[0,T]}|Z(t,\cdot)|^{p}_{C_{x_{0},\theta}}\Big]\leq(4C_{1})^{p}\,C_{p,d,\lambda,\gamma,\theta_{1},T}\,A_{T},
\end{equation}
which is \eqref{eq:weighted-conv} with $C_{T}:=(4C_{1})^{p}C_{p,d,\lambda,\gamma,\theta_{1},T}A_{T}/N^{p}$. This constant depends on $T$, $p$, $\theta_{0}$ and $\theta_{1}$ only (the exponents $\lambda,\chik,\varepsilon,\gamma$ having been chosen as functions of $p$, $\theta_{1}$ and $\eta$ in Step~1), it does not depend on $x_{0}$ nor on $h$, and it satisfies \eqref{eq:CT-to-zero} by the remark following \eqref{eq:AT-def}, the constant of Proposition~\ref{prop: Salins} being nondecreasing in $T$. The proof is complete.
\end{proof}

\subsection{A weighted-center norm}\label{subsec:X-weighted}

The contraction argument of Section~\ref{sec:extend-proof} runs in the translate $\mathcal{Y}_{\tau}$ of \eqref{eq:Yspace} of the space $\mathcal{X}_{\tau}$ of \cite[Section~6]{S}. Its distance is the norm \eqref{eq:Xnorm}, which places a supremum over all the centers $x_{0}$ outside the expectation. The continuity argument of Section~\ref{subsec:continuity} below runs in a variant of this space, where the centers far from the origin are penalized by a polynomial weight. We shall see that this penalization is what absorbs the loss of exponent $\theta_{1}$ of Proposition~\ref{prop:weighted-conv}.

\begin{definition}\label{def:X-weighted}
Let $\theta>0$, $p\geq 1$, $\theta_{0}\in[0,\theta]$ and $\tau>0$. For an adapted random field $Y$ on $[0,\tau]\times\R^{d}$ with paths in $C_{0,\theta}([0,\tau]\times\R^{d})$, we set
\begin{equation}\label{eq:X-weighted}
\|Y\|_{\mathcal{X}^{\theta_{0}}_{\tau}}:=\sup_{x_{0}\in\R^{d}}w^{\theta_{0}}_{0}(x_{0})\,\Big(\mathbf{E}\Big[\sup_{t\in[0,\tau]}|Y(t,\cdot)|^{p}_{C_{x_{0},\theta}}\Big]\Big)^{1/p},
\end{equation}
where $w_{0}^{\theta_{0}}(x_{0})=(1+|x_{0}|^{\theta_{0}})^{-1}$ is the weight of \eqref{eq: weight1} centered at the origin. For a random variable $\psi$ with values in $C_{\theta}$ we also set
\begin{equation}\label{eq:X-weighted-0}
\|\psi\|_{\mathcal{X}^{\theta_{0}}_{0}}:=\sup_{x_{0}\in\R^{d}}w^{\theta_{0}}_{0}(x_{0})\,\big(\mathbf{E}\big[|\psi|^{p}_{C_{x_{0},\theta}}\big]\big)^{1/p}.
\end{equation}
The quantities \eqref{eq:X-weighted} and \eqref{eq:X-weighted-0} take their values in $[0,\infty]$.
\end{definition}

For $\theta_{0}=0$, relation \eqref{eq:X-weighted} is the norm \eqref{eq:Xnorm} used in the proof of Theorem~\ref{thm:extend}. For $\theta_{0}>0$ it is a weaker norm, in which the moments at a center $x_{0}$ are divided by $1+|x_{0}|^{\theta_{0}}$. Since we are manipulating a substantial number of different norms, let us record how they compare. Item \ref{it:X-pointwise} below is the key to the Gronwall type argument of Section~\ref{subsec:continuity}. Namely it converts the weighted pointwise moments of $Y$, which are the input of Proposition~\ref{prop:weighted-conv}, into the norm \eqref{eq:X-weighted}. Let us also stress that the exponent $\theta$ of the spaces is present in every item, through the norms $|\cdot|_{C_{x_{0},\theta}}$ in \eqref{eq:X-weighted} and \eqref{eq:X-weighted-0}. In item \ref{it:X-initial}, this exponent is traded for the smaller exponent $\theta_{0}$ of the space where the data live. This is what allows Theorem~\ref{thm:continuity} to state its conclusion for data in $C_{\theta_{0}}$.

\begin{lemma}\label{lem:X-weighted}
Let $\theta,p,\theta_{0},\tau$ be as in Definition~\ref{def:X-weighted}, and let $c_{\theta_{0}}=\max\{1,2^{\theta_{0}-1}\}$ be the constant of Lemma~\ref{lem:weight-elem}. Then the following holds true.
\begin{enumerate}
\item\label{it:X-compare} For every $Y$ as in Definition~\ref{def:X-weighted},
\begin{equation}\label{eq:X-compare}
\Big(\mathbf{E}\Big[\sup_{t\in[0,\tau]}|Y(t,\cdot)|^{p}_{C_{\theta}}\Big]\Big)^{1/p}\leq\|Y\|_{\mathcal{X}^{\theta_{0}}_{\tau}}\leq\|Y\|_{\mathcal{X}_{\tau}}.
\end{equation}
\item\label{it:X-initial} If $\psi$ is a random variable with values in $C_{\theta_{0}}$, then
\begin{equation}\label{eq:X-initial}
\|\psi\|_{\mathcal{X}^{\theta_{0}}_{0}}\leq 2c_{\theta_{0}}\,\big(\mathbf{E}\big[|\psi|^{p}_{C_{\theta_{0}}}\big]\big)^{1/p}.
\end{equation}
\item\label{it:X-pointwise} For every $Y$ as in Definition~\ref{def:X-weighted} and every $(s,y)\in[0,\tau]\times\R^{d}$,
\begin{equation}\label{eq:X-pointwise}
\|Y(s,y)\|_{L^{p}(\Omega)}\leq\big(1+|y|^{\theta_{0}}\big)\,\|Y\|_{\mathcal{X}^{\theta_{0}}_{\tau}}.
\end{equation}
\item\label{it:X-terminal} For every $Y$ as in Definition~\ref{def:X-weighted}, the terminal value satisfies $\|Y(\tau,\cdot)\|_{\mathcal{X}^{\theta_{0}}_{0}}\leq\|Y\|_{\mathcal{X}^{\theta_{0}}_{\tau}}$.
\end{enumerate}
\end{lemma}

\begin{proof}
The left inequality in \eqref{eq:X-compare} is just a consequence of the definition of $\|\cdot\|_{\mathcal{X}_{\tau}^{\theta_{0}}}$, and the right one holds because $w_{0}^{\theta_{0}}\leq 1$. As for \eqref{eq:X-initial}, since $\theta_{0}\leq\theta$ we have $w^{\theta}_{x_{0}}\leq 2w^{\theta_{0}}_{x_{0}}$ by \eqref{eq:weight-two-exp}, and relation \eqref{eq:weight-triangle} with $(x,y,x_{0})$ replaced by $(x_{0},x,0)$ gives $1+|x|^{\theta_{0}}\leq c_{\theta_{0}}(1+|x_{0}|^{\theta_{0}})(1+|x-x_{0}|^{\theta_{0}})$. Hence
\begin{equation*}
|\psi|_{C_{x_{0},\theta}}\leq 2\,|\psi|_{C_{x_{0},\theta_{0}}}\leq 2c_{\theta_{0}}\big(1+|x_{0}|^{\theta_{0}}\big)\,|\psi|_{C_{\theta_{0}}},\qquad x_{0}\in\R^{d},
\end{equation*}
and \eqref{eq:X-initial} follows by multiplying by $w_{0}^{\theta_{0}}(x_{0})$ and taking the $L^{p}(\Omega)$ norm. For \eqref{eq:X-pointwise}, the pointwise bound \eqref{eq: pointwiseNorm} with center $y$ gives, for every $s\in[0,\tau]$,
\begin{equation*}
|Y(s,y)|\leq|Y(s,\cdot)|_{C_{y,\theta}}\leq\sup_{t\in[0,\tau]}|Y(t,\cdot)|_{C_{y,\theta}} .
\end{equation*}
Taking the $L^{p}(\Omega)$ norm, and then multiplying and dividing by $w_{0}^{\theta_{0}}(y)$, we get
\begin{equation*}
\|Y(s,y)\|_{L^{p}(\Omega)}\leq\big(1+|y|^{\theta_{0}}\big)\,w_{0}^{\theta_{0}}(y)\Big(\mathbf{E}\Big[\sup_{t\in[0,\tau]}|Y(t,\cdot)|^{p}_{C_{y,\theta}}\Big]\Big)^{1/p}\leq\big(1+|y|^{\theta_{0}}\big)\,\|Y\|_{\mathcal{X}^{\theta_{0}}_{\tau}}.
\end{equation*}
Item \ref{it:X-terminal} is obvious.
\end{proof}
\subsection{Continuity with respect to the initial condition}\label{subsec:continuity}

We now combine Proposition~\ref{prop:weighted-conv} with the fixed point relation \eqref{eq:fixedpoint}. The point is that the loss of exponent $\theta_{1}$ in \eqref{eq:weighted-conv} is absorbed by the norm of Definition~\ref{def:X-weighted}, whose exponent $\theta_{0}$ on the centers is smaller than the exponent $\theta$ of the spaces.

\begin{theorem}\label{thm:continuity}
Let Assumptions~\ref{assump:noise}, \ref{assump:reac} and \ref{assump:diff} hold, let $\theta>0$ and $\theta_{0}\in[0,\theta)$, and let $p\geq p^{*}$. Let $\xi_{1},\xi_{2}$ be two initial conditions satisfying Assumption~\ref{assump: initial-cond-fixed}. Let $u_{i}:=u(\cdot\,;\xi_{i})$ be the solutions given by Theorem~\ref{thm:extend}, and set $D:=u_{1}-u_{2}$. Then for every $T>0$ there exists a constant $C_{T}$, depending on $T,p,d,\theta,\theta_{0},\eta,\specm,\kappa$ and $L_{\sigma}$ but not on $\xi_{1},\xi_{2}$, such that
\begin{equation}\label{eq:continuity}
\|D\|_{\mathcal{X}^{\theta_{0}}_{T}}\leq C_{T}\,\|\xi_{1}-\xi_{2}\|_{\mathcal{X}^{\theta_{0}}_{0}},
\end{equation}
where the norms are defined by \eqref{eq:X-weighted} and \eqref{eq:X-weighted-0}. In particular, if $\xi_{1}-\xi_{2}$ takes its values in $C_{\theta_{0}}$, then
\begin{equation}\label{eq:continuity-fixed}
\mathbf{E}\Big[\sup_{t\in[0,T]}|u(t,\cdot\,;\xi_{1})-u(t,\cdot\,;\xi_{2})|^{p}_{C_{\theta}}\Big]\leq\big(2c_{\theta_{0}}C_{T}\big)^{p}\,\mathbf{E}\big[|\xi_{1}-\xi_{2}|^{p}_{C_{\theta_{0}}}\big].
\end{equation}
\end{theorem}

Before the proof, let us comment on the two extreme values of $\theta_{0}$. For $\theta_{0}=0$ and for deterministic data, the norm $\|\phi_{1}-\phi_{2}\|_{\mathcal{X}^{0}_{0}}$ is the supremum norm $|\phi_{1}-\phi_{2}|_{C_{b}}$. This is the continuity given by the contraction argument, as mentioned in item (ii) of Remark~\ref{rmk:markov-property}. The value $\theta_{0}=\theta$ is excluded whatever the value of $p$ is, since the loss $\theta_{1}=\theta-\theta_{0}$ of Proposition~\ref{prop:weighted-conv} is positive, and this is where the Feller property on $C_{\theta}$ would sit. Relation \eqref{eq:continuity-fixed} is thus a Lipschitz continuity from $C_{\theta_{0}}$ to $L^{p}(\Omega;C([0,T];C_{\theta}))$, with a loss of exponent between the input and the output.

\begin{proof}
If the right hand side of \eqref{eq:continuity} is infinite, there is nothing to prove. We thus assume throughout the proof that $\|\xi_{1}-\xi_{2}\|_{\mathcal{X}^{\theta_{0}}_{0}}<\infty$. The value of the constant $C$ may change from line to line, but always depends only on the parameters listed in the statement. We fix $T>0$. For $i=1,2$ we write $U_{0}^{i}:=H\xi_{i}$ for the field \eqref{eq:U0} with initial condition $\xi_{i}$, and $Z_{i}:=\mathcal{Z}^{\sigma(u_{i})}$ for the stochastic convolution \eqref{eq:StochConv} with integrand $\sigma(u_{i})$. The fixed point relation \eqref{eq:fixedpoint} then reads
\begin{equation}\label{eq:cont-fixedpoint}
u_{i}=\mathcal{M}\big(U_{0}^{i}+Z_{i}\big),\qquad\text{that is}\qquad u_{i}=\Gamma_{i}(u_{i}),
\end{equation}
where $\Gamma_{i}$ is the Picard map \eqref{eq:Gamma} associated with the initial condition $\xi_{i}$. The proof is divided in three steps. We first check that $\|D\|_{\mathcal{X}^{\theta_{0}}_{T}}$ is finite, then we prove \eqref{eq:continuity} on a short interval, and finally we iterate thanks to the restart property of Lemma~\ref{lem:restart}.

\smallskip

\noindent
{\it Step 1: A priori finiteness.} We show that $\|D\|_{\mathcal{X}^{\theta_{0}}_{T}}<\infty$. Notice that each $\|u_{i}\|_{\mathcal{X}^{\theta_{0}}_{T}}$ may be infinite, since no moment is assumed on $\xi_{i}$. We thus decompose $D$ in such a way that the initial conditions only enter through their difference, namely
\begin{equation*}
D=\big(u_{1}-\mathcal{M}(U_{0}^{1})\big)-\big(u_{2}-\mathcal{M}(U_{0}^{2})\big)+\big(\mathcal{M}(U_{0}^{1})-\mathcal{M}(U_{0}^{2})\big).
\end{equation*}
For the first two terms, recall that $u_{i}=\Gamma_{i}(u_{i})$ by \eqref{eq:cont-fixedpoint}. Hence relation \eqref{eq:Gamma-affine} gives, for every $x_{0}\in\R^{d}$, almost surely
\begin{equation*}
\sup_{t\in[0,T]}\big|u_{i}(t,\cdot)-\mathcal{M}(U_{0}^{i})(t,\cdot)\big|_{C_{x_{0},\theta}}\leq C_{\theta,\kappa}e^{C_{\theta,\kappa}T}\sup_{t\in[0,T]}|Z_{i}(t,\cdot)|_{C_{x_{0},\theta}} .
\end{equation*}
Since $\sigma$ is bounded, the bound \eqref{eq:StochConvSigma} on stochastic convolutions with bounded integrands controls the $p$-th moment of the right hand side uniformly in $x_{0}$, whatever the initial condition $\xi_{i}$ is. Therefore $u_{i}-\mathcal{M}(U_{0}^{i})$ belongs to $\mathcal{X}_{T}$, and its $\mathcal{X}^{\theta_{0}}_{T}$ norm is finite by \eqref{eq:X-compare}. As for the last term, the Lipschitz bound \eqref{eq:MLip} and the bound \eqref{eq:cont-free} below on the heat semigroup give, for every $x_{0}\in\R^{d}$, almost surely
\begin{equation*}
\sup_{t\in[0,T]}\big|\mathcal{M}(U_{0}^{1})(t,\cdot)-\mathcal{M}(U_{0}^{2})(t,\cdot)\big|_{C_{x_{0},\theta}}\leq C\,|\xi_{1}-\xi_{2}|_{C_{x_{0},\theta}} .
\end{equation*}
Multiplying by $w_{0}^{\theta_{0}}(x_{0})$, taking the $L^{p}(\Omega)$ norm and the supremum over $x_{0}$, the $\mathcal{X}^{\theta_{0}}_{T}$ norm of the last term is bounded by $C\,\|\xi_{1}-\xi_{2}\|_{\mathcal{X}^{\theta_{0}}_{0}}$, which is finite. By the triangle inequality, we conclude that $\|D\|_{\mathcal{X}^{\theta_{0}}_{T}}<\infty$. Let us also mention that the boundedness of $\sigma$ is only used here, to get the finiteness of $\|D\|_{\mathcal{X}^{\theta_{0}}_{T}}$.

\smallskip

\noindent
{\it Step 2: Short time estimate.} Let $\tau\in(0,T]$. The fixed point relation \eqref{eq:cont-fixedpoint} on the interval $[0,\tau]$, the Lipschitz bound \eqref{eq:MLip} with center $x_{0}$ and the triangle inequality give almost surely
\begin{equation}\label{eq:cont-D-basic}
\sup_{t\in[0,\tau]}|D(t,\cdot)|_{C_{x_{0},\theta}}\leq C_{\theta,\kappa}e^{C_{\theta,\kappa}\tau}\Big(\sup_{t\in[0,\tau]}|U_{0}^{1}(t,\cdot)-U_{0}^{2}(t,\cdot)|_{C_{x_{0},\theta}}+\sup_{t\in[0,\tau]}|Z_{1}(t,\cdot)-Z_{2}(t,\cdot)|_{C_{x_{0},\theta}}\Big).
\end{equation}
We multiply \eqref{eq:cont-D-basic} by $w_{0}^{\theta_{0}}(x_{0})$, take the $L^{p}(\Omega)$ norm and the supremum over $x_{0}$. For the first term, we invoke the bound \eqref{eq:U0Bound}, which holds for every center $x_{0}$ with the same constants since the proof of Lemma~\ref{lem:fc-heat} is invariant under translations. Applied to $U_{0}^{1}-U_{0}^{2}=H(\xi_{1}-\xi_{2})$, it gives almost surely
\begin{equation}\label{eq:cont-free}
\sup_{t\in[0,\tau]}|U_{0}^{1}(t,\cdot)-U_{0}^{2}(t,\cdot)|_{C_{x_{0},\theta}}\leq c_{\theta}\big(1+m_{\theta}T^{\theta/2}\big)\,|\xi_{1}-\xi_{2}|_{C_{x_{0},\theta}} .
\end{equation}
Multiplying \eqref{eq:cont-free} by $w_{0}^{\theta_{0}}(x_{0})$, taking the $L^{p}(\Omega)$ norm and the supremum over $x_{0}$, the first term of \eqref{eq:cont-D-basic} is thus bounded by $C\,\|\xi_{1}-\xi_{2}\|_{\mathcal{X}^{\theta_{0}}_{0}}$, and we obtain
\begin{equation}\label{eq:cont-short-1}
\|D\|_{\mathcal{X}^{\theta_{0}}_{\tau}}
\leq C\,\|\xi_{1}-\xi_{2}\|_{\mathcal{X}^{\theta_{0}}_{0}}
+C\sup_{x_{0}\in\R^{d}}w_{0}^{\theta_{0}}(x_{0})\,\Big(\mathbf{E}\Big[\sup_{t\in[0,\tau]}|Z_{1}(t,\cdot)-Z_{2}(t,\cdot)|^{p}_{C_{x_{0},\theta}}\Big]\Big)^{1/p} .
\end{equation}
We now bound the stochastic term in \eqref{eq:cont-short-1} for a fixed center $x_{0}$. By item \ref{it:conv-lin} of Remark~\ref{rmk:conv-bounded}, $Z_{1}-Z_{2}=\mathcal{Z}^{h}$ with $h:=\sigma(u_{1})-\sigma(u_{2})$, which is an adapted random field. We apply Proposition~\ref{prop:weighted-conv} to $h$ with the exponents $\theta_{0}$ and $\theta_{1}:=\theta-\theta_{0}$, which is allowed since $p\geq p^{*}$. We get
\begin{equation}\label{eq:cont-short-2}
\Big(\mathbf{E}\Big[\sup_{t\in[0,\tau]}|Z_{1}(t,\cdot)-Z_{2}(t,\cdot)|^{p}_{C_{x_{0},\theta}}\Big]\Big)^{1/p}\leq C_{\tau}^{1/p}\,\mathcal{N}^{x_{0}}_{p,\theta_{0},\tau}(h),
\end{equation}
where $C_{\tau}$ is the constant of \eqref{eq:weighted-conv}, which does not depend on $x_{0}$ and vanishes as $\tau\to 0$ by~\eqref{eq:CT-to-zero}. Next we bound the weighted functional \eqref{eq:Nfunctional} of $h$. By the Lipschitz property of $\sigma$ in Assumption~\ref{assump:diff} and then by \eqref{eq:X-pointwise}, for every $(s,y)\in[0,\tau]\times\R^{d}$ we have
\begin{equation}\label{eq:cont-short-3}
w_{x_{0}}^{\theta_{0}}(y)\,\|h(s,y)\|_{L^{p}(\Omega)}
\leq L_{\sigma}\,w_{x_{0}}^{\theta_{0}}(y)\,\|D(s,y)\|_{L^{p}(\Omega)}
\leq L_{\sigma}\,w_{x_{0}}^{\theta_{0}}(y)\big(1+|y|^{\theta_{0}}\big)\,\|D\|_{\mathcal{X}^{\theta_{0}}_{\tau}} .
\end{equation}
Relation \eqref{eq:weight-triangle} with $(x,y,x_{0})$ replaced by $(x_{0},y,0)$ gives $1+|y|^{\theta_{0}}\leq c_{\theta_{0}}(1+|x_{0}|^{\theta_{0}})(1+|y-x_{0}|^{\theta_{0}})$, that is
\begin{equation}\label{eq:cont-short-4}
w_{0}^{\theta_{0}}(x_{0})\,w_{x_{0}}^{\theta_{0}}(y)\big(1+|y|^{\theta_{0}}\big)\leq c_{\theta_{0}},\qquad x_{0},y\in\R^{d}.
\end{equation}
Taking the supremum over $(s,y)$ in \eqref{eq:cont-short-3}, multiplying by $w_{0}^{\theta_{0}}(x_{0})$ and using \eqref{eq:cont-short-4}, we obtain
\begin{equation}\label{eq:cont-short-5}
w_{0}^{\theta_{0}}(x_{0})\,\mathcal{N}^{x_{0}}_{p,\theta_{0},\tau}(h)\leq c_{\theta_{0}}L_{\sigma}\,\|D\|_{\mathcal{X}^{\theta_{0}}_{\tau}},
\end{equation}
uniformly in $x_{0}$. This is the step where the weighted-center norm is needed. Now plugging~\eqref{eq:cont-short-5} into \eqref{eq:cont-short-2}, then into \eqref{eq:cont-short-1}, we end up with
\begin{equation}\label{eq:cont-short-6}
\|D\|_{\mathcal{X}^{\theta_{0}}_{\tau}}
\leq C\,\|\xi_{1}-\xi_{2}\|_{\mathcal{X}^{\theta_{0}}_{0}}+C\,C_{\tau}^{1/p}c_{\theta_{0}}L_{\sigma}\,\|D\|_{\mathcal{X}^{\theta_{0}}_{\tau}} .
\end{equation}
Owing to \eqref{eq:CT-to-zero}, we can choose an integer $N\geq 1$ such that $T_{1}:=T/N$ satisfies $C\,C_{T_{1}}^{1/p}c_{\theta_{0}}L_{\sigma}\leq 1/2$. Since $T_{1}$ only depends on the parameters listed in the statement, it does not depend on $\xi_{1},\xi_{2}$. By Step~1, the last term in \eqref{eq:cont-short-6} with $\tau=T_{1}$ is finite and can be absorbed into the left hand side, so that
\begin{equation}\label{eq:cont-short-7}
\|D\|_{\mathcal{X}^{\theta_{0}}_{T_{1}}}\leq C\,\|\xi_{1}-\xi_{2}\|_{\mathcal{X}^{\theta_{0}}_{0}} .
\end{equation}

\smallskip

\noindent
{\it Step 3: Iteration.} This is a standard patching argument, and we only sketch it. Let $k\in\{1,\ldots,N-1\}$. By Lemma~\ref{lem:restart}, the fields $u_{i}(kT_{1}+\cdot)$ are the solutions given by Theorem~\ref{thm:extend} on $[0,T_{1}]$ for the shifted noise $W^{kT_{1}}$ of \eqref{eq:shifted-noise}, with the $\mathcal{F}_{kT_{1}}$-measurable initial conditions $u_{i}(kT_{1},\cdot)$. Steps~1 and 2 apply to them and give
\begin{equation}\label{eq:cont-iter}
\|D(kT_{1}+\cdot)\|_{\mathcal{X}^{\theta_{0}}_{T_{1}}}\leq C\,\|D(kT_{1},\cdot)\|_{\mathcal{X}^{\theta_{0}}_{0}}\leq C\,\|D((k-1)T_{1}+\cdot)\|_{\mathcal{X}^{\theta_{0}}_{T_{1}}},
\end{equation}
where the second inequality is item \ref{it:X-terminal} of Lemma~\ref{lem:X-weighted}. Let us stress that this is the second place where the weighted-center norm is needed. Namely, the output \eqref{eq:cont-short-7} of Step~2 is a weighted-center quantity, and it is accepted as input by Step~2 itself, with the same exponents $\theta_{0}$ and $\theta$. Iterating \eqref{eq:cont-iter} from \eqref{eq:cont-short-7}, and then bounding the supremum over $[0,T]$ in \eqref{eq:X-weighted} by the sum of the suprema over the $N$ intervals $[kT_{1},(k+1)T_{1}]$, we get \eqref{eq:continuity}. Relation \eqref{eq:continuity-fixed} follows from \eqref{eq:continuity}, the left inequality in \eqref{eq:X-compare} and \eqref{eq:X-initial}. The proof is complete.
\end{proof}

\subsection{Measurability and a Feller property between two exponents}\label{subsec:feller-two}

We now draw the consequences of Theorem~\ref{thm:continuity} for the transition semigroup $\mathcal{P}_{t}$ of Definition~\ref{def: transition}. Recall from item (ii) of Remark~\ref{rmk:markov-property} that relation \eqref{eq: transition} requires the map $\phi\mapsto u(t,\cdot\,;\phi)$ to be measurable, and that the Feller property requires its continuity in probability. Theorem~\ref{thm:continuity} only gives the continuity from a space $C_{\theta_{0}}$ to the space $C_{\theta}$, with $\theta_{0}<\theta$. We shall thus obtain both properties on the smaller space, in a form involving the two exponents. To avoid a conflict with the exponent $\theta$ of the state space, we call these two exponents $\rho<\rho'$ in this subsection. In Section~\ref{sec:invariance}, $\rho$ will be the exponent of the tightness, and $\rho'$ the exponent of the test functions. Throughout this subsection, we fix two weight exponents and an integrability exponent such that
\begin{equation}\label{eq:two-exponents}
0<\nrho<\nrhop
\qquad\text{and}\qquad
p\geq p^{*},
\end{equation}
where $p^{*}$ is such that Theorem~\ref{thm:continuity} applies with $(\theta_{0},\theta)$ replaced by $(\rho,\rho')$. Furthermore, we want that both couples $(\nrho,p)$ and $(\nrhop,p)$ satisfy \eqref{assump: theta-p}, so that Theorem~\ref{thm:extend} holds in the spaces $C_{\nrho}$ and $C_{\nrhop}$. This fact is used in the following Lemma~\ref{lem:two-spaces}, which also records how the two spaces $C_{\nrho}$ and $C_{\nrhop}$ are related. Its proof, which rests on the Lusin-Souslin theorem for item \ref{it:ts-borel}, is deferred to Appendix~\ref{app:measurability}.

\begin{lemma}\label{lem:two-spaces}
Let $0<\nrho<\nrhop$, and let $C_{\nrho}$ and $C_{\nrhop}$ be the spaces defined by \eqref{eq:Ctheta}. Then the following holds true.
\begin{enumerate}
\item\label{it:ts-inclusion} The space $C_{\nrho}$ is contained in $C_{\nrhop}$, and $|\phi|_{C_{\nrhop}}\leq 2\,|\phi|_{C_{\nrho}}$ for every $\phi\in C_{\nrho}$.
\item\label{it:ts-borel} The space $C_{\nrho}$ is a Borel subset of $C_{\nrhop}$. Moreover, a subset of $C_{\nrho}$ is a Borel subset of $C_{\nrho}$ if and only if it is a Borel subset of $C_{\nrhop}$.
\item\label{it:ts-solution} Let Assumptions~\ref{assump:noise}, \ref{assump:reac} and \ref{assump:diff} hold, and let $\nrho$, $\nrhop$ and $p$ satisfy~\eqref{eq:two-exponents}. For $\phi\in C_{\nrho}$, let $u(\cdot\,;\phi)$ be the solution given by Theorem~\ref{thm:extend} in the space $C_{\nrhop}$. Then $u(\cdot\,;\phi)$ has paths in $C([0,T];C_{\nrho})$ for every $T>0$, and it coincides with the solution given by Theorem~\ref{thm:extend} in the space $C_{\nrho}$. In particular, for every $t\geq 0$, $u(t,\cdot\,;\phi)$ is a random variable with values in $C_{\nrho}$.
\end{enumerate}
\end{lemma}

We can now state the main result of this subsection. Recall that the spaces $C_{b}(C_{\nrhop})$ and $B_{b}(C_{\nrhop})$ alluded to below are respectively defined by \eqref{eq:normCbE} and Definition~\ref{def: transition}. Since the statement compares functions on $C_{\nrhop}$ with functions on $C_{\nrho}$, we also introduce a notation for restrictions. For a function $F$ defined on $C_{\nrhop}$, we write $F|_{C_{\nrho}}$ for its restriction to $C_{\nrho}$, and for a set $\mathcal{A}$ of functions defined on $C_{\nrhop}$ we set
\begin{equation}\label{eq:restriction}
\mathcal{A}|_{C_{\nrho}}:=\big\{F|_{C_{\nrho}}\,;\,F\in\mathcal{A}\big\} .
\end{equation}
Owing to Lemma~\ref{lem:two-spaces}, these restrictions keep their regularity. Namely, the inclusion of $C_{\nrho}$ in $C_{\nrhop}$ is continuous by item \ref{it:ts-inclusion}, and the Borel subsets of $C_{\nrho}$ are the Borel subsets of $C_{\nrhop}$ contained in $C_{\nrho}$ by item \ref{it:ts-borel}, so that
\begin{equation}\label{eq:restriction-incl}
C_{b}(C_{\nrhop})|_{C_{\nrho}}\subset C_{b}(C_{\nrho}),
\qquad\text{and}\qquad
B_{b}(C_{\nrhop})|_{C_{\nrho}}\subset B_{b}(C_{\nrho}) .
\end{equation}

\begin{proposition}\label{prop:feller-two}
Let Assumptions~\ref{assump:noise}, \ref{assump:reac} and \ref{assump:diff} hold, and let $\nrho$, $\nrhop$ and $p$ satisfy \eqref{eq:two-exponents}. For $\phi\in C_{\nrho}$ and $t\geq 0$, let $u(t,\cdot\,;\phi)$ be the solution given by Theorem~\ref{thm:extend}, which takes its values in $C_{\nrho}$ by item~\ref{it:ts-solution} of Lemma~\ref{lem:two-spaces}. We consider the quantity $\mathcal{P}_{t}F(\phi)$ defined by \eqref{eq: transition}. Then the following holds true.
\begin{enumerate}
\item\label{it:ft-lip} For every $T>0$ there exists a constant $C_{T}$, depending on $T,p,d,\nrhop,\nrho,\eta,\specm,\kappa$ and $L_{\sigma}$, such that for all $\phi_{1},\phi_{2}\in C_{\nrho}$ we have
\begin{equation}\label{eq:feller-lip}
\mathbf{E}\Big[\sup_{t\in[0,T]}|u(t,\cdot\,;\phi_{1})-u(t,\cdot\,;\phi_{2})|^{p}_{C_{\nrhop}}\Big]\leq C_{T}\,|\phi_{1}-\phi_{2}|^{p}_{C_{\nrho}} .
\end{equation}
\item\label{it:ft-feller} For every $t\geq 0$ and every $F\in C_{b}(C_{\nrhop})$, we have $(\mathcal{P}_{t}F)|_{C_{\nrho}}\in C_{b}(C_{\nrho})$. In other words, with the notation \eqref{eq:restriction}, we have
\begin{equation}\label{eq:feller-two}
\mathcal{P}_{t}\big(C_{b}(C_{\nrhop})\big)\big|_{C_{\nrho}}\subset C_{b}(C_{\nrho}),\qquad t\geq 0 .
\end{equation}
\item\label{it:ft-meas} For every $t\geq 0$ and every $F\in B_{b}(C_{\nrhop})$, we have $(\mathcal{P}_{t}F)|_{C_{\nrho}}\in B_{b}(C_{\nrho})$, that is
\begin{equation}\label{eq:feller-two-meas}
\mathcal{P}_{t}\big(B_{b}(C_{\nrhop})\big)\big|_{C_{\nrho}}\subset B_{b}(C_{\nrho}),\qquad t\geq 0 .
\end{equation}
Moreover, for $F\in B_{b}(C_{\nrho})$, relation \eqref{eq: transition} defines $\mathcal{P}_{t}F(\phi)$ for $\phi\in C_{\nrho}$, and the function $\mathcal{P}_{t}F$ thus obtained belongs to $B_{b}(C_{\nrho})$. That is, $\mathcal{P}_{t}$ maps $B_{b}(C_{\nrho})$ to itself. Furthermore, for $t\geq 0$, $\phi\in C_{\nrho}$ and a Borel subset $A$ of $C_{\nrho}$, the quantity
\begin{equation}\label{eq:kernel-def}
P_{t}(\phi,A):=\mathbf{P}\big(u(t,\cdot\,;\phi)\in A\big)
\end{equation}
defines a transition kernel on $C_{\nrho}$. Namely, $A\mapsto P_{t}(\phi,A)$ is a probability measure for every $\phi\in C_{\nrho}$, and $\phi\mapsto P_{t}(\phi,A)$ is Borel on $C_{\nrho}$ for every Borel set $A$. In addition, for $F\in B_{b}(C_{\nrho})$ and $\phi\in C_{\nrho}$ we have
\begin{equation}\label{eq:Pt-kernel}
\mathcal{P}_{t}F(\phi)=\int_{C_{\nrho}}F(\psi)\,P_{t}(\phi,d\psi) .
\end{equation}
\end{enumerate}
\end{proposition}

\begin{proof}
We prove the three items in the order of the statement, since item \ref{it:ft-meas} relies on item \ref{it:ft-feller}.

\smallskip

\noindent
{\it Item \ref{it:ft-lip}.} Let $\phi_{1},\phi_{2}\in C_{\nrho}$. By item \ref{it:ts-inclusion} of Lemma~\ref{lem:two-spaces}, these are deterministic elements of $C_{\nrhop}$, which thus satisfy Assumption~\ref{assump: initial-cond-fixed}, and their difference takes its values in $C_{\nrho}$. Moreover $p\geq p^{*}$ by \eqref{eq:two-exponents}. Hence relation \eqref{eq:continuity-fixed} of Theorem~\ref{thm:continuity} applies with $(\theta_{0},\theta)=(\rho,\rho')$, and it reads as \eqref{eq:feller-lip} for deterministic initial conditions, with $C_{T}:=(2c_{\nrho}C_{T}')^{p}$ where $C_{T}'$ is the constant in \eqref{eq:continuity}.

\smallskip

\noindent
{\it Item \ref{it:ft-feller}.} Let $t\geq 0$ and $F\in C_{b}(C_{\nrhop})$. The function $\mathcal{P}_{t}F$ is bounded by $|F|_{C_{b}(C_{\nrhop})}$, so that we only have to prove its continuity on $C_{\nrho}$. Let $(\phi_{n})_{n\geq 1}$ be a sequence converging to $\phi$ in $C_{\nrho}$. By \eqref{eq:feller-lip} with $T=t$ (the case $t=0$ being trivial), we have
\begin{equation}\label{eq:feller-two-1}
\lim_{n\to\infty}\mathbf{E}\big[|u(t,\cdot\,;\phi_{n})-u(t,\cdot\,;\phi)|^{p}_{C_{\nrhop}}\big]=0 .
\end{equation}
Hence $u(t,\cdot\,;\phi_{n})$ converges to $u(t,\cdot\,;\phi)$ in probability in $C_{\nrhop}$, and every subsequence of $(\phi_{n})_{n\geq 1}$ admits a further subsequence $(\phi_{n_{k}})_{k\geq 1}$ along which $u(t,\cdot\,;\phi_{n_{k}})$ converges to $u(t,\cdot\,;\phi)$ almost surely in $C_{\nrhop}$. Since $F$ is continuous and bounded on $C_{\nrhop}$, the dominated convergence theorem gives
\begin{equation}\label{eq:feller-two-2}
\lim_{k\to\infty}\mathcal{P}_{t}F(\phi_{n_{k}})
=\lim_{k\to\infty}\mathbf{E}\big[F(u(t,\cdot\,;\phi_{n_{k}}))\big]
=\mathbf{E}\big[F(u(t,\cdot\,;\phi))\big]
=\mathcal{P}_{t}F(\phi) .
\end{equation}
The limit in \eqref{eq:feller-two-2} does not depend on the subsequence. Therefore the whole sequence $\mathcal{P}_{t}F(\phi_{n})$ converges to $\mathcal{P}_{t}F(\phi)$, which proves the continuity of $(\mathcal{P}_{t}F)|_{C_{\nrho}}$ on $C_{\nrho}$.

\smallskip

\noindent
{\it Item \ref{it:ft-meas}.} Let $t\geq 0$ be fixed, and let $\mathcal{H}$ be the set of functions $F\in B_{b}(C_{\nrhop})$ such that $\phi\mapsto\mathcal{P}_{t}F(\phi)$ is Borel on $C_{\nrho}$. The set $\mathcal{H}$ is a vector space containing the constant functions. It is closed under bounded monotone limits, by the monotone convergence theorem and since a pointwise limit of Borel functions is Borel. By item \ref{it:ft-feller}, $\mathcal{H}$ contains $C_{b}(C_{\nrhop})$, which is closed under multiplication and generates the Borel $\sigma$-algebra of the metric space $C_{\nrhop}$. By the monotone class theorem, $\mathcal{H}$ is the whole of $B_{b}(C_{\nrhop})$. This is relation \eqref{eq:feller-two-meas}.

Next, let $F\in B_{b}(C_{\nrho})$, and let $\bar{F}$ be the function on $C_{\nrhop}$ which coincides with $F$ on $C_{\nrho}$ and vanishes outside $C_{\nrho}$. By item \ref{it:ts-borel} of Lemma~\ref{lem:two-spaces}, the Borel subsets of $C_{\nrho}$ are Borel in $C_{\nrhop}$, and so is the complement of $C_{\nrho}$. Hence $\bar{F}\in B_{b}(C_{\nrhop})$. Moreover, $u(t,\cdot\,;\phi)$ takes its values in $C_{\nrho}$ for $\phi\in C_{\nrho}$, by item \ref{it:ts-solution} of Lemma~\ref{lem:two-spaces}. Thus
\begin{equation}\label{eq:F-bar}
\mathcal{P}_{t}F(\phi)=\mathbf{E}\big[F(u(t,\cdot\,;\phi))\big]=\mathbf{E}\big[\bar{F}(u(t,\cdot\,;\phi))\big]=\mathcal{P}_{t}\bar{F}(\phi),\qquad\phi\in C_{\nrho},
\end{equation}
which is Borel on $C_{\nrho}$ by \eqref{eq:feller-two-meas}. Taking $F=\mathbf{1}_{A}$ in \eqref{eq:F-bar} for a Borel subset $A$ of $C_{\nrho}$, we get that $\phi\mapsto P_{t}(\phi,A)=\mathcal{P}_{t}\mathbf{1}_{A}(\phi)$ is Borel on $C_{\nrho}$. For a fixed $\phi\in C_{\nrho}$, the set function $A\mapsto P_{t}(\phi,A)$ is the law of the $C_{\nrho}$-valued random variable $u(t,\cdot\,;\phi)$, which is a probability measure on $C_{\nrho}$. Hence $P_{t}$ is a transition kernel on $C_{\nrho}$, and relation \eqref{eq:Pt-kernel} is the expression of the expectation $\mathbf{E}[F(u(t,\cdot\,;\phi))]$ in terms of this law. The proof is complete.
\end{proof}

The transition kernel of Proposition~\ref{prop:feller-two} is what gives a meaning to the measure $\mathcal{P}_{t}^{*}\mu$ of Definition~\ref{def: transition}. We gather here the measurability facts that this requires, and which are used again in Section~\ref{sec:exist}. The proof of the following lemma is deferred to Appendix~\ref{app:measurability}.
\begin{lemma}\label{lem:plug-in}
Let the conditions of Proposition~\ref{prop:feller-two} prevail, and let $\mu$ be a probability measure on $C_{\nrho}$. Then the following holds true.
\begin{enumerate}
\item\label{it:pi-kernel} For every $t\geq 0$, the formula
\begin{equation}\label{eq:Pt-star-kernel}
\mathcal{P}_{t}^{*}\mu(A):=\int_{C_{\nrho}}P_{t}(\phi,A)\,\mu(d\phi), \quad \text{where } A \text{ is a Borel subset of }C_{\nrho},
\end{equation} and where
$P_{t}$ is the transition kernel \eqref{eq:kernel-def}, defines a probability measure on $C_{\nrho}$. Moreover the first identity in \eqref{eq: Pt} holds for every $F\in B_{b}(C_{\nrho})$.
\item\label{it:pi-plug} Let $\xi$ be an initial condition satisfying Assumption~\ref{assump: initial-cond-fixed}, with values in $C_{\nrho}$ and $\mathrm{Law}(\xi)=\mu$, and let $u(\cdot\,;\xi)$ be the solution given by Theorem~\ref{thm:extend}. Then for every $t\geq 0$, $u(t,\cdot\,;\xi)$ is a $C_{\nrho}$-valued random variable and
\begin{equation}\label{eq:plug-in}
\mathrm{Law}\big(u(t,\cdot\,;\xi)\big)=\mathcal{P}_{t}^{*}\mu ,
\qquad\text{that is}\qquad
\mathbf{E}\big[F(u(t,\cdot\,;\xi))\big]=\int_{C_{\nrho}}\mathcal{P}_{t}F(\phi)\,\mu(d\phi),\quad F\in B_{b}(C_{\nrho}).
\end{equation}
In particular, the last identity in \eqref{eq: Pt} holds, and $\mathcal{P}_{t}^{*}\mu$ does not depend on the choice of $\xi$ with law $\mu$.
\end{enumerate}
\end{lemma}

We close this subsection with some comments on Proposition~\ref{prop:feller-two}.

\begin{remark}\label{rmk:feller-two}
Let the conditions of Proposition~\ref{prop:feller-two} prevail.
\begin{enumerate}[label=\textnormal{\textbf{(\alph*)}}]
\item Relation \eqref{eq:feller-two} is weaker than the Feller property on either of the two spaces. On the space $C_{\nrho}$, the Feller property would read $\mathcal{P}_{t}(C_{b}(C_{\nrho}))\subset C_{b}(C_{\nrho})$, and the set $C_{b}(C_{\nrhop})|_{C_{\nrho}}$ of \eqref{eq:restriction-incl} is strictly smaller than $C_{b}(C_{\nrho})$, since the norm $|\cdot|_{C_{\nrho}}$ is not continuous for the topology of $C_{\nrhop}$. On the space $C_{\nrhop}$, the Feller property would require the continuity of $\phi\mapsto u(t,\cdot\,;\phi)$ for the norm $|\cdot|_{C_{\nrhop}}$ on the initial conditions. This is the case $\theta_{0}=\theta$ of Theorem~\ref{thm:continuity}, which is excluded whatever the value of $p$ is, as discussed after the statement of that theorem.
\item Nevertheless, relation \eqref{eq:feller-two} is sufficient to obtain the existence of an invariant measure. Namely, the classical Krylov-Bogoliubov argument tests the limit of the averaged laws against functions of the form $\mathcal{P}_{s}F$. It needs the restriction $(\mathcal{P}_{s}F)|_{C_{\nrho}}$ to be continuous on the space $C_{\nrho}$ where the tightness is established, while $F$ ranges over a class of functions generating the Borel $\sigma$-algebra of that space. This is how \eqref{eq:feller-two} will be used in Section~\ref{sec:invariance}, with $C_{\nrho}$ as the space of the tightness.
\item Owing to item \ref{it:ft-meas} of Proposition~\ref{prop:feller-two}, the semigroup relation $\mathcal{P}_{t+h}=\mathcal{P}_{t}\mathcal{P}_{h}$ of item (i) of Remark~\ref{rmk:markov-property} makes sense on $B_{b}(C_{\nrho})$.
\item For a probability measure $\mu$ on $C_{\nrho}$, the measure $\mathcal{P}_{t}^{*}\mu$ of Definition~\ref{def: transition} is the one given by the kernel formula \eqref{eq:Pt-star-kernel} of Lemma~\ref{lem:plug-in}, and by item \ref{it:pi-plug} of that lemma it is the law of the solution started from any $C_{\nrho}$-valued initial condition with law $\mu$.
\end{enumerate}
\end{remark}

\section{Uniform bounds}\label{sec:bounds}

This section gathers the moment bounds on the solution of \eqref{NSHE} on which the existence of an invariant measure rests. These bounds are uniform in time, and they are the input of the tightness argument of Section~\ref{sec:tightness}. We proceed in four steps. Section~\ref{subsec:prelim} collects two prerequisites, namely the factorization formula of Da Prato and Zabczyk for the stochastic convolution $\mathcal{S}^{\beta}u$ of the penalized equation, and a moment bound on the factorization process $Y_{\alpha}^{\beta}$. In Section~\ref{sec: verification} we apply the Kolmogorov criterion of Proposition~\ref{prop: Salins} to $\mathcal{S}^{\beta}u$, and we obtain a bound which decays as $\beta\to\infty$. In Section~\ref{sec: UniformBoundsU} we derive a uniform-in-time $L^{p}$ bound on $u$. This is achieved through the decomposition $u=w+\mathcal{S}^{\beta}u$, where $w$ solves a random PDE which is controlled pathwise, thanks to the dissipativity Assumption~\ref{assump: dissip}. In Section~\ref{subsec: tight} we upgrade this $L^{p}$ control to a uniform bound in a weighted Besov space (Proposition~\ref{prop: est2}). This is the bound which the compact embedding of Lemma~\ref{lem: compact} turns into tightness in Section~\ref{sec:tightness}.

Throughout this section, we work under the conditions of Theorem~\ref{thm:extend}, and we let $u(t,x)=u(t,x;u_{0})$ denote the solution of $(\ref{NSHE})$ given by that theorem, for an initial condition $u_{0}$ satisfying Assumption~\ref{assump: initial-cond-fixed}. The moment conditions on $u_{0}$ which some of the results below require are stated in those results.

\subsection{Preliminaries}\label{subsec:prelim}
Before we turn to the moment bounds, we need two prerequisites. Namely we recall the factorization formula established by Da Prato and Zabczyk (see their monograph \cite{DZ}), and we collect the moment estimates on the factorization process $Y_{\alpha}^{\beta}$ which the bounds of Section~\ref{sec: verification} require.

\subsubsection{Factorization}
Here we lay out a well known method developed by Da Prato and Zabczyk, which will be useful for us in our preliminary estimates in section \ref{subsec: tight}. Although the proof can be seen as a slight elaboration of \cite[Theorem 5.10]{DZ}, we sketch it in Appendix~\ref{app:factorization} due to the extra exponential term $e^{-\beta(t-s)}$.
\begin{lemma}\label{lem: fact} Under the conditions of Theorem~\ref{thm:extend}, let $\{u(t,x)=u(t,x;u_{0})\}_{(t,x)\in [0,T]\times \R^d}$ denote the unique mild solution of $(\ref{NSHE})$. We let $\beta>0$, $\alpha \in (0,\frac{\eta}{2})$, where $\eta \in (0,1)$ is from the strong Dalang condition $(\ref{eq:strongDalang})$. For $(t,x)\in [0,T]\times \R^d$ we define 
\begin{equation} \label{eq:Sbeta}
\mathcal{S}^{\beta}u(t,x):=\int_{0}^{t}\int_{\R^d}e^{-\beta(t-s)}G(t-s,x-y)\sigma(u(s,y))W(ds,dy).
\end{equation}
Then $\mathcal{S}^{\beta}u$ can also be written as 
\begin{equation}\label{eq:factorization}
    \mathcal{S}^{\beta}u(t,x)= \frac{\sin(\alpha \pi)}{\pi}\int_{0}^{t}\int_{\R^d}e^{-\beta(t-s)}G(t-s,x-y)(t-s)^{\alpha-1}Y_{\alpha}^{\beta}(s,y)ds\,dy,
\end{equation}
where the process $Y_{\alpha}^{\beta}$ is given by
\begin{equation}\label{eq:Yprocess}
    Y_{\alpha}^{\beta}(s,x):=\int_{0}^{s}\int_{\R^d}e^{-\beta(s-r)}G(s-r,x-y)\sigma(u(r,y))(s-r)^{-\alpha}W(dr,dy).
\end{equation}
\end{lemma}

\subsubsection{Tools and moment estimates for $Y_{\alpha}^{\beta}$}
The moment bounds of this section rest on two tools. The first one is the spectral bound on powers of the heat kernel of Lemma~\ref{lem: spectral}, stated in Section~\ref{subsec:conv-tools} and proved in Appendix~\ref{app:conv}. For $\vk\in(1-\eta,1)$, it produces a singularity $\tau^{-\vk}$ at $\tau=0$. The second one is an elementary bound on the time integrals $\int_{0}^{t}e^{-\beta u}\,u^{-\aK}\,du$ in which this singularity ends up, with an explicit decay in $\beta$. We state it as Lemma~\ref{lem:int-exp} below, with a proof deferred to Appendix~\ref{app:factorization}, and then apply it to the random field $Y^{\beta}_{\alpha}$ of \eqref{eq:Yprocess}, for which Section~\ref{subsec: tight} needs a moment bound which is uniform in space and time. Let us also recall that $\sigma$ is bounded by $M_{\sigma}$, according to Assumption~\ref{assump:diff}. All the moment bounds of Sections~\ref{subsec:prelim} and~\ref{sec: verification} are expressed in terms of $M_{\sigma}$.

\begin{lemma}\label{lem:int-exp}
Let $\aK\in[0,1)$, $\ze\in[0,1-\aK)$, $\beta>0$, and $t\ge 0$. Then the following inequality holds true:
\begin{equation}\label{eq:int-exp}
\int_{0}^{t}e^{-\beta u}\,u^{-\aK}\,du
\;\leq\; \frac{\Gamma(1-\aK-\ze)}{\beta^{1-\aK-\ze}}\,t^{\ze}.
\end{equation}
\end{lemma}

As announced above, we now apply Lemma~\ref{lem:int-exp} to establish a uniform bound on the $p$-th moment of the random field $Y^{\beta}_{\alpha}$ appearing in the factorization lemma. The proof combines the BDG inequality with Lemmas~\ref{lem: spectral} and \ref{lem:int-exp}, and it is deferred to Appendix~\ref{app:factorization}.

\begin{lemma}\label{lem:Ybeta-moment}
Under the conditions of Theorem~\ref{thm:extend}, let $\alpha \in(0,\eta/2)$, $\hat{\vk}\in (1-\eta,\,1-2\alpha)$ and $p\geq 1$. Then there exists a constant $C_{p,\specm,\eta,\hat{\vk}}$, independent of $T$, such that the random field $Y^{\beta}_{\alpha}$ defined in \eqref{eq:Yprocess} satisfies
\begin{equation}\label{eq:Ybeta-moment}
\sup_{s\in[0,T]}\sup_{y\in \R^d}\mathbf{E}\left[ |Y_{\alpha}^{\beta}(s,y)|^{p}\right]\leq C_{p,\specm,\eta, \hat{\vk}}\,M_{\sigma}^{p}\left(\frac{1}{\beta^{1-2\alpha}} +\frac{1}{\beta^{1-\hat{\vk}-2\alpha}} \right)^{p/2}.
\end{equation}
\end{lemma}

\subsection{Moment bounds for the stochastic convolution $\mathbf{\mathcal{S}^{\beta}u}$}\label{sec: verification}
In this subsection we bound the moments of the stochastic convolution $\mathcal{S}^{\beta}u$, uniformly in space and time. Our tool is the Kolmogorov-type criterion on unbounded domains due to Salins~\cite{S}, recalled as Proposition~\ref{prop: Salins} in Section~\ref{subsec:conv-tools}. Since $\sigma$ is bounded by $M_{\sigma}$ (see Assumption~\ref{assump:diff}), the bounds below involve $M_{\sigma}$ only, and not the solution $u$. Their important feature is that they decay as $\beta\to\infty$. This is what allows us to close the bound on $u$ in Section~\ref{sec: UniformBoundsU}, by choosing $\beta$ large enough. Let us first recall the definition of $\mathcal{S}^{\beta}u$ and the equation it solves. For a parameter $\beta>0$, the process $\mathcal{S}^{\beta}u$ of $(\ref{eq:Sbeta})$ is the stochastic convolution
\begin{equation}\label{eq: Sbeta-recall}
    \mathcal{S}^{\beta}u(t,x) = \int_{0}^{t}\int_{\R^d}e^{-\beta(t-s)}G(t-s,x-y)\sigma(u(s,y))W(ds,dy).
\end{equation}
One readily checks that $\mathcal{S}^{\beta}u$ is the mild solution of the penalized SPDE
   \begin{equation}\label{pNSHE}
        \begin{cases}
            \partial_{t}\mathcal{S}^{\beta}u(t,x)=\frac{1}{2}\Delta \mathcal{S}^{\beta}u(t,x)-\beta \mathcal{S}^{\beta}u(t,x)+\sigma(u(t,x))\dot{W}(dtdx),\\
            \mathcal{S}^{\beta}u(0,x)=0.
        \end{cases}
    \end{equation}
In Lemmas~\ref{lem: SbetaBound}, \ref{lem: SpaceIncrem} and \ref{lem: TimeIncrem} below, we verify the hypotheses $(\ref{eq: propSalins1})$, $(\ref{eq: propSalins2})$ and $(\ref{eq: propSalins3})$ of Proposition~\ref{prop: Salins} for the random field $X(t,x)=\mathcal{S}^{\beta}u(t,x)$, with the aid of Lemmas~\ref{lem: spectral} and~\ref{lem:int-exp}. The criterion is then applied in Proposition~\ref{prop: SbetaDecay}. The proofs of the three lemmas are deferred to Appendix~\ref{app:factorization}. We begin with the pointwise bound $(\ref{eq: propSalins1})$.
\begin{lemma}\label{lem: SbetaBound} Under the conditions of Theorem~\ref{thm:extend}, let $\mathcal{S}^{\beta}u(t,x)$ be the stochastic convolution defined in $(\ref{eq: Sbeta-recall})$, for $\beta>0$. For any $p\in [1,\infty)$ and $\vk \in (1-\eta,1)$, there exists a constant $C_{p,\specm,\eta,\vk}$, depending only on $p$, the spectral measure $\specm$, the strong Dalang exponent $\eta$ and $\vk$, such that for all $t\geq 0$ and $x\in \R^d$ we have
\begin{equation}\label{eq: prelimLemmaBound1}
\mathbf{E}[|\mathcal{S}^{\beta}u(t,x)|^{p}]\leq C_{p,\specm,\eta,\vk}\,M_{\sigma}^{p}\left(\frac{1}{\beta}+\frac{1}{\beta^{1-\vk}}\right)^{p/2}.
\end{equation}
\end{lemma}

\begin{remark}
Estimates of the type $(\ref{eq: prelimLemmaBound1})$--$(\ref{eq: prelimLemmaBound3})$ are standard in the SPDE literature (see \cite{SS} for instance). We include the proofs here in order to track the explicit decay of the bounds in $\beta$. This will be an essential ingredient for our ergodicity argument.
\end{remark}

Lemma~\ref{lem: SbetaBound} has verified hypothesis $(\ref{eq: propSalins1})$ of Proposition~\ref{prop: Salins}. We now turn to hypothesis~(\ref{eq: propSalins2}), which requires a bound on the spatial increments of $\mathcal{S}^{\beta}u$.

\begin{lemma}\label{lem: SpaceIncrem} Under the conditions of Theorem~\ref{thm:extend}, let $\mathcal{S}^{\beta}u(t,x)$ be the stochastic convolution defined in $(\ref{eq: Sbeta-recall})$, for $\beta>0$. For any $p\in [1,\infty)$, $\lambda\in(0,\eta)$ and $\vk\in(1-\eta,1-\lambda)$, there exists a constant $C_{p,\lambda,d,\specm,\eta,\vk}$, depending only on $p$, $\lambda$, $d$, the spectral measure $\specm$, the strong Dalang exponent $\eta$ and $\vk$, such that for all $t\geq 0$ and $x,y\in \R^d$ we have
    \begin{equation} \label{eq: prelimLemmaBound2}
        \mathbf{E}[|\mathcal{S}^{\beta}u(t,x)-\mathcal{S}^{\beta}u(t,y)|^{p}]
        \leq C_{p,\lambda,d,\specm,\eta,\vk}\,M_{\sigma}^{p} \left(\frac{1}{\beta^{1-\lambda}}+\frac{1}{\beta^{1-\lambda-\vk}}\right)^{p/2}|x-y|^{\lambda p}.
    \end{equation}
\end{lemma}

Lemma~\ref{lem: SpaceIncrem} has verified hypothesis $(\ref{eq: propSalins2})$ of Proposition~\ref{prop: Salins}. We now turn to the time-increment hypothesis $(\ref{eq: propSalins3})$. The proof follows the same BDG strategy as Lemma~\ref{lem: SpaceIncrem}, replacing the spatial gradient bound $(\ref{eq: gradG})$ by the time-derivative bound $(\ref{eq: dtG})$.

\begin{lemma}\label{lem: TimeIncrem}
Under the conditions of Theorem~\ref{thm:extend}, let $\mathcal{S}^{\beta}u(t,x)$ be the stochastic convolution defined in $(\ref{eq: Sbeta-recall})$, for $\beta>0$. For any $p\geq 1$, $\lambda\in(0,\eta/(d+2))$, and $\vk\in(1-\eta,1-\lambda(d+2))$, there exists a constant $C=C_{d,p,\lambda,\specm,\eta,\vk}$ such that
    \begin{equation}\label{eq: prelimLemmaBound3}
        \mathbf{E}[|\mathcal{S}^{\beta}u(t_{2},x)-\mathcal{S}^{\beta}u(t_{1},x)|^{p}]
        \leq C\,M_{\sigma}^{p}
        \left(\frac{1}{\beta^{1-\lambda-\vk}}+\frac{1}{\beta^{1-2\lambda-\vk}}\right)^{p/2}|t_{2}-t_{1}|^{\lambda p/2}
    \end{equation}
    for every $x\in \R^d$ and $t_{1},t_{2}\geq 0$. In particular, $(\ref{eq: prelimLemmaBound3})$ verifies hypothesis $(\ref{eq: propSalins3})$ of Proposition~\ref{prop: Salins}.
\end{lemma}

Lemmas~\ref{lem: SbetaBound}, \ref{lem: SpaceIncrem}, and \ref{lem: TimeIncrem} have verified, respectively, the pointwise bound $(\ref{eq: propSalins1})$, the spatial-increment bound $(\ref{eq: propSalins2})$, and the time-increment bound $(\ref{eq: propSalins3})$ required by Proposition~\ref{prop: Salins} for $X = \mathcal{S}^{\beta}u$. We now apply that Kolmogorov-type result, as well as its space-only version Proposition~\ref{prop: Salins-space}, and obtain two bounds on the moments of $\mathcal{S}^{\beta}u$ in $C_{x_{0},\theta}$. The first one controls the supremum over a time interval of fixed length, and the second one is uniform in time. Both are uniform in the center $x_{0}$ and decay in $\beta$. Let us recall from Proposition~\ref{prop: Salins} that the constant of the Kolmogorov criterion depends on the length of the time interval, which is why the first bound is stated on intervals $[t,t+T]$ and not on $[0,\infty)$.

\begin{proposition}\label{prop: SbetaDecay}
Under the conditions of Theorem~\ref{thm:extend}, let $\mathcal{S}^{\beta}u(t,x)$ be the stochastic convolution defined in $(\ref{eq: Sbeta-recall})$, for $\beta>0$, and let $\theta>0$ and $p\geq p^{*}$.
\begin{enumerate}
\item\label{it:Sbeta-window} For every $T>0$ there exists a constant $C_{p}(\beta,T)$, independent of $x_{0}$, nondecreasing in $T$, with $C_{p}(\beta,T)\to 0$ as $\beta\to \infty$, such that for every $x_{0}\in\R^{d}$ and every $t\geq 0$ we have
    \begin{equation}\label{eq: SigmadecayBound}
        \mathbf{E}\left[\sup_{s\in[t,t+T]}\left|\mathcal{S}^{\beta}u(s,\cdot) \right|_{C_{x_{0},\theta}}^{p} \right]\leq C_{p}(\beta,T)\,M_{\sigma}^{p}.
    \end{equation}
\item\label{it:Sbeta-fixed} There exists a constant $C_{p}(\beta)$, independent of $x_{0}$, with $C_{p}(\beta)\to 0$ as $\beta\to \infty$, such that for every $x_{0}\in\R^{d}$ we have
    \begin{equation}\label{eq: SigmadecayUnif}
        \sup_{t\geq 0}\mathbf{E}\left[\left|\mathcal{S}^{\beta}u(t,\cdot) \right|_{C_{x_{0},\theta}}^{p} \right]\leq C_{p}(\beta)\,M_{\sigma}^{p}.
    \end{equation}
\end{enumerate}
\end{proposition}
 \begin{proof}
We first prove item~\ref{it:Sbeta-window}, and start with the choice of the exponents. We shall use Lemma~\ref{lem: SpaceIncrem} with a space exponent $\lambda_{1}$ and Lemma~\ref{lem: TimeIncrem} with a time exponent $\lambda_{2}$. We fix
  \begin{equation}\label{eq: paramConstraints}
      \lambda_{1}\in(0,\min\{\eta,\theta\}), \qquad
      \lambda_{2}\in\Bigl(0,\frac{\eta}{d+2}\Bigr), \qquad
      \vk\in\bigl(1-\eta,\,1-\max\{\lambda_{1},\lambda_{2}(d+2)\}\bigr),
  \end{equation}
independently of $p$. Since $p\geq p^{*}$, we may assume that $p\lambda_{1}>d+1$, $p\lambda_{2}>2(d+1)$ and $\theta>p\lambda_{1}/(p-(d+1))$, where the last condition holds for $p$ large since $\lambda_{1}<\theta$. With these choices, $\vk$ is admissible in Lemma~\ref{lem: SbetaBound}, the pair $(\lambda_{1},\vk)$ is admissible in Lemma~\ref{lem: SpaceIncrem}, and the pair $(\lambda_{2},\vk)$ is admissible in Lemma~\ref{lem: TimeIncrem}.

Now fix $t\geq 0$ and $T>0$, and consider the random field $X(s,x):=\mathcal{S}^{\beta}u(t+s,x)$ on $[0,T]\times\R^{d}$. The constants of Lemmas~\ref{lem: SbetaBound}, \ref{lem: SpaceIncrem} and~\ref{lem: TimeIncrem} do not depend on time. Hence those three lemmas assert that $X$ satisfies the hypotheses $(\ref{eq: propSalins1})$--$(\ref{eq: propSalins3})$ of Proposition~\ref{prop: Salins}, with the space exponent $\lambda=\lambda_{1}$, the time exponent $\gamma=\lambda_{2}/2$, and
\begin{align}\label{eq: AT-Sbeta}
&A_{T}=A(\beta)  \\
&:=C M_{\sigma}^{p}\left[\left(\frac{1}{\beta}+\frac{1}{\beta^{1-\vk}}\right)^{p/2}+\left(\frac{1}{\beta^{1-\lambda_{1}}}+\frac{1}{\beta^{1-\lambda_{1}-\vk}}\right)^{p/2}\right. \nonumber
\left.\qquad+\left(\frac{1}{\beta^{1-\lambda_{2}-\vk}}+\frac{1}{\beta^{1-2\lambda_{2}-\vk}}\right)^{p/2}\right],
\end{align}
where $C$ depends on $p$, $d$, $\eta$, $\specm$, $\lambda_{1}$, $\lambda_{2}$ and $\vk$, and neither on $t$, $T$ nor $x_{0}$. The conditions \eqref{eq: propSalins4} are satisfied as well, by the three conditions obtained from $p\geq p^{*}$ above. Proposition~\ref{prop: Salins} thus yields
\begin{equation*}
\mathbf{E}\Big[\sup_{s\in[0,T]}|X(s,\cdot)|^{p}_{C_{x_{0},\theta}}\Big]\leq C_{p,d,\lambda_{1},\lambda_{2}/2,\theta,T}\,A(\beta),
\end{equation*}
which is $(\ref{eq: SigmadecayBound})$ with $C_{p}(\beta,T):=C_{p,d,\lambda_{1},\lambda_{2}/2,\theta,T}\,A(\beta)/M_{\sigma}^{p}$. This constant is independent of $t$ and $x_{0}$, and nondecreasing in $T$, as stated in Proposition~\ref{prop: Salins}. It converges to $0$ as $\beta\to\infty$, since all the exponents of $\beta$ in \eqref{eq: AT-Sbeta} are positive by the constraint on $\vk$ in \eqref{eq: paramConstraints}.

We now turn to item~\ref{it:Sbeta-fixed}. Fix $t\geq 0$ and consider the random field $X(x):=\mathcal{S}^{\beta}u(t,x)$ on $\R^{d}$. This time we pick
\begin{equation}\label{eq: paramConstraints-space}
\lambda\in(0,\min\{\eta,\theta\}),\qquad
\vk\in(1-\eta,1-\lambda),
\end{equation}
and since $p\geq p^{*}$ we may assume that $p\lambda>d$ and $\theta>p\lambda/(p-d)$. By Lemmas~\ref{lem: SbetaBound} and~\ref{lem: SpaceIncrem}, the field $X$ satisfies the hypotheses \eqref{eq: propSalins-space1} of Proposition~\ref{prop: Salins-space} with the exponent $\lambda$ and
\begin{equation}\label{eq: A-Sbeta-space}
A=A'(\beta):=C M_{\sigma}^{p}\left[\left(\frac{1}{\beta}+\frac{1}{\beta^{1-\vk}}\right)^{p/2}+\left(\frac{1}{\beta^{1-\lambda}}+\frac{1}{\beta^{1-\lambda-\vk}}\right)^{p/2}\right],
\end{equation}
where $C$ depends on $p$, $d$, $\eta$, $\specm$, $\lambda$ and $\vk$, and neither on $t$ nor $x_{0}$. The conditions \eqref{eq: propSalins-space2} hold by the choice of $\lambda$ and $p^{*}$. Proposition~\ref{prop: Salins-space} thus gives
\begin{equation*}
\mathbf{E}\Big[|\mathcal{S}^{\beta}u(t,\cdot)|^{p}_{C_{x_{0},\theta}}\Big]\leq C_{p,d,\lambda,\theta}\,A'(\beta),
\end{equation*}
with a constant independent of $t$ and $x_{0}$. Taking the supremum over $t\geq 0$, we obtain $(\ref{eq: SigmadecayUnif})$ with $C_{p}(\beta):=C_{p,d,\lambda,\theta}\,A'(\beta)/M_{\sigma}^{p}$, which converges to $0$ as $\beta\to\infty$ since $\vk<1-\lambda$.
 \end{proof}

\subsection{Uniform bounds for $u$}\label{sec: UniformBoundsU}
Proposition~\ref{prop: SbetaDecay} provides a bound on the moments of the stochastic convolution $\mathcal{S}^{\beta}u$ in $C_{x_{0},\theta}$, which is uniform in time and in the center $x_{0}$, and which decays in $\beta$. We now turn to the full solution $u$, and establish a uniform-in-time bound on its moments in $C_{x_{0},\theta}$. As mentioned above, the strategy is to write $u = w + \mathcal{S}^{\beta}u$, where $w := u - \mathcal{S}^{\beta}u$ satisfies a random PDE driven by $f$ and controlled via the subdifferential estimate of Proposition~\ref{prop: subdiffEstimate} in Section~\ref{sec:extend-tools}. Combining both parts then yields the following result.

\begin{proposition}\label{prop: est1}
Under the conditions of Theorem~\ref{thm:extend} and Assumption~\ref{assump: dissip}, let $\theta>0$ and $p\geq p^{*}$, let $\nu>1$ be the exponent of \eqref{assump: growth} and Assumption~\ref{assump: dissip}, and let $\{u(t,x;u_{0}); \, t\ge 0, \, x\in\R^{d}\}$ be the solution of $(\ref{NSHE})$ with initial condition $u_{0}$. We assume that $u_{0}$ satisfies the moment condition \eqref{eq:initial_bound}, that is $\mathbf{E}[|u_{0}|^{p}_{C_{\theta}}]<\infty$. Then there exists a constant $C$, independent of $x_{0}$, such that for every $x_{0}\in\R^{d}$ we have
\begin{equation}\label{est0}
\sup_{t\geq 0}\mathbf{E}\left[|u(t,\cdot)|^{p}_{C_{x_{0},\theta}}\right]\leq C\left(1+\mathbf{E}\left[|u_{0}|^{p}_{C_{x_{0},\theta}}\right]\right).
\end{equation}
Note that the right hand side of \eqref{est0} is finite for every $x_{0}\in\R^{d}$, owing to \eqref{eq:initial_bound} and to the equivalence of the norms $|\cdot|_{C_{x_{0},\theta}}$ over the centers (see Remark~\ref{rem: equiv-norms}).
\end{proposition}
\begin{proof}
Throughout this proof, the value of the constant $C$ may change from line to line, but always depends only on the parameters explicitly indicated. Let us fix $\beta>0$. We shall not need $\beta$ to be large, and the constant $C$ in \eqref{est0} is allowed to depend on it. Recall that $\mathcal{S}^{\beta}u$ is the solution of $(\ref{pNSHE})$, whose mild form is given by
    \begin{equation*}
        \mathcal{S}^{\beta}u(t,x)=\int_{0}^{t}\int_{\R^d}e^{-\beta(t-s)}G(t-s,x-y)\sigma(u(s,y))W(ds,dy).
    \end{equation*}
   Next, we define $w(t,x):=u(t,x)-\mathcal{S}^{\beta}u(t,x)$. Referring to equations \eqref{NSHE} and \eqref{pNSHE}, it can be seen that $w$ solves the following random PDE:
   \begin{equation}\label{eq: w}
   \begin{cases}
       \partial_{t}w(t,x)=\frac{1}{2}\Delta w(t,x)+\beta\,\mathcal{S}^{\beta}u(t,x)+f(u(t,x))\\
       w(0,x)=u_{0}(x).
   \end{cases}
   \end{equation}
Let us stress that the penalization term in \eqref{eq: w} is $\beta\mathcal{S}^{\beta}u$, and not $\beta u$. This is what allows us to control $w$ by the stochastic convolution alone, without any a priori bound on $u$ itself. Our strategy to get \eqref{est0} relies on properties of the process $w$. We divide the proof in several steps. Steps 1 and 2 derive a differential inequality for a weighted supremum norm of $w$, Steps 3 to 5 bound the three terms of its right hand side, Step 6 reduces the resulting inequality to a linear one and solves it by variation of constants, and Step 7 takes expectations.

\smallskip

   \noindent
   \textit{Step 1: An equation for a weighted version of $w$.}
Following in the same spirit of \cite[Theorem 5.4]{S}, we define the following weighted version of $w$:
\begin{equation}\label{eq:tilde-w}
\widetilde{w}(t,x):=\frac{w(t,x)}{(1+|x-x_{0}|^2)^{\theta/2}},
\end{equation}
where the weight $(1+|x-x_{0}|^{2})^{\theta/2}$ is chosen here in place of $(1+|x-x_{0}|^{\theta})$ because the former is differentiable in $x$. This property  will be needed to derive the equation for $\widetilde{w}$ below (see Remark~\ref{rmk: weightEquiv} and the equivalence $(\ref{eq: universalC1})$ for a comparison of the two weights).
Multiplying \eqref{eq: w} by the weight and invoking the conjugation identity \eqref{eq:conjugation} of Lemma~\ref{lem:conjugation}, exactly as in Step~2 of the proof of Lemma~\ref{lem:Muniq} (see \eqref{eq:Muniq-weighted-eq}, with the forcing $\beta\mathcal{S}^{\beta}u+f(u)$ of \eqref{eq: w} in place of $g$ and the differentiation justified as in \cite{S}), we get that $\widetilde{w}$ solves
\begin{equation*}
\partial_{t}\widetilde{w}(t,x)=\mathcal{Q}^{(1)}(t,x)+\mathcal{Q}^{(2)}(t,x)+\mathcal{Q}^{(3)}(t,x),
\end{equation*}
with initial condition
\begin{equation*}
\widetilde{w}(0,x)=\frac{u_{0}(x)}{(1+|x-x_{0}|^{2})^{\theta/2}},
\end{equation*}
and where the terms $\mathcal{Q}^{(1)},\mathcal{Q}^{(2)},\mathcal{Q}^{(3)}$ are respectively defined by
\begin{align}
\mathcal{Q}^{(1)}(t,x)&=\mathcal{L}_{\theta}\widetilde{w}(t,x),\label{eq: Q1} \\
\mathcal{Q}^{(2)}(t,x)&=\frac{f(\mathcal{S}^{\beta}u(t,x)+w(t,x))-f(\mathcal{S}^{\beta}u(t,x))}{(1+|x-x_{0}|^{2})^{\theta/2}},\label{eq: Q2}\\
\mathcal{Q}^{(3)}(t,x)&= \beta\frac{\mathcal{S}^{\beta}u(t,x)}{(1+|x-x_{0}|^{2})^{\theta/2}}+\frac{f(\mathcal{S}^{\beta}u(t,x))}{(1+|x-x_{0}|^{2})^{\theta/2}}.\label{eq: Q3}
\end{align}
In \eqref{eq: Q1}, $\mathcal{L}_{\theta}$ is the differential operator defined in \eqref{eq:Ltheta-def}, with the center $x_{0}$ of the weight.

\smallskip

\noindent
   \textit{Step 2: Applying Proposition $\ref{prop: subdiffEstimate}$.}
   The solution $u$ has paths in $C([0,T];C_{\theta})=C_{0,\theta}([0,T]\times\R^{d})$ by Theorem~\ref{thm:extend} and Lemma~\ref{lem:fc-paths}, and so does $\mathcal{S}^{\beta}u$ by Proposition~\ref{prop:StochConvS} applied to the bounded integrand $e^{-\beta(t-s)}\sigma(u)$ (see Remark~\ref{rmk:conv-bounded}). Since the spaces $C_{x_{0},\theta}([0,T]\times\R^{d})$ coincide as sets for all centers (Remark~\ref{rem: equiv-norms}), $w\in C_{x_{0},\theta}([0,T]\times\R^d)$ for every $T>0$, and the weighted version $\widetilde{w}$ defined by \eqref{eq:tilde-w} is an element of $C_{0}([0,T]\times\R^d)$, and for any $t\geq 0$ there exists at least one point $x_{t}\in \R^d$ where $|\widetilde{w}(t,\cdot)|$ attains its maximum, that is $|\widetilde{w}(t,x_{t})|=|\widetilde{w}(t,\cdot)|_{C_{0}}$. We are thus in the situation of Step~3 of the proof of Lemma~\ref{lem:Muniq}: the measure $\mu_{t}:=\delta_{x_{t}}\,\mathrm{sign}(\widetilde{w}(t,x_{t}))$ belongs to $\partial|\widetilde{w}(t,\cdot)|_{C_{0}}$ by item~\ref{item: subdiffChar} of Proposition~\ref{prop: subdiffEstimate}, and the upper-left-derivative estimate \eqref{eq: derivEst} applied with this $\mu_{t}$ yields, exactly as \eqref{eq:Muniq-dini},
\begin{equation}
\frac{d^{-}}{dt}|\widetilde{w}(t,\cdot)|_{C_{0}}\leq \sum_{i=1}^{3}\mathcal{Q}^{(i)}(t,x_{t})\mathrm{sign}(\widetilde{w}(t,x_{t}))=\sum_{i=1}^{3}{Q}^{(i)}(t,x_{t}),\label{eq: SumQi}
\end{equation}
where each $Q^{(i)}$ is defined by
\begin{equation}\label{g1}
{Q}^{(i)}(t,x_{t}):=\mathcal{Q}^{(i)}(t,x_{t}) \, \mathrm{sign}\lp \tilde{w}(t,x_{t})\rp .
\end{equation}
We will now bound the terms ${Q}^{(i)}$ separately. The bounds rest on the following two elementary relations, valid for every $x\in\R^{d}$ and every function $\phi\in C_{x_{0},\theta/\nu}(\R^{d})$,
\begin{equation}\label{eq:est1-elem}
(1+|x-x_{0}|^{2})^{\theta/2}\leq (1+|x-x_{0}|^{2})^{\theta\nu/2}, \qquad
\frac{|\phi(x)|^{\nu}}{(1+|x-x_{0}|^{2})^{\theta/2}}\leq C_{1}^{\nu}\,|\phi|^{\nu}_{C_{x_{0},\theta/\nu}}.
\end{equation}
The first relation holds since $\nu>1$, and the second one is the equivalence of weights \eqref{eq: universalC1}, written with the exponent $\theta/\nu$ in place of $\theta$ and raised to the power $\nu$.

\smallskip

\noindent
   \textit{Step 3: Bounding ${Q}^{(1)}$.}
   Recall that $\mathcal{Q}^{(1)}=\mathcal{L}_{\theta}\widetilde{w}$ according to \eqref{eq: Q1}, and that $x_{t}$ is a point where $|\widetilde{w}(t,\cdot)|$ attains its maximum. We can thus apply item~\ref{it:conj-max} of Lemma~\ref{lem:conjugation} to $\psi=\widetilde{w}(t,\cdot)$ and $x^{*}=x_{t}$. Owing to \eqref{eq:Ltheta-max} and to the definition \eqref{g1} of $Q^{(1)}$, this yields
\begin{equation}\label{eq:q1}
Q^{(1)}(t,x_{t})\leq |c_{\theta}|_{C_{b}}\,|\widetilde{w}(t,\cdot)|_{C_{0}},
\end{equation}
where $c_{\theta}$ is the bounded coefficient \eqref{eq:Ltheta-coef}. This linear term will be absorbed in Step~6 by the dissipative term coming from Step~4.

\smallskip

\noindent
   \textit{Step 4: Bounding ${Q}^{(2)}$.}
   Recall that $\mathcal{Q}^{(2)}$ is given by~\eqref{eq: Q2} and ${Q}^{(2)}=\mathcal{Q}^{(2)}\, \mathrm{sign}( \tilde{w})$ as specified in~\eqref{g1}.
For bounding this term, we invoke Assumption~\ref{assump: dissip} with $u=\mathcal{S}^{\beta}u(t,x_{t})$ and $h=w(t,x_{t})$, and then the two relations in \eqref{eq:est1-elem}. We obtain
\begin{align}\label{eq:q2}
    Q^{(2)}(t,x_{t})&\leq \frac{-c_{2}|w(t,x_{t})|^\nu+c_{1}(1+|\mathcal{S}^{\beta}u(t,x_{t})|^{\nu})}{(1+|x_{t}-x_{0}|^{2})^{\theta/2}}\nonumber\\
    &\leq -c_{2}\,|\widetilde{w}(t,\cdot)|_{C_0}^{\nu}+C\left(1+|\mathcal{S}^{\beta}u(t,\cdot)|_{C_{x_{0},\theta/\nu}}^{\nu}\right),
\end{align}
where $C$ depends on $c_{1}$, $C_{1}$ and $\nu$. Observe that the first relation in \eqref{eq:est1-elem} has been used to lower bound the dissipative term, thanks to the identity $|\widetilde{w}(t,\cdot)|_{C_{0}}=|w(t,x_{t})|(1+|x_{t}-x_{0}|^{2})^{-\theta/2}$.

\smallskip

\noindent
   \textit{Step 5: Bounding ${Q}^{(3)}$.}
Recall that $\cq^{(3)}$ is defined in~\eqref{eq: Q3} as the sum of the penalization term $\beta\mathcal{S}^{\beta}u$ and the nonlinear drift $f(\mathcal{S}^{\beta}u)$, both weighted and evaluated at the maximizer $x_{t}$, and that ${Q}^{(3)}=\mathcal{Q}^{(3)}\, \mathrm{sign}( \tilde{w})$ according to~\eqref{g1}. The penalization term is bounded by $\beta C_{1}|\mathcal{S}^{\beta}u(t,\cdot)|_{C_{x_{0},\theta}}$, thanks to the equivalence of weights \eqref{eq: universalC1}. The drift term is handled by the growth condition \eqref{assump: growth}, which gives $|f(\mathcal{S}^{\beta}u(t,x_{t}))|\leq K(1+|\mathcal{S}^{\beta}u(t,x_{t})|^{\nu})$, and then by the second relation in \eqref{eq:est1-elem}. We summarize those considerations in the bound below:
\begin{equation}\label{eq:q3}
Q^{(3)}(t,x_{t})
\le
 \beta C_{1}|\mathcal{S}^{\beta}u(t,\cdot)|_{C_{x_{0},\theta}}+C\left(1+|\mathcal{S}^{\beta}u(t,\cdot)|^{\nu}_{C_{x_{0},\theta/\nu}}\right),
\end{equation}
where $C$ depends on $K$, $C_{1}$ and $\nu$.

\smallskip

\noindent
   \textit{Step 6: Linear comparison.}
We now plug the estimates~\eqref{eq:q1}, \eqref{eq:q2} and~\eqref{eq:q3} of Steps~3--5 into \eqref{eq: SumQi}, which gives a single differential inequality for $|\widetilde{w}(t,\cdot)|_{C_{0}}$. To this aim, we set $\varphi(t):=|\widetilde{w}(t,\cdot)|_{C_{0}}$ and we gather the forcing terms of \eqref{eq:q2} and \eqref{eq:q3} into
\begin{equation}\label{eq: forcing}
F(t):=C\left(1+|\mathcal{S}^{\beta}u(t,\cdot)|^{\nu}_{C_{x_{0},\theta/\nu}}+\beta\,|\mathcal{S}^{\beta}u(t,\cdot)|_{C_{x_{0},\theta}}\right),
\end{equation}
where $C$ depends on $c_{1}$, $K$, $C_{1}$ and $\nu$. Then \eqref{eq: SumQi} reads
\begin{equation}\label{eq: diffineq-raw}
\frac{d^{-}}{dt}\varphi(t)\leq |c_{\theta}|_{C_{b}}\,\varphi(t)-c_{2}\,\varphi(t)^{\nu}+F(t).
\end{equation}
Since $\nu>1$ by \eqref{assump: growth}, Young's inequality yields $(1+|c_{\theta}|_{C_{b}})\,\varphi\leq c_{2}\varphi^{\nu}+C$, with a constant $C$ depending on $|c_{\theta}|_{C_{b}}$, $c_{2}$ and $\nu$. Enlarging the constant in \eqref{eq: forcing} accordingly, we end up with the linear differential inequality
\begin{equation}\label{eq: diffineq}
\frac{d^{-}}{dt}\varphi(t)\leq -\varphi(t)+F(t), \qquad t\geq 0.
\end{equation}
Let us stress two features of \eqref{eq: diffineq}. First, the forcing term $F$ only involves the stochastic convolution $\mathcal{S}^{\beta}u$, and not the solution $u$. Second, the dissipativity of $f$ has been used to produce the negative linear term $-\varphi$, and this is the only place where it enters. We now apply the comparison theorem~\cite[Theorem 1.3]{T} to \eqref{eq: diffineq} and to the linear equation $\dot{z}=-z+F$ with $z(0)=\varphi(0)$. By the variation of constants formula, this yields
\begin{equation}\label{eq: W-est}
\varphi(t)\leq \varphi(0)\,e^{-t}+\int_{0}^{t}e^{-(t-s)}F(s)\,ds, \qquad t\geq 0,
\end{equation}
where $\varphi(0)=|\widetilde{w}(0,\cdot)|_{C_{0}}\leq C_{1}|u_{0}|_{C_{x_{0},\theta}}$ by \eqref{eq:tilde-w} and \eqref{eq: universalC1}.

\smallskip

\noindent
   \textit{Step 7: Conclusion.}
   Let us summarize our considerations so far. Invoking the decomposition $u=w+\mathcal{S}^{\beta}u$, then the equivalence of weights \eqref{eq: universalC1} and the definition \eqref{eq:tilde-w} of $\widetilde{w}$, we have for every $t\geq 0$
  \begin{equation*}
 |u(t,\cdot)|_{C_{x_{0},\theta}}
 \leq 
 C_{1}\,|\widetilde{w}(t,\cdot)|_{C_{0}}+|\mathcal{S}^{\beta}u(t,\cdot)|_{C_{x_{0},\theta}} .
\end{equation*}
Plugging \eqref{eq: W-est} in this relation, raising to the power $p$ and taking expectations, we obtain
\begin{eqnarray}
    \mathbf{E}\left[|u(t,\cdot)|_{C_{x_{0},\theta}}^{p}\right]
    &\leq& C\left(\mathbf{E}\left[|u_{0}|^{p}_{C_{x_{0},\theta}}\right]+\mathbf{E}\left[\left(\int_{0}^{t}e^{-(t-s)}F(s)\,ds\right)^{p}\right]\right. \nonumber\\
    &&\left.\qquad+\mathbf{E}\left[|\mathcal{S}^{\beta}u(t,\cdot)|^{p}_{C_{x_{0},\theta}}\right]\right),
    \label{eq: NormuBound1}
\end{eqnarray}
where $C$ depends on $p$ and $C_{1}$. The measure $e^{-(t-s)}ds$ on $[0,t]$ has total mass at most $1$. Hence Jensen's inequality gives
\begin{equation}\label{eq: jensen-F}
\mathbf{E}\left[\left(\int_{0}^{t}e^{-(t-s)}F(s)\,ds\right)^{p}\right]\leq \int_{0}^{t}e^{-(t-s)}\,\mathbf{E}\left[F(s)^{p}\right]ds\leq \sup_{s\geq 0}\mathbf{E}\left[F(s)^{p}\right].
\end{equation}
Next we go back to the expression \eqref{eq: forcing} of $F$. Raising it to the power $p$, we get
\begin{eqnarray}
\lefteqn{\sup_{s\geq 0}\mathbf{E}\left[F(s)^{p}\right]+\mathbf{E}\left[|\mathcal{S}^{\beta}u(t,\cdot)|^{p}_{C_{x_{0},\theta}}\right]} \nonumber\\
&\leq& C\left(1+\sup_{s\geq 0}\mathbf{E}\left[|\mathcal{S}^{\beta}u(s,\cdot)|^{p\nu}_{C_{x_{0},\theta/\nu}}\right]+\left(1+\beta^{p}\right)\sup_{s\geq 0}\mathbf{E}\left[|\mathcal{S}^{\beta}u(s,\cdot)|^{p}_{C_{x_{0},\theta}}\right]\right),
\label{eq: F-moment}
\end{eqnarray}
where $C$ depends on $p$ and on the constant in \eqref{eq: forcing}. The two suprema in \eqref{eq: F-moment} are handled by item~\ref{it:Sbeta-fixed} of Proposition~\ref{prop: SbetaDecay}, applied with the pairs $(\theta,p)$ and $(\theta/\nu,p\nu)$ of weight and integrability exponents. Both pairs are admissible there since $p\nu\geq p\geq p^{*}$, where $p^{*}$ is allowed to depend on $\theta$ and $\nu$. Invoking \eqref{eq: SigmadecayUnif} twice, the right hand side of \eqref{eq: F-moment} is thus bounded by $C(1+C_{p\nu}'(\beta)M_{\sigma}^{p\nu}+(1+\beta^{p})C_{p}(\beta)M_{\sigma}^{p})$, where $C_{p}(\beta)$ and $C_{p\nu}'(\beta)$ are the constants of \eqref{eq: SigmadecayUnif} for the pairs $(\theta,p)$ and $(\theta/\nu,p\nu)$ respectively, both independent of $x_{0}$. Plugging this bound and \eqref{eq: jensen-F} into \eqref{eq: NormuBound1}, we obtain
\begin{equation*}
    \mathbf{E}\left[|u(t,\cdot)|_{C_{x_{0},\theta}}^{p}\right]
    \leq C\left(1+\mathbf{E}\left[|u_{0}|^{p}_{C_{x_{0},\theta}}\right]\right),
\end{equation*}
where $C$ depends on $\beta$, $p$, $\nu$, $d$, $\eta$, $\theta$, $M_{\sigma}$ and on the constants $c_{1}$, $c_{2}$, $K$, $C_{1}$, $|c_{\theta}|_{C_{b}}$, but neither on $t$ nor on $x_{0}$. Taking the supremum over $t\geq 0$, we obtain $(\ref{est0})$.
\end{proof}

\subsection{A weighted Besov bound}\label{subsec: tight}
Armed with the moment bound of Proposition~\ref{prop: est1}, we now establish the Besov bound on which the tightness argument of Section~\ref{sec:tightness} rests. The key step is Proposition~\ref{prop: est2}, which upgrades the $L^{p}$ control on $u$ to a uniform bound in a weighted Besov--H\"{o}lder space.

Throughout this subsection, we work under the conditions of Theorem~\ref{thm:extend} and the dissipativity Assumption~\ref{assump: dissip}, and we let $u(t,x)=u(t,x;u_{0})$ denote the solution of $(\ref{NSHE})$, for an initial condition $u_{0}$ satisfying the moment condition \eqref{eq:initial-pnu} below.

\begin{proposition}\label{prop: est2}
Under the conditions of Theorem~\ref{thm:extend} and Assumption~\ref{assump: dissip}, let $\theta>0$, and let $\nu$ be the exponent from the growth condition $(\ref{assump: growth})$ of the drift term $f$. Let $p\geq p^{*}$, and let $\tilde{s}\in(0,\eta-2/p)$ and $\tilde{\theta}>\nu\theta+d/p$ be two exponents. Furthermore, assume that the initial condition satisfies
\begin{equation}\label{eq:initial-pnu}
\mathbf{E}\big[|u_{0}|^{p\nu}_{C_{\theta}}\big]<\infty .
\end{equation}
Then, for every $x_{0}\in\R^{d}$ and for any $T^{*}>0$, the solution $\{u(t,x;u_{0});\, t\ge 0,\, x\in\R^{d}\}$ of $(\ref{NSHE})$ satisfies
    \begin{equation}\label{est2}
       C_{T^{*}}:=\sup_{t \ge T^*} \mathbf{E}\left[ |u(t,\cdot)|_{\mathcal{B}^{\tilde{s}}_{p,p}(w_{x_{0}}^{\tilde{\theta}})}^{p} \right] < \infty,
    \end{equation}
    where $\mathcal{B}^{\tilde{s}}_{p,p}(w_{x_{0}}^{\tilde{\theta}})$ is the Besov space of Definition~\ref{def: Besov}.
\end{proposition}

\begin{remark}
Proposition~\ref{prop: est2} requires the moment of order $p\nu$ of the initial condition in \eqref{eq:initial-pnu}, whereas Proposition~\ref{prop: est1} only needs the moment of order $p$ in \eqref{eq:initial_bound}. The additional factor of $\nu$ arises from the nonlinear drift term $f$: in the decomposition \eqref{eq:decomp-2}, the term $D$ involves $f(u)$. This term satisfies the growth bound $|f(u)|\leq K(1+|u|^{\nu})$ by \eqref{assump: growth}. Controlling the Besov seminorm of $D$ therefore requires $\mathbf{E}[|u(t,\cdot)|^{p\nu}_{C_{x_{0},\theta}}]$ to be finite, which in turn forces \eqref{eq:initial-pnu}. For the same reason, the bound on $D$ only holds with the weight $(w^{\theta}_{x_{0}})^{\nu}$, see \eqref{eq: D-pathwise}. Integrating it against the weight $w^{\tilde\theta}_{x_{0}}$ of the Besov norm in \eqref{eq: uBesovseminorm} requires $\tilde{\theta}>\nu\theta+d/p$, which is why this condition appears in the statement.
\end{remark}

\begin{proof}[Proof of Proposition~\ref{prop: est2}]
We fix a parameter $\alpha\in(1/p,\eta/2)$ for the factorization identity of Lemma~\ref{lem: fact} such that $\tilde{s}<2\alpha-2/p$, which is possible since $\tilde{s}<\eta-2/p$.
 Recall that we decomposed the solution of \eqref{NSHE} as
 \begin{equation}\label{eq: recall-decomp}
 u(t,x)=w(t,x)+\mathcal{S}^{\beta}u(t,x),
 \end{equation}
 where $w(t,x)$ is a mild solution of $(\ref{eq: w})$ and $\mathcal{S}^{\beta}u$ is defined by~\eqref{eq: Sbeta-recall}. 
 We further decompose $w$ into
 \begin{equation}\label{eq:decomp-2}
w(t,x)=I(t,x)+D(t,x),
\quad\text{so that}\quad
u(t,x)=I(t,x)+D(t,x)+\mathcal{S}^{\beta}u(t,x) \, ,
 \end{equation}
 where we define
    \begin{align}\label{eq: I}
I(t,x)&:=\int_{\R^d}e^{-\beta t}G(t,x-z)u_{0}(z)\,dz\\ \label{eq: D}
D(t,x)&:=\int_{0}^{t}\int_{\R^d}e^{-\beta(t-s)}G(t-s,x-z)\{\beta u(s,z)+f(u(s,z))\}\,dz\,ds.
    \end{align}
   In the remainder of the proof we bound the $\mathcal{B}^{\tilde{s}}_{p,p}(w_{x_{0}}^{\tilde{\theta}})$-norm of Definition~\ref{def: Besov} for each term of the decomposition \eqref{eq:decomp-2}. To this aim, in Steps~1 and~2 we will  directly estimate the spatial increments $\delta_{h}I$ and $\delta_{h}D$, where for $h\in \R^d$, $\delta_{h}$ is the finite difference operator acting on space-time functions $\varphi=\varphi(t,x)$, defined by
   \[\delta_{h}\varphi(t,x):=\varphi(t,x+h)-\varphi(t,x).\]
   Finally in Step~3 we bound $\mathcal{S}^{\beta}u$ by the factorization method, Schauder estimates, Besov embeddings and the moment bounds of $Y_{\alpha}^{\beta}$ from Lemma~\ref{lem:Ybeta-moment}. Step~4 combines these bounds.
   
   \noindent
 \textit{Step 1: Bounding $\delta_{h}I(r,x)$.}
Let $t>0$ and $T^{*}\leq t$. We estimate the $p$-th moment of $\delta_{h}I(r,x)$, following the same strategy as in the proof of Lemma~\ref{lem: SpaceIncrem}.
Namely we first apply the $\lambda/(1-\lambda)$ decomposition~\eqref{eq: decompH} from the proof of Lemma~\ref{lem: SpaceIncrem} to the heat kernel difference in $\delta_{h}I(r,x)$. Using the sub-additivity of $t\mapsto t^{1-\lambda}$, we have
\begin{equation*}
|G(r,x+h-z)-G(r,x-z)| \leq |G(r,x+h-z)-G(r,x-z)|^{\lambda}
\bigl\{G(r,x+h-z)^{1-\lambda}+G(r,x-z)^{1-\lambda}\bigr\}.
\end{equation*}
Employing this bound and writing $|u_{0}(z)|=|u_{0}(z)|\,w_{x_{0}}^{\theta}(z)\,w_{x_{0}}^{\theta}(z)^{-1}$, we obtain
\begin{multline*}
    \sup_{r\in[T^{*},t]}\mathbf{E}\!\left[|\delta_{h}I(r,\cdot)(x)|^{p}\right]
    \leq \sup_{r\in[T^{*},t]}\mathbf{E}\Bigg[\Bigg(
        \int_{\mathbb{R}^d} e^{-\beta r}
        \left|G(r,x{+}h{-}z)-G(r,x{-}z)\right|^{\lambda}
        \\
\times
        \bigl\{G(r,x{+}h{-}z)^{1-\lambda}+G(r,x{-}z)^{1-\lambda}\bigr\}
        \,|u_{0}(z)|\,w_{x_{0}}^{\theta}(z)\,w_{x_{0}}^{\theta}(z)^{-1}\,dz
    \Bigg)^{\!p}\,\Bigg].
\end{multline*}
Next, we apply the mean value theorem to the heat kernel difference and invoke the gradient bound~\eqref{eq: gradG}, exactly as in~\eqref{eq: MVTforH}, to obtain
\begin{equation*}
|G(r,x+h-z)-G(r,x-z)|^{\lambda} \leq C_{\lambda}\,r^{-\frac{\lambda(d+1)}{2}}\,|h|^{\lambda}.
\end{equation*}
Substituting this bound, bounding $|u_{0}(z)|w_{x_{0}}^{\theta}(z)$ by $|u_{0}|_{C_{x_{0},\theta}}$ and taking expectations, we arrive at
\begin{align}\label{eq: supI-prelim}
    \sup_{r\in[T^{*},t]}\mathbf{E}\!\left[|\delta_{h}I(r,\cdot)(x)|^{p}\right]
    &\lesssim_{\lambda}\;
        \mathbf{E}\big[|u_{0}|^{p}_{C_{x_{0},\theta}}\big]
        \sup_{r\in[T^{*},t]}\Bigg(
        \int_{\mathbb{R}^d} e^{-\beta r}\,r^{-\frac{\lambda(d+1)}{2}}\,|h|^{\lambda}
        \nonumber\\
    &\hspace{1.0cm}\times
        \bigl\{G(r,x{+}h{-}z)^{1-\lambda}+G(r,x{-}z)^{1-\lambda}\bigr\}
        \,w_{x_{0}}^{\theta}(z)^{-1}\,dz
    \Bigg)^{\!p}.
\end{align}
It remains to bound the spatial integral in~\eqref{eq: supI-prelim}. Writing $G(r,y)^{1-\lambda}=C_{\lambda,d}\,r^{d\lambda/2}G(r/(1-\lambda),y)$ and applying the weight inequality \eqref{eq:weight-triangle} to $w_{x_{0}}^{\theta}(z)^{-1}=1+|z-x_{0}|^{\theta}$, we get
\begin{equation}\label{eq: Gauss-split}
\int_{\mathbb{R}^d} G(r,x\pm h-z)^{1-\lambda}\,w_{x_{0}}^{\theta}(z)^{-1}\,dz
\;\leq\; C_{\lambda,d}\,r^{d\lambda/2}\left(w_{x_{0}}^{\theta}(x)^{-1}+\int_{\mathbb{R}^d}G\!\left(\tfrac{r}{1-\lambda},y\right)|y|^{\theta}\,dy\right).
\end{equation}
Using the Gaussian moment $\int_{\R^{d}}G(\tau,y)|y|^{\theta}\,dy=m_{\theta}\tau^{\theta/2}$ of \eqref{eq:G-moment} with $\tau=r/(1-\lambda)$ in~\eqref{eq: Gauss-split}, and absorbing $(1-\lambda)^{-\theta/2}$ into the implicit constant, we conclude
\begin{equation}\label{eq: Gauss-wt-conv}
\int_{\mathbb{R}^d} G(r,x\pm h-z)^{1-\lambda}\,w_{x_{0}}^{\theta}(z)^{-1}\,dz
\lesssim_{\lambda,\theta,d} r^{d\lambda/2}(1+r^{\theta/2})\,w_{x_{0}}^{\theta}(x)^{-1}.
\end{equation}
Inserting~\eqref{eq: Gauss-wt-conv} into~\eqref{eq: supI-prelim}, the spatial integral contributes the factor $r^{-\lambda(d+1)/2}\cdot r^{d\lambda/2}(1+r^{\theta/2})=r^{-\lambda/2}(1+r^{\theta/2})$, which grows at most polynomially in $r$. Since the exponential $e^{-\beta r}$ in~\eqref{eq: supI-prelim} dominates any polynomial, we have $\sup_{r\geq T^{*}}e^{-\beta r}\,r^{-\lambda/2}(1+r^{\theta/2})\leq C_{\beta,\lambda,\theta}<\infty$, and we conclude
\begin{equation}\label{eq: deltaI}
\mathbf{E}\!\left[\sup_{r\in[T^{*},t]}|\delta_{h}I(r,\cdot)(x)|^{p}\right]
\lesssim_{\lambda,p,T^{*}}
\mathbf{E}\big[|u_{0}|^{p}_{C_{x_{0},\theta}}\big]\,|h|^{\lambda p}\,w_{x_{0}}^{\theta}(x)^{-p}.
\end{equation}

\smallskip
\noindent
\textit{Step 2: Bounding $\delta_{h}D(r,x)$.}
We proceed pathwise. The argument parallels Step~1, so we only detail the two new ingredients: the growth condition~\eqref{assump: growth} on $f$, and the uniform-in-time bound~\eqref{est0}.
Specifically, applying the $\lambda/(1-\lambda)$ decomposition~\eqref{eq: decompH} and the gradient bound~\eqref{eq: MVTforH} to the heat kernel difference in $\delta_{h}D(r,x)$, exactly as in Step~1, we obtain
\begin{multline*}
|\delta_{h}D(r,x)|^{p}
\leq\Bigg(\int_{0}^{r}\int_{\R^d}e^{-\beta(r-s)}(r-s)^{-\frac{\lambda(d+1)}{2}}|h|^{\lambda}\\
\times\bigl\{G(r-s,x{+}h{-}z)^{1-\lambda}+G(r-s,x{-}z)^{1-\lambda}\bigr\}
\bigl\{\beta|u(s,z)|+|f(u(s,z))|\bigr\}\,dz\,ds \Bigg)^{p}.
\end{multline*}
We now bound the factor $\beta|u(s,z)|+|f(u(s,z))|$. The definition of $|\cdot|_{C_{x_{0},\theta}}$ gives $|u(s,z)|\leq|u(s,\cdot)|_{C_{x_{0},\theta}}\,w_{x_{0}}^{\theta}(z)^{-1}$. The growth condition~\eqref{assump: growth} on $f$ then yields $|f(u(s,z))|\leq K(1+|u(s,\cdot)|_{C_{x_{0},\theta}}^{\nu}\,w_{x_{0}}^{\theta}(z)^{-\nu})$. Absorbing the linear term $\beta|u(s,z)|$ by Young's inequality, we get
\begin{equation*}
\beta|u(s,z)|+|f(u(s,z))| \lesssim_{\beta}\, 1 + |u(s,\cdot)|_{C_{x_{0},\theta}}^{\nu}\,w_{x_{0}}^{\theta}(z)^{-\nu}.
\end{equation*}
Substituting this bound into the previous estimate, we arrive at
\begin{multline*}
|\delta_{h}D(r,\cdot)(x)|^{p}
\lesssim_{\beta} \Bigg(\int_{0}^{r}\int_{\R^d}e^{-\beta(r-s)}(r-s)^{-\frac{\lambda(d+1)}{2}}|h|^{\lambda}\\
\times\bigl\{G(r-s,x{+}h{-}z)^{1-\lambda}+G(r-s,x{-}z)^{1-\lambda}\bigr\}
\bigl(1+|u(s,\cdot)|_{C_{x_{0},\theta}}^{\nu}\bigr)\,w_{x_{0}}^{\theta}(z)^{-\nu}\,dz\,ds \Bigg)^{p}.
\end{multline*}
Next, we set
\[H:=\int_{0}^{r}\int_{\R^d}e^{-\beta(r-s)}(r-s)^{-\frac{\lambda(d+1)}{2}}\bigl\{G(r-s,x{+}h{-}z)^{1-\lambda}+G(r-s,x{-}z)^{1-\lambda}\bigr\}\,w_{x_{0}}^{\theta}(z)^{-\nu}\,dz\,ds,\]
and apply Jensen's inequality alongside the uniform bound \eqref{est0} with $p$ replaced by $p\nu$ (valid since $u_{0}$ satisfies \eqref{eq:initial-pnu} by assumption) we arrive at
\[\mathbf{E}[|\delta_{h}D(r,\cdot)(x)|^{p}]\lesssim_{\beta,p} H^{p}\big(1+\mathbf{E}\big[|u_{0}|^{p\nu}_{C_{x_{0},\theta}}\big]\big)|h|^{\lambda p}.\]
It remains to bound $H^{p}$. The argument follows the same lines as~\eqref{eq: Gauss-wt-conv}: writing $G(r-s,y)^{1-\lambda}=C_{\lambda,d}(r-s)^{d\lambda/2}G((r-s)/(1-\lambda),y)$ and applying the weight inequality \eqref{eq:weight-triangle} with the exponent $\theta\nu$ to 
\begin{equation*}
w_{x_{0}}^{\theta}(z)^{-\nu}=\left(1+|z-x_{0}|^{\theta}\right)^{\nu}\leq C_{\nu}\left(1+|z-x_{0}|^{\theta\nu}\right) ,
\end{equation*}
we obtain
\begin{multline}\label{eq: Gauss-split-nu}
\int_{\R^d} G(r-s,x\pm h-z)^{1-\lambda}\,w_{x_{0}}^{\theta}(z)^{-\nu}\,dz \\
\;\leq\; C_{\lambda,\nu,d}\,(r-s)^{d\lambda/2}\left(w_{x_{0}}^{\theta}(x)^{-\nu}+\int_{\R^d}G\!\left(\tfrac{r-s}{1-\lambda},y\right)|y|^{\theta\nu}\,dy\right).
\end{multline}
Using the Gaussian moment $\int_{\R^{d}}G(\tau,y)|y|^{\theta\nu}\,dy=m_{\theta\nu}\tau^{\theta\nu/2}$ of \eqref{eq:G-moment} with $\tau=(r-s)/(1-\lambda)$ in~\eqref{eq: Gauss-split-nu}, we thus get the analogue of~\eqref{eq: Gauss-wt-conv}:
\begin{equation}\label{eq: Gauss-wt-conv-nu}
\int_{\R^d} G(r-s,x\pm h-z)^{1-\lambda}\,w_{x_{0}}^{\theta}(z)^{-\nu}\,dz
\;\lesssim_{\lambda,\theta,\nu,d}\;(r-s)^{d\lambda/2}\bigl(1+(r-s)^{\theta\nu/2}\bigr)\,w_{x_{0}}^{\theta}(x)^{-\nu}.
\end{equation}
After inserting~\eqref{eq: Gauss-wt-conv-nu} into $H$, the product $(r-s)^{-\lambda(d+1)/2}\cdot(r-s)^{d\lambda/2}(1+(r-s)^{\theta\nu/2})=(r-s)^{-\lambda/2}(1+(r-s)^{\theta\nu/2})$ grows at most polynomially in $r-s$ and is absorbed by the exponential $e^{-\beta(r-s)}$, giving
\begin{equation}\label{eq: time-int-nu}
\int_{0}^{r}e^{-\beta(r-s)}(r-s)^{-\lambda/2}\bigl(1+(r-s)^{\theta\nu/2}\bigr)\,ds
\;\leq\;
C_{\beta,\lambda,\theta,\nu}.
\end{equation}
Therefore inserting~\eqref{eq: Gauss-wt-conv-nu} and~\eqref{eq: time-int-nu} into $H$ yields the bound
\begin{equation}\label{eq: D-pathwise}
\sup_{r\in[0,t]}\mathbf{E}[|\delta_{h}D(r,x)|^{p}]
\;\lesssim_{\beta,p}\;
C_{\beta,\lambda,\theta,\nu}\,|h|^{\lambda p}\left(1+\mathbf{E}\big[|u_{0}|^{p\nu}_{C_{x_{0},\theta}}\big]\right)\,w_{x_{0}}^{\theta}(x)^{-\nu p}.
\end{equation}
\smallskip

\noindent
\textit{Step 3: Bounding $\mathcal{S}^{\beta}u(t,\cdot)$ in Besov norms.} The proof of this step is inspired by the methodology from Proposition 4.5 and Remark 4.6 in~\cite{C}. The factorization method stated in Lemma~\ref{lem: fact} will be applied at this point. Together with Minkowski's inequality, it enables to write
\begin{equation*}
|\mathcal{S}^{\beta}u(t,\cdot)|_{\mathcal{B}^{\tilde{s}}_{p,p}(w_{x_{0}}^{\tilde\theta})}
    \leq c \int_{0}^{t}e^{-\beta(t-r)}(t-r)^{\alpha-1}|H_{t-r}Y_{\alpha}^{\beta}(r,\cdot)|_{\mathcal{B}^{\tilde{s}}_{p,p}(w_{x_{0}}^{\tilde\theta})}dr .
    \end{equation*}
Fix $\varepsilon\in(0,2\alpha-\frac{2}{p}-\tilde{s})$, considered as a small parameter. We now apply a Schauder estimate taken from \cite[\textsection 3.2]{MW} as well as a Besov embedding borrowed from \cite[\textsection 2.3.3, \textsection 4.3.2]{ET}:
\begin{equation}\label{schauder1}
\|H_{t}f\|_{\mathcal{B}^{\tilde{s}}_{p,p}(w_{x_{0}}^{\tilde\theta})}\lesssim t^{-(\tilde{s}+\varepsilon)/2}\|f\|_{\mathcal{B}_{p,p}^{-\varepsilon}(w_{x_{0}}^{\tilde{\theta}})} ,
\quad \text{and}\quad
\|f\|_{\mathcal{B}_{p,p}^{-\varepsilon}(w_{x_{0}}^{\tilde\theta})}\lesssim \|f\|_{L^{p}(w_{x_{0}}^{\tilde\theta})},
\end{equation}
to end up with:
\begin{align}\label{eq: III3}
    & |\mathcal{S}^{\beta}u(t,\cdot)|_{\mathcal{B}^{\tilde{s}}_{p,p}(w_{x_{0}}^{\tilde\theta})}
    \leq c\int_{0}^{t}e^{-\beta(t-r)}(t-r)^{\alpha-(\tilde{s}+\varepsilon)/2-1}|Y_{\alpha}^{\beta}(r,\cdot)|_{L^{p}\left(w_{x_{0}}^{\tilde\theta}\right)}dr\nonumber\\
    &\leq c \lp\int_{0}^{t}e^{-\frac{\beta p }{2(p-1)}(t-r)}(t-r)^{(\alpha-(\tilde{s}+\varepsilon)/2-1)\frac{p}{p-1}}dr \rp^{\frac{p-1}{p}}
    \lp\int_{0}^{t}e^{-\frac{\beta p}{2}(t-r)}|Y_{\alpha}^{\beta}(r,\cdot)|^{p}_{L^{p}\left(w_{x_{0}}^{\tilde\theta}\right)}dr \rp^{\frac{1}{p}},
    \end{align}
where we have used H\"older's inequality for the second upper bound. Note that the first integral in \eqref{eq: III3} is a constant with no time-dependence thanks to \eqref{eq:int-exp} and the fact that $\tilde{s}+\varepsilon\in(0,2\alpha-2/p)$.
Next, by Fubini and making use of the moment estimate of $Y_{\alpha}^{\beta}$ from Lemma \ref{lem:Ybeta-moment} for the second integral in \eqref{eq: III3},  we arrive at
\begin{align}\label{h2}
\sup_{t\geq 0}\mathbf{E}[|\mathcal{S}^{\beta}u(t,\cdot)|_{\mathcal{B}^{\tilde{s}}_{p,p}(w_{x_{0}}^{\tilde\theta})}^{p}]&\lesssim \sup_{t\geq0}\int_{0}^{t}\int_{\R^d}e^{-\frac{\beta p}{2}(t-r)}w_{x_{0}}^{\tilde\theta}(x)^{p}\mathbf{E}[|Y_{\alpha}^{\beta}(r,x)|^{p}]dx\,dr\\
    &\lesssim M_{\sigma}^{p}\sup_{t\geq 0}\int_{0}^{t}\int_{\R^d}e^{-\frac{\beta p}{2}(t-r)}w_{x_{0}}^{\tilde\theta}(x)^{p}\,dx\,dr. \notag
\end{align}
Then, since $\int_{\R^{d}}w_{x_{0}}^{\tilde\theta}(x)^{p}\,dx<\infty$ owing to $\tilde{\theta}p>d$, we obtain
\begin{equation}\label{eq:Sbeta-Besov}
\sup_{t\geq 0}\mathbf{E}[|\mathcal{S}^{\beta}u(t,\cdot)|_{\mathcal{B}^{\tilde{s}}_{p,p}(w_{x_{0}}^{\tilde\theta})}^{p}]\lesssim M_{\sigma}^{p}\sup_{t\geq 0}\int_{0}^{t}e^{-\frac{\beta p}{2}(t-r)}\,dr<+\infty,
\end{equation}
which is exactly what we wanted to show in this step.


\noindent
\textit{Step 4: Regularity.} We first combine the bounds of Steps~1 and~2. Invoking the estimates \eqref{eq: deltaI} and \eqref{eq: D-pathwise}, and recalling the decomposition of $w$ from \eqref{eq:decomp-2}-\eqref{eq: I}-\eqref{eq: D}, we arrive at
\begin{equation}\label{eq: step12}
\sup_{r\in[T^{*},t]}\mathbf{E}\left[|\delta_{h}w(r,x)|^{p}\right]\lesssim \big(1+\mathbf{E}\big[|u_{0}|^{p\nu}_{C_{x_{0},\theta}}\big]\big)w_{x_{0}}^{\theta}(x)^{-\nu p} |h|^{\lambda p},
\end{equation}
where we have used the trivial bound $w_{x_{0}}^{\theta}(x)^{-p}\leq w_{x_{0}}^{\theta}(x)^{-\nu p}$. We recall the $\mathcal{B}^{\tilde{s}}_{p,p}(w_{x_{0}}^{\tilde\theta})$-norm from \eqref{eq: Besov-ppnorm} and we proceed to bound it like so:
\begin{equation}\label{eq: uBesov}\sup_{r\in [T^{*},t]}\mathbf{E}\left[ |u(r,\cdot)|^{p}_{\mathcal{B}^{\tilde{s}}_{p,p}(w_{x_0}^{\tilde\theta})}\right]\lesssim \sup_{r \in [T^{*},t]}
\mathbf{E}\left[
|u(r,\cdot)|_{L^{p}(w_{x_{0}}^{\tilde\theta})}^{p}
\right]
\;+\;
\sup_{r \in [T^{*},t]} \mathbf{E}\left[ 
\left[ u(r,\cdot) \right]_{\mathcal{B}^{\tilde{s}}_{p,p}(w_{x_0}^{\tilde\theta})}^{p}
\right].
\end{equation}
We start by bounding the $L^{p}(w_{x_{0}}^{\tilde\theta})$-term of \eqref{eq: uBesov}. Using the decomposition $\tilde\theta=\theta+(\tilde\theta-\theta)$, $\tilde\theta>\theta +\frac{d}{p}$ (which holds since $\tilde\theta>\nu\theta+\frac{d}{p}$ and $\nu>1$) and \eqref{est0} we arrive at
\begin{align}
\sup_{r \in [T^{*},t]}\mathbf{E}\left[ |w(r,\cdot)|_{L^{p}(w_{x_{0}}^{\tilde\theta})}^{p}
\right]&\lesssim \sup_{r \in [T^{*},t]}\mathbf{E}\left[\int_{\R^d}|u(r,x)|^{p}(w_{x_{0}}^{\theta+(\tilde\theta-\theta)}(x))^{p}\,dx \right]\nonumber\\
&\lesssim \sup_{r\in[T^{*},t]}\mathbf{E}\left[|u(r,\cdot)|_{C_{x_{0},\theta}}^{p}\right]\left(\int_{\R^d}\frac{dx}{(1+|x-x_{0}|^{2})^{\frac{(\tilde\theta-\theta)p}{2}}} \right)\nonumber \\
&\lesssim 1+\mathbf{E}\big[|u_{0}|^{p}_{C_{x_{0},\theta}}\big]<\infty.\label{eq: uLpnu}
\end{align}
Next, we bound the Besov seminorm term of \eqref{eq: uBesov}. We start with the  triangle inequality, $[u(r,\cdot)]_{\mathcal{B}_{p,p}^{\tilde{s}}(w_{x_{0}}^{\tilde{\theta}})}^{p}\lesssim[w(r,\cdot)]^{p}_{\mathcal{B}_{p,p}^{\tilde{s}}(w_{x_{0}}^{\tilde{\theta}})} +[\mathcal{S}^{\beta}u(r,\cdot)]^{p}_{\mathcal{B}_{p,p}^{\tilde{s}}(w_{x_{0}}^{\tilde{\theta}})}$ and follow by invoking the estimate \eqref{eq: step12} and choosing $\lambda>\tilde{s}$ alongside the assumption $\tilde{\theta}>\nu\theta +\frac{d}{p}$ of the statement to arrive at
\begin{align}\label{eq: uBesovseminorm}
\sup_{r \in [T^{*},t]} \mathbf{E}\left[
\left[ w(r,\cdot) \right]_{\mathcal{B}^{\tilde{s}}_{p,p}(w_{x_0}^{\tilde\theta})}^{p}\right]&\leq \sup_{r \in [T^{*},t]}  \int_{|h|\leq 1}\int_{\R^d}\mathbf{E}\left[|\delta_{h}w(r,\cdot)(x)|^{p}\right]\frac{(w_{x_{0}}^{\tilde\theta}(x))^{p}}{|h|^{d+\tilde{s}p}}dx\,dh \nonumber\\
&\lesssim\left(\int_{|h|\leq 1} \frac{dh}{|h|^{d+\tilde{s}p-\lambda p}}\right) \left(\int_{\R^d}\frac{dx}{(1+|x-x_{0}|^{2})^{\frac{\tilde\theta p}{2}-\frac{\nu\theta p}{2}}} \right)\nonumber\\
&<+\infty.
\end{align}
Also, notice that the second term $[\mathcal{S}^{\beta}u(r,\cdot)]^{p}_{\mathcal{B}_{p,p}^{\tilde{s}}(w_{x_{0}}^{\tilde{\theta}})}$ has already been bounded in \eqref{eq:Sbeta-Besov}. Putting together \eqref{eq: uLpnu} and \eqref{eq: uBesovseminorm} in \eqref{eq: uBesov} we obtain
\begin{equation}\label{eq:w-Besov}
\sup_{r\in [T^{*},t]}\mathbf{E}\left[ |w(r,\cdot)|^{p}_{\mathcal{B}^{\tilde{s}}_{p,p}(w_{x_0}^{\tilde\theta})}\right]<+\infty.
\end{equation}
Finally, combining \eqref{eq:Sbeta-Besov} and \eqref{eq:w-Besov} into the decomposition of $u$, we have
\begin{equation*}
\sup_{r\in[T^{*},t]}\mathbf{E}\left[|u(r,\cdot)|^{p}_{\mathcal{B}_{p,p}^{\tilde{s}}(w_{x_{0}}^{\tilde\theta})} \right]\leq C(T^{*},\mathbf{E}[|u_{0}|^{p\nu}_{C_{x_{0},\theta}}])<+\infty,
\end{equation*}
where $C(T^{*},\mathbf{E}[|u_{0}|^{p\nu}_{C_{x_{0},\theta}}])$ is a constant independent of $t$. Therefore, by letting $t\to \infty$, we have
\begin{equation*}
\sup_{r\in[T^{*},\infty)}\mathbf{E}\left[|u(r,\cdot)|^{p}_{\mathcal{B}^{\tilde{s}}_{p,p}(w_{x_0}^{\tilde\theta})} \right]<\infty,
\end{equation*}
which is what we wanted to show.
\end{proof}
\section{Existence of Invariant Measure}\label{sec:exist}

This section is devoted to the proof of Theorem~\ref{theorem: main}, by the classical Krylov-Bogoliubov strategy. In Section~\ref{sec:tightness} we combine the uniform Besov bound of Proposition~\ref{prop: est2} with the compact embedding of Lemma~\ref{lem: compact}. This gives the tightness of the time-averaged laws $\{R_{T}^{*}\mu\}_{T>0}$ in a space $C_{\rho}$, hence a weak limit $\mu_{\infty}$ along a sequence of times. In Section~\ref{sec:invariance} we prove that $\mu_{\infty}$ is invariant for the transition semigroup $\mathcal{P}_{t}$ of Definition~\ref{def: transition}. The key ingredient is the Feller property between the exponent $\rho$ and a larger exponent $\rho'$, established in Section~\ref{subsec:feller-two}. The existence of an invariant measure is then Theorem~\ref{thm:existence-invariant}.

\subsection{Tightness}\label{sec:tightness}
Armed with the uniform Besov bound of Proposition~\ref{prop: est2}, we now establish tightness of the time-averaged laws of the solution. Let us summarize the setting and the assumptions in force throughout this section.
\begin{enumerate}
\item We work under the conditions of Theorem~\ref{thm:extend} and Assumption~\ref{assump: dissip}, with $\theta>0$, $p\geq p^{*}$ in the sense of Remark~\ref{rmk:pstar}, and $\nu$ the exponent of \eqref{assump: growth}.
\item The initial condition $u_{0}$ satisfies Assumption~\ref{assump: initial-cond-fixed} together with the moment condition~\eqref{eq:initial-pnu} of Proposition~\ref{prop: est2}, that is $\mathbf{E}[|u_{0}|^{p\nu}_{C_{\theta}}]<\infty$. Every deterministic $\phi\in C_{\theta}$ qualifies.
\item The solution $u(t,\cdot\,;u_{0})$ is the one given by Theorem~\ref{thm:extend}.
\item All the spaces are taken with the center $x_{0}=0$, that is in the fixed-center spaces $C_{\theta}$, $C_{\rho}$ of~\eqref{eq:Ctheta}. Recall from Remark~\ref{rem: equiv-norms} that the spaces $C_{x_{0},\rho}$ coincide as sets for all centers, with equivalent norms, so that nothing below depends on this choice.
\end{enumerate}
We first fix the duality notation, then introduce the time-averaged laws.
\begin{definition}\label{def:duality-f-mu}
For $\rho>0$, let $\tfF\colon C_{\rho}\to \mathbb{R}$ be a bounded Borel measurable function and let $\mu$ be a bounded Borel measure on $C_{\rho}$. We write
\begin{equation}
\left\langle \tfF,\mu\right\rangle:=\int_{C_{\rho}}\tfF(v)\,\mu(dv).
\end{equation}
The same notation is used on the spaces $C_{x_{0},\rho}$ with a generic center $x_{0}\in\R^{d}$.
\end{definition}
\begin{definition}\label{def:average-meas}
Under the conditions of Theorem~\ref{thm:extend} and Assumption~\ref{assump: dissip}, let $u_{0}$ satisfy \eqref{eq:initial-pnu} and let $u(t,\cdot;u_{0})$ be the solution of \eqref{NSHE} given by Theorem~\ref{thm:extend}. Let $\mu$ be the law of $u_{0}$ and $\mu_{t}$ be the law of $u(t,\cdot;u_{0})$, both probability measures on $C_{\theta}$. For $T>0$, the \emph{time-averaged law} of $u(t,\cdot;u_{0})$ is the probability measure on $C_{\theta}$ defined by
\begin{equation}\label{eq:RT-def}
R_{T}^{*}\mu:=\frac{1}{T}\int_{0}^{T}\mu_{t}\,dt,
\qquad\text{that is}\qquad
\langle \tfF,R_{T}^{*}\mu\rangle=\frac{1}{T}\int_{0}^{T}\mathbf{E}\big[\tfF(u(t,\cdot\,;u_{0}))\big]\,dt
\end{equation}
for every bounded Borel $\tfF$ on $C_{\theta}$.
\end{definition}
\begin{remark}\label{rmk:RT-measurable}
The integral in \eqref{eq:RT-def} makes sense since $t\mapsto\mathbf{E}[\tfF(u(t,\cdot\,;u_{0}))]$ is continuous for $\tfF\in C_{b}(C_{\theta})$, by the continuity of the paths of $u$ and dominated convergence. We can then get a measurable path $t\mapsto\mathbf{E}[\tfF(u(t,\cdot\,;u_{0}))]$ for every bounded Borel $\tfF$ by a monotone class argument.
\end{remark}
To pass from the uniform Besov bound of Proposition~\ref{prop: est2} to tightness in $C_{x_{0},\rho}$, we invoke the compact embedding of Lemma~\ref{lem: compact}: for $\tilde{s}>d/p$, the space $\mathcal{B}_{p,p}^{\tilde{s}}(w_{x_{0}}^{\tilde{\theta}})$ embeds compactly into $C_{x_{0},\rho}$ for every $\rho>\tilde{\theta}$. Note that the strict inequality $\rho>\tilde{\theta}$ cannot be dispensed with. Indeed, the compactness of Lemma~\ref{lem: compact} comes with a loss in the weight exponent, and the balls of $\mathcal{B}_{p,p}^{\tilde{s}}(w_{x_{0}}^{\tilde{\theta}})$ are bounded, but not compact, in $C_{x_{0},\tilde{\theta}}$. The choice of the Besov parameters $\tilde{s}$ and $\tilde{\theta}$ in the proof of Proposition~\ref{prop:tight-average} is tailored to this requirement.

\begin{proposition}\label{prop:tight-average}
Under the conditions of Theorem~\ref{thm:extend} and Assumption~\ref{assump: dissip}, let $\nu$ be the exponent from the growth condition \eqref{assump: growth}, let $\theta>0$ and $\rho>\nu\theta$, and let $p\geq p^{*}$. Let $u_{0}$ satisfy \eqref{eq:initial-pnu} and let $\mu$ be the law of $u_{0}$. Then the family of time-averaged laws $\{R_{T}^{*}\mu\}_{T\geq 0}$ from Definition~\ref{def:average-meas} is tight in $C_{\rho}$.
\end{proposition}

\begin{proof}
Since $\rho>\nu\theta$ and $p\geq p^{*}$, we have $\nu\theta+d/p<\rho$ and $d/p<\eta-2/p$. We may thus fix $\tilde{\theta}\in(\nu\theta+d/p,\rho)$ and $\tilde{s}\in(d/p,\eta-2/p)$. Then Proposition~\ref{prop: est2} applies with these exponents, and $\tilde{s}>d/p$ allows us to apply the compact embedding of Lemma~\ref{lem: compact}.
We fix the center $x_{0}=0$, we write $w_{0}^{\theta}$ for the weight $w_{x_{0}}^{\theta}|_{x_{0}=0}$, we let $\varepsilon>0$ and consider the closed ball 
\begin{equation}
\overline{B}_{\mathcal{B}^{\tilde{s}}_{p,p}(w_{0}^{\tilde{\theta}})}(0,R):
=\left\{v\in \mathcal{B}^{\tilde{s}}_{p,p}(w_{0}^{\tilde{\theta}}): |v|_{\mathcal{B}^{\tilde{s}}_{p,p}(w_{0}^{\tilde{\theta}})}\leq R
\right\},
\end{equation}
which is compact in $C_{x_{0},\rho}$ for $\rho>\tilde{\theta}$, by Lemma~\ref{lem: compact}. Moreover, by estimate~\eqref{est2} of Proposition~\ref{prop: est2}, the constant
\begin{equation*}
C_{T^{*}}:=\sup_{t\geq T^{*}}\mathbf{E}\!\left[|u(t,\cdot)|^{p}_{\mathcal{B}^{\tilde{s}}_{p,p}(w_{0}^{\tilde{\theta}})}\right]
\end{equation*}
is finite. Applying Markov's inequality to the event $\{|u(t,\cdot)|_{\mathcal{B}^{\tilde{s}}_{p,p}(w_{0}^{\tilde{\theta}})}>R\}$, we arrive at
\begin{equation}\label{j1}
\mathbf{P}\left(u(t)\notin \overline{B}_{\mathcal{B}^{\tilde{s}}_{p,p}(w_{0}^{\tilde{\theta}})}(0,R)\right)\leq \frac{1}{R^{p}}\sup_{t\geq T^{*}}\mathbf{E}\!\left[|u(t,\cdot)|^{p}_{\mathcal{B}^{\tilde{s}}_{p,p}(w_{0}^{\tilde{\theta}})} \right]=\frac{C_{T^{*}}}{R^{p}}.
\end{equation}
Choosing $R=R(\varepsilon)$ large enough implies tightness for $\{\mu_{t}\}_{t\geq T^{*}}$, that is
\begin{equation}\label{eq: u-tight}
\mathbf{P}\left(u(t)\in \overline{B}_{\mathcal{B}^{\tilde{s}}_{p,p}(w_{0}^{\tilde{\theta}})}(0,R)\right)\geq 1-\varepsilon,
\quad\text{for all}\quad t\geq T^{*}.
\end{equation}
The passage from the tightness \eqref{eq: u-tight} of the laws $\mu_{t}$, $t\geq T^{*}$, to the tightness of the averaged laws $R_{T}^{*}\mu$, $T>0$, is a matter of standard considerations, see e.g. \cite[Section~6.1]{DPZ96} or \cite[Chapter~8]{cerrai:01:second}.
\end{proof}
The following lemma draws the two main consequences of Proposition~\ref{prop:tight-average}. It is standard, see e.g. \cite[Section~3.1]{DPZ96} or \cite[Section~2.1]{cerrai:01:second}, and we include the short proof for the reader's convenience.

\begin{lemma}\label{lem:convergence-mu-t}
Under the conditions of Proposition~\ref{prop:tight-average} (in particular $\rho>\nu\theta$ and $p\geq p^{*}$), let $u_{0}$ satisfy \eqref{eq:initial-pnu} and let $\mu$ be the law of $u_{0}$. Then:
\begin{enumerate}
\item\label{it:j2}
There exist a subsequence $\{T_{n}\}_{n\geq 1}$ with $T_{n}\to\infty$ and a probability measure $\mu_{\infty}$ on $C_{\rho}$ such that $R_{T_{n}}^{*}\mu\Rightarrow\mu_{\infty}$ weakly in $C_{\rho}$.
\item\label{it:j3} Every such limiting measure satisfies the moment bound 
\begin{equation}\label{eq:mu-infty-moment}
\int_{C_{\rho}}|v|^{p\nu}_{C_{\theta}}\,\mu_{\infty}(dv)\leq C\Big(1+\mathbf{E}\big[|u_{0}|^{p\nu}_{C_{\theta}}\big]\Big)<\infty,
\end{equation}
where $C$ is the constant of \eqref{est0} for the exponent $p\nu$, and where $|v|_{C_{\theta}}:=\sup_{x\in\R^{d}}|v(x)|w_{0}^{\theta}(x)$ is understood as $+\infty$ when this supremum is infinite. In particular, one has $\mu_{\infty}(\{v\in C_{\rho}:|v|_{C_{\theta}}<\infty\})=1$.
\end{enumerate}
\end{lemma}

\begin{proof}
\noindent\textit{Proof of item~\ref{it:j2}.}
By Proposition~\ref{prop:tight-average}, the family $\{R_{T}^{*}\mu\}_{T>0}$ is tight in $C_{\rho}$. By Prokhorov's theorem, tightness implies relative compactness for the weak topology, so there exist a subsequence $\{T_{n}\}_{n\geq 1}$ with $T_{n}\to\infty$ and a probability measure $\mu_{\infty}$ on $C_{\rho}$ such that $R_{T_{n}}^{*}\mu\Rightarrow\mu_{\infty}$ weakly in $C_{\rho}$.

\smallskip
\noindent\textit{Proof of item~\ref{it:j3}.}
Set $\tfPsi(v):=|v|^{p\nu}_{C_{\theta}}\in[0,\infty]$ for $v\in C_{\rho}$. The function $\tfPsi$ is lower semicontinuous on $C_{\rho}$: it is the supremum over $x\in\R^{d}$ of the functions $v\mapsto(|v(x)|w_{0}^{\theta}(x))^{p\nu}$, each of which is continuous on $C_{\rho}$ since $|v(x)|\leq|v|_{C_{\rho}}/w_{0}^{\rho}(x)$. Hence, by the weak convergence $R_{T_{n}}^{*}\mu\Rightarrow\mu_{\infty}$ of item~\ref{it:j2} and the Portmanteau theorem (applied to the bounded lower semicontinuous functions $\tfPsi\wedge M$, and then letting $M\to\infty$ by monotone convergence), we have
\begin{equation*}
\int_{C_{\rho}}\tfPsi(v)\,\mu_{\infty}(dv)\leq\liminf_{n\to\infty}\left<\tfPsi,R_{T_{n}}^{*}\mu\right>
=\liminf_{n\to\infty}\frac{1}{T_{n}}\int_{0}^{T_{n}}\mathbf{E}\left[|u(t,\cdot)|^{p\nu}_{C_{\theta}}\right]dt ,
\end{equation*}
where the equality is the definition \eqref{eq:RT-def} of the averaged measure. The uniform-in-time bound \eqref{est0} of Proposition~\ref{prop: est1}, written with $p\nu$ in place of $p$ and with the center $x_{0}=0$, bounds the expectation by $C(1+\mathbf{E}[|u_{0}|^{p\nu}_{C_{\theta}}])$ for every $t\geq 0$, which is \eqref{eq:mu-infty-moment}. This completes the proof of item \ref{it:j3}.
\end{proof}

\subsection{Existence of an invariant measure}\label{sec:invariance}
We now identify the limit $\mu_{\infty}$ of Lemma~\ref{lem:convergence-mu-t} as an invariant measure for \eqref{NSHE}. We first fix the exponents, then run the classical Krylov-Bogoliubov argument with the two-exponent Feller property of Section~\ref{subsec:feller-two}.

\smallskip

\noindent
{\it The exponents.} Three weight exponents are involved. The exponent $\theta$ of the state space $C_{\theta}$ is the one of the moment bounds of Section~\ref{sec:bounds}. Proposition~\ref{prop:tight-average} gives the tightness in $C_{\rho}$ for every $\rho>\nu\theta$, and the Feller property of Proposition~\ref{prop:feller-two} is available between the exponent $\rho$ and a larger exponent $\rho'$. We thus fix, once and for all,
\begin{equation}\label{eq:exponent-chain}
\nu\theta<\rho<\rho'
\qquad\text{and}\qquad
p\geq p^{*},
\end{equation}
where $p^{*}$ is such that Proposition~\ref{prop:tight-average} and Lemma~\ref{lem:convergence-mu-t} apply with the exponent $\rho$, and that the triple $(\rho,\rho',p)$ satisfies \eqref{eq:two-exponents}, so that the results of Section~\ref{subsec:feller-two} apply as stated and Theorem~\ref{thm:extend} holds in the spaces $C_{\rho}$ and $C_{\rho'}$. Observe also that the solution $u(\cdot\,;u_{0})$ given by Theorem~\ref{thm:extend}, which has paths in $C([0,T];C_{\theta})$, is also 
a solution living in  $C_{\rho}$ and  $C_{\rho'}$. In particular $u(t,\cdot\,;u_{0})$ is a $C_{\rho}$-valued random variable, and $\mu_{t}=\mathrm{Law}(u(t,\cdot\,;u_{0}))$ is a probability measure on $C_{\rho}$ for every $t\geq 0$, and so is $R_{T}^{*}\mu$.

\smallskip

\noindent
{\it The transition semigroup on $C_{\rho}$.} With these exponents, Proposition~\ref{prop:feller-two} provides the transition kernel $P_{t}$ on $C_{\rho}$ of \eqref{eq:kernel-def}, the semigroup $\mathcal{P}_{t}$ on $B_{b}(C_{\rho})$, and the Feller property between the two exponents
\begin{equation}\label{eq:feller-rho}
\mathcal{P}_{t}\big(C_{b}(C_{\rho'})\big)\big|_{C_{\rho}}\subset C_{b}(C_{\rho}),\qquad t\geq 0 ,
\end{equation}
which is \eqref{eq:feller-two}. For a probability measure $\mu$ on $C_{\rho}$, the measure $\mathcal{P}_{t}^{*}\mu$ on $C_{\rho}$ is given by the kernel formula \eqref{eq:Pt-star-kernel} of Lemma~\ref{lem:plug-in}, and it satisfies the duality relation
\begin{equation}\label{eq:duality-rho}
\langle F,\mathcal{P}_{t}^{*}\mu\rangle=\langle\mathcal{P}_{t}F,\mu\rangle,\qquad F\in B_{b}(C_{\rho}),
\end{equation}
which is the first identity in \eqref{eq: Pt}. Moreover, by item~\ref{it:pi-plug} of Lemma~\ref{lem:plug-in} applied to the $C_{\rho}$-valued initial condition $u_{0}$, the laws $\mu_{t}$ of Definition~\ref{def:average-meas} are $\mu_{t}=\mathcal{P}_{t}^{*}\mu$, so that the averaged measure \eqref{eq:RT-def} reads
\begin{equation}\label{eq:RT-semigroup}
R_{T}^{*}\mu=\frac{1}{T}\int_{0}^{T}\mathcal{P}_{t}^{*}\mu\,dt,
\qquad\text{that is}\qquad
\langle F,R_{T}^{*}\mu\rangle=\frac{1}{T}\int_{0}^{T}\langle\mathcal{P}_{t}F,\mu\rangle\,dt,\quad F\in B_{b}(C_{\rho}).
\end{equation}
Finally, the semigroup property $\mathcal{P}_{t+s}=\mathcal{P}_{t}\mathcal{P}_{s}$ holds on $B_{b}(C_{\rho})$, by item (i) of Remark~\ref{rmk:markov-property} and item (c) of Remark~\ref{rmk:feller-two}. We say that a probability measure $\mu_{\infty}$ on $C_{\rho}$ is \emph{invariant} for the semigroup if $\mathcal{P}_{s}^{*}\mu_{\infty}=\mu_{\infty}$ for every $s\geq 0$.

\smallskip

We can now state and prove the main result of this paper, which gives a precise form to Theorem~\ref{theorem: main}.
\begin{theorem}\label{thm:existence-invariant}
Let Assumptions~\ref{assump:noise}, \ref{assump:reac}, \ref{assump: dissip} and \ref{assump:diff} hold, let $\theta>0$ and $\rho>\nu\theta$ be two weight exponents, where $\nu$ is the exponent of the growth condition \eqref{assump: growth}, and let $p\geq p^{*}$. Then the following holds true.
\begin{enumerate}
\item\label{it:inv-limits} Let $u_{0}$ be an initial condition satisfying Assumption~\ref{assump: initial-cond-fixed} and the moment condition \eqref{eq:initial-pnu}, and let $\mu=\mathrm{Law}(u_{0})$. By Proposition~\ref{prop:tight-average} and Lemma~\ref{lem:convergence-mu-t}, the family $\{R_{T}^{*}\mu\}_{T>0}$ of averaged measures \eqref{eq:RT-semigroup} is tight in $C_{\rho}$, and admits weak limits $\mu_{\infty}$ in $C_{\rho}$ along sequences $T_{n}\to\infty$. Every such limit $\mu_{\infty}$ is an invariant probability measure for the transition semigroup $(\mathcal{P}_{t})_{t\geq 0}$ of \eqref{NSHE} on $C_{\rho}$.
\item\label{it:inv-exist} In particular, for every deterministic initial condition $u_{0}\in C_{\theta}$, the averaged measures $R_{T}^{*}\delta_{u_{0}}$ admit weak limits in $C_{\rho}$, and every such limit is an invariant probability measure. Hence the transition semigroup of \eqref{NSHE} admits an invariant probability measure on $C_{\rho}$, for every $\theta>0$ and every $\rho>\nu\theta$.
\end{enumerate}
\end{theorem}

\begin{proof}
We take $\rho':=\rho+1$ in \eqref{eq:exponent-chain}, so that $p^{*}$ does not depend on $\rho'$. Item~\ref{it:inv-exist} follows from item~\ref{it:inv-limits}, since a deterministic $u_{0}\in C_{\theta}$ satisfies Assumption~\ref{assump: initial-cond-fixed} and \eqref{eq:initial-pnu}. We thus focus on the proof of item~\ref{it:inv-limits}. To this aim, we fix a weak limit $\mu_{\infty}$ of $R_{T_{n}}^{*}\mu$ in $C_{\rho}$. We have to prove that $\mathcal{P}_{s}^{*}\mu_{\infty}=\mu_{\infty}$ for every $s\geq 0$. We fix $s\geq 0$ and proceed in three steps, following the classical argument (see \cite[Section~3.1]{DPZ96}), the only difference being that the test functions are continuous on the larger space $C_{\rho'}$.

\smallskip

\noindent
{\it Step 1: Test functions.} Let $F\in C_{b}(C_{\rho'})$, which we restrict to $C_{\rho}$. By item~\ref{it:ts-inclusion} of Lemma~\ref{lem:two-spaces}, the inclusion of $C_{\rho}$ into $C_{\rho'}$ is continuous, so that $F|_{C_{\rho}}\in C_{b}(C_{\rho})$, and by the Feller property \eqref{eq:feller-rho} of Proposition~\ref{prop:feller-two} we also have $\mathcal{P}_{s}F|_{C_{\rho}}\in C_{b}(C_{\rho})$. Note that in the sequel of the proof we avoid writing the restrictions $F|_{C_{\rho}}$ and $(P_{t}F)|_{C_{\rho}}$ for aesthetic reasons. The weak convergence $R_{T_{n}}^{*}\mu\Rightarrow\mu_{\infty}$ in $C_{\rho}$ thus applies to both functions, and together with the duality \eqref{eq:duality-rho} of Lemma~\ref{lem:plug-in} it gives
\begin{equation}\label{eq:KB-1}
\langle F,\mathcal{P}_{s}^{*}\mu_{\infty}\rangle=\langle\mathcal{P}_{s}F,\mu_{\infty}\rangle=\lim_{n\to\infty}\langle\mathcal{P}_{s}F,R_{T_{n}}^{*}\mu\rangle,
\qquad
\langle F,\mu_{\infty}\rangle=\lim_{n\to\infty}\langle F,R_{T_{n}}^{*}\mu\rangle .
\end{equation}

\smallskip

\noindent
{\it Step 2: The shift.} Let $T>0$. Applying \eqref{eq:RT-semigroup} to the function $\mathcal{P}_{s}F\in B_{b}(C_{\rho})$, and then by the semigroup property of item (i) of Remark~\ref{rmk:markov-property}, we have
\begin{equation}\label{eq:KB-shift}
\langle\mathcal{P}_{s}F,R_{T}^{*}\mu\rangle=\frac{1}{T}\int_{0}^{T}\langle\mathcal{P}_{t}\mathcal{P}_{s}F,\mu\rangle\,dt
=\frac{1}{T}\int_{0}^{T}\langle\mathcal{P}_{t+s}F,\mu\rangle\,dt
=\frac{1}{T}\int_{s}^{T+s}\langle\mathcal{P}_{t}F,\mu\rangle\,dt .
\end{equation}
We now split the last integral in \eqref{eq:KB-shift} according to $[s,T+s]=[0,T]\cup[T,T+s]\setminus[0,s)$, which gives
\begin{equation}\label{eq:KB-split}
\langle\mathcal{P}_{s}F,R_{T}^{*}\mu\rangle
=\frac{1}{T}\int_{0}^{T}\langle\mathcal{P}_{t}F,\mu\rangle\,dt
+\frac{1}{T}\int_{T}^{T+s}\langle\mathcal{P}_{t}F,\mu\rangle\,dt
-\frac{1}{T}\int_{0}^{s}\langle\mathcal{P}_{t}F,\mu\rangle\,dt .
\end{equation}
The first term in the right hand side of \eqref{eq:KB-split} is $\langle F,R_{T}^{*}\mu\rangle$, by \eqref{eq:RT-semigroup} applied to $F$. In the two other terms, the integrand is bounded by $|\mathcal{P}_{t}F|_{C_{b}(C_{\rho})}\leq|F|_{C_{b}(C_{\rho'})}$, since $\mathcal{P}_{t}$ is a contraction on bounded functions and $|F|_{C_{b}(C_{\rho})}\leq|F|_{C_{b}(C_{\rho'})}$. In addition, the integrals run over intervals of length $s$. Hence we obtain
\begin{equation}\label{eq:KB-2}
\big|\langle\mathcal{P}_{s}F,R_{T}^{*}\mu\rangle-\langle F,R_{T}^{*}\mu\rangle\big|\leq\frac{2s}{T}\,|F|_{C_{b}(C_{\rho'})},
\qquad T>0 .
\end{equation}
Taking $T=T_{n}$ in and letting $n\to\infty$, the right hand side of \eqref{eq:KB-2} vanishes, while by \eqref{eq:KB-1} the two terms in the left hand side converge to $\langle F,\mathcal{P}_{s}^{*}\mu_{\infty}\rangle$ and $\langle F,\mu_{\infty}\rangle$ respectively. We thus get
\begin{equation}\label{eq:KB-3}
\langle F,\mathcal{P}_{s}^{*}\mu_{\infty}\rangle=\langle F,\mu_{\infty}\rangle,
\qquad\text{for every }F\in C_{b}(C_{\rho'})\text{ restricted to } C_{\rho}.
\end{equation}

\smallskip

\noindent
{\it Step 3: From test functions to Borel sets.} The passage from \eqref{eq:KB-3} to the equality of $\mathcal{P}_{s}^{*}\mu_{\infty}$ and $\mu_{\infty}$ on the bounded Borel functions is a classical application of the monotone class theorem. We give the details, because of the delicate point that the test functions in \eqref{eq:KB-3} are the restrictions to $C_{\rho}$ of the functions of $C_{b}(C_{\rho'})$, a class which is strictly smaller than $C_{b}(C_{\rho})$ (see \eqref{eq:restriction-incl}). We proceed in three parts.
\begin{enumerate}[label=\textnormal{(\alph*)}]
\item \emph{The class of test functions.} Let $\mathcal{C}:=\{F|_{C_{\rho}}\,;\,F\in C_{b}(C_{\rho'})\}$. It is a vector space of bounded functions on $C_{\rho}$, closed under multiplication and containing the constants, and by \eqref{eq:KB-3} the two probability measures $\mathcal{P}_{s}^{*}\mu_{\infty}$ and $\mu_{\infty}$ on $C_{\rho}$ integrate every element of $\mathcal{C}$ in the same way.
\item \emph{The $\sigma$-algebra generated by $\mathcal{C}$.} The functions of $C_{b}(C_{\rho'})$ generate the Borel $\sigma$-algebra of the metric space $C_{\rho'}$, since every closed set $K\subset C_{\rho'}$ is the zero set of the bounded continuous function $\mathrm{dist}(\cdot,K)\wedge 1$. Hence the class $\mathcal{C}$ generates the trace on $C_{\rho}$ of the Borel $\sigma$-algebra of $C_{\rho'}$, which is the Borel $\sigma$-algebra of $C_{\rho}$ by item~\ref{it:ts-borel} of Lemma~\ref{lem:two-spaces}.
\item \emph{The monotone class argument.} The set $\mathcal{H}$ of bounded Borel functions $G$ on $C_{\rho}$ such that $\langle G,\mathcal{P}_{s}^{*}\mu_{\infty}\rangle=\langle G,\mu_{\infty}\rangle$ is a vector space containing the constants, closed under bounded monotone limits by the monotone convergence theorem, and containing $\mathcal{C}$ by part (a). Since $\mathcal{C}$ is closed under multiplication and generates the Borel $\sigma$-algebra of $C_{\rho}$ by part (b), the monotone class theorem asserts that $\mathcal{H}$ is the whole of $B_{b}(C_{\rho})$. Taking $G=\mathbf{1}_{A}$ for a Borel set $A$ gives $\mathcal{P}_{s}^{*}\mu_{\infty}(A)=\mu_{\infty}(A)$.
\end{enumerate}
This proves $\mathcal{P}_{s}^{*}\mu_{\infty}=\mu_{\infty}$, and the proof is complete.
\end{proof}

\begin{remark}\label{rmk:invariant-comments}
Let us comment on Theorem~\ref{thm:existence-invariant}.
\begin{enumerate}[label=\textnormal{\textbf{(\alph*)}}]
\item The invariant measure lives on $C_{\rho}$, a space larger than the state space $C_{\theta}$ of the moment bounds, and its invariance is for the semigroup of Proposition~\ref{prop:feller-two} on $C_{\rho}$. The condition $\rho>\nu\theta$ has two sources. The factor $\nu$ comes from the polynomial growth of $f$, through the Besov bound of Proposition~\ref{prop: est2}, which requires $\tilde{\theta}>\nu\theta+d/p$. The strict inequality $\rho>\tilde{\theta}$ comes from the compact embedding of Lemma~\ref{lem: compact}, which loses an arbitrarily small amount of weight exponent (see the discussion before Proposition~\ref{prop:tight-average}). The larger exponent $\rho'$ only enters through the class of test functions used in the proof. By \eqref{eq:mu-infty-moment}, $\mu_{\infty}$ is nevertheless concentrated on the functions $v$ with $|v|_{C_{\theta}}<\infty$.
\item The Feller property on a single space, namely $\mathcal{P}_{t}(C_{b}(C_{\rho}))\subset C_{b}(C_{\rho})$, is not available (see item (a) of Remark~\ref{rmk:feller-two}), and the proof shows that it is not needed: the weak convergence of the averaged measures is only tested on functions which are continuous for the weaker topology of $C_{\rho'}$, and these functions still generate the Borel $\sigma$-algebra of $C_{\rho}$.
\item The moment condition \eqref{eq:initial-pnu} on the initial condition in item~\ref{it:inv-limits} is only used to produce the tightness, and item~\ref{it:inv-exist} is the case of a deterministic initial condition, for which it is void. The random initial conditions covered by item~\ref{it:inv-limits} are thus those with a finite moment of order $p\nu$ in $C_{\theta}$ for some $p\geq p^{*}$, for instance those with moments of all orders. 
\end{enumerate}
\end{remark}

\section*{Acknowledgements}
\noindent
EA would like to thank Le Chen for discussions on comparison principles and Mickey Salins for discussions on weighted stochastic convolution estimates. Both authors were partially supported by the NSF
grant DMS-2450734.

\section*{AI declaration}

Most of this work was carried out by the authors themselves. This was the case up until a gap appeared in the proof of the Feller property (Section~\ref{sec:feller}), where the parameters did not match. At that point Claude (Anthropic) helped the authors mathematically to close this gap.
Furthermore, during the preparation of this work, the authors used Claude to assist with proofreading, formatting, checking
mathematical cross-references, reviewing mathematical arguments, and proposing
corrections to proofs. The authors take full responsibility for the content
of the article.

\appendix

\section{Weighted spaces and the heat semigroup}\label{app:weights}

This section gathers the proofs of the elementary lemmas on the weights \eqref{eq: weight1} and on the space $C_{\theta}$ of \eqref{eq:Ctheta} which are used throughout the paper. Namely, we prove the conjugation identity of Lemma~\ref{lem:conjugation}, the two properties of the fixed-center space stated in Lemmas~\ref{lem:fc-paths} and~\ref{lem:fc-heat}, and the weight inequalities of Lemmas~\ref{lem:weight-absorb} and~\ref{lem:weights}.

\begin{proof}[Proof of Lemma~\ref{lem:conjugation}]
We start with item \ref{it:conj-identity}. Set $\rho(x):=(1+|x-x_{0}|^{2})^{\theta/2}$, so that $\rho=1/\hat{w}_{x_{0}}^{\theta}$, and write $\phi=\rho\psi$ with $\psi:=\hat{w}_{x_{0}}^{\theta}\phi$. A direct computation gives
\begin{equation*}
\nabla\rho(x)=\theta\,(x-x_{0})\,(1+|x-x_{0}|^{2})^{\frac{\theta}{2}-1},
\quad
\Delta\rho(x)=\Big(\theta d+\theta(\theta+d-2)|x-x_{0}|^{2}\Big)(1+|x-x_{0}|^{2})^{\frac{\theta}{2}-2}.
\end{equation*}
Plugging those expressions into the product rule $\Delta(\rho\psi)=\rho\Delta\psi+2\nabla\rho\cdot\nabla\psi+\psi\Delta\rho$ and multiplying by $\hat{w}_{x_{0}}^{\theta}/2=1/(2\rho)$, we get
\begin{equation*}
\hat{w}_{x_{0}}^{\theta}\,\frac{1}{2}\Delta\phi
=\frac{1}{2}\Delta\psi+\frac{\nabla\rho}{\rho}\cdot\nabla\psi+\frac{\Delta\rho}{2\rho}\,\psi
=\frac{1}{2}\Delta\psi+\frac{\theta\,(x-x_{0})\cdot\nabla\psi}{1+|x-x_{0}|^{2}}+c_{\theta}\,\psi ,
\end{equation*}
which is \eqref{eq:conjugation}.

We now turn to item \ref{it:conj-max}. The numerator of $c_{\theta}$ in \eqref{eq:Ltheta-coef} grows like $|x|^{2}$ and its denominator like $|x|^{4}$, so that $c_{\theta}$ is bounded. Next let $\psi$ and $x^{*}$ be as in the statement. Since $\psi$ is nonzero, we have $\psi(x^{*})\neq 0$. If $\psi(x^{*})>0$, then $x^{*}$ is a global maximizer of $\psi$, so that $\nabla\psi(x^{*})=0$ and the Hessian of $\psi$ at $x^{*}$ is nonpositive, which implies $\Delta\psi(x^{*})\leq 0$. If $\psi(x^{*})<0$, then $x^{*}$ is a global minimizer and the two inequalities are reversed. In both cases,
\begin{equation}\label{eq:max-conditions}
\nabla\psi(x^{*})=0,\qquad\text{and}\qquad \mathrm{sign}(\psi(x^{*}))\,\Delta\psi(x^{*})\leq 0 .
\end{equation}
Multiplying \eqref{eq:Ltheta-def} at $x=x^{*}$ by $\mathrm{sign}(\psi(x^{*}))$ and using \eqref{eq:max-conditions}, the only remaining term is $c_{\theta}(x^{*})|\psi(x^{*})|$, which is bounded by the right hand side of \eqref{eq:Ltheta-max}.
\end{proof}

\begin{proof}[Proof of Lemma~\ref{lem:fc-paths}]
The equality of norms is the definition \eqref{eq: wsnorm2} of the norm of $C_{0,\theta}([0,T]\times\R^{d})$, so we only have to check that the two sets coincide. We prove the two inclusions in two steps.

\smallskip

\noindent
{\it Step 1: Space-time functions are continuous paths.} Let $\psi\in C_{0,\theta}([0,T]\times\mathbb{R}^d)$ and fix $\varepsilon>0$. By \eqref{eq: weightSpace2} there exists $R>0$ such that $|\psi(t,x)|\le\varepsilon(1+|x|^\theta)$ for all $t\in[0,T]$ and $|x|\ge R$. In particular $\psi(t,\cdot)\in C_\theta$ for every $t$. Moreover, since the weight in \eqref{eq:Ctheta} is bounded by $1$, for $s,t\in[0,T]$ we get
\begin{equation*}
|\psi(t,\cdot)-\psi(s,\cdot)|_{C_\theta}\le\sup_{|x|\le R}|\psi(t,x)-\psi(s,x)|+2\varepsilon .
\end{equation*}
The function $\psi$ is uniformly continuous on the compact set $[0,T]\times\{|x|\le R\}$. Hence $\limsup_{s\to t}|\psi(t,\cdot)-\psi(s,\cdot)|_{C_\theta}\le2\varepsilon$, and since $\varepsilon$ is arbitrary, $\psi\in C([0,T];C_\theta)$.

\smallskip

\noindent
{\it Step 2: Continuous paths are space-time functions.} Conversely, let $\psi\in C([0,T];C_\theta)$. For $(t,x),(s,y)\in[0,T]\times\mathbb{R}^d$ we have
\begin{equation*}
|\psi(t,x)-\psi(s,y)|\le(1+|x|^\theta)\,|\psi(t,\cdot)-\psi(s,\cdot)|_{C_\theta}+|\psi(s,x)-\psi(s,y)|,
\end{equation*}
which shows that $\psi\in C([0,T]\times\mathbb{R}^d)$. Next, consider the set $K:=\{\psi(t,\cdot):t\in[0,T]\}$. It is the image of the compact interval $[0,T]$ under the map $t\mapsto\psi(t,\cdot)$, which is continuous with values in $C_{\theta}$. Hence $K$ is a compact subset of $C_{\theta}$. Given $\varepsilon>0$, we can thus choose a finite $\varepsilon$-net $\{\psi(t_i,\cdot)\}_{i\le N}$ of $K$. That is, for every $t\in[0,T]$ there exists $i\le N$ such that $|\psi(t,\cdot)-\psi(t_{i},\cdot)|_{C_{\theta}}\le\varepsilon$. Owing to the definition \eqref{eq:Ctheta} of the norm, this yields, for every $x\in\R^{d}$,
\begin{equation*}
\sup_{t\in[0,T]}\frac{|\psi(t,x)|}{1+|x|^\theta}\le\varepsilon+\max_{i\le N}\frac{|\psi(t_i,x)|}{1+|x|^\theta} .
\end{equation*}
The right hand side converges to $\varepsilon$ as $|x|\to\infty$, since each $\psi(t_i,\cdot)$ belongs to $C_\theta$. As $\varepsilon$ is arbitrary, $\psi$ satisfies \eqref{eq: weightSpace2} with $x_0=0$, which finishes the proof.
\end{proof}

\begin{proof}[Proof of Lemma~\ref{lem:fc-heat}]
We divide the proof in two steps. The first one establishes the bound \eqref{eq:U0Bound}, and the second one the continuity of $t\mapsto H_{t}\phi$ in $C_{\theta}$, by approximating $\phi$ with compactly supported functions.

\smallskip

\noindent
{\it Step 1: The bound \eqref{eq:U0Bound}.} Using the weight inequality \eqref{eq:weight-triangle} with center $x_{0}=0$, that is $1+|y|^\theta\le c_\theta\,(1+|x|^\theta)(1+|x-y|^\theta)$, in \eqref{eq:heatsemigroup}, together with $|\phi(y)|\le|\phi|_{C_\theta}(1+|y|^\theta)$, we get
\begin{equation}\label{eq:U0Bound-pre}
|H_t\phi(x)|\le c_\theta(1+|x|^\theta)\,|\phi|_{C_\theta}\int_{\mathbb{R}^d}G(t,w)(1+|w|^\theta)\,dw .
\end{equation}
The integral in \eqref{eq:U0Bound-pre} equals $1+m_{\theta}t^{\theta/2}$ by item \ref{it:G-moment} of Definition~\ref{def: heatkernel}. Plugging this identity into \eqref{eq:U0Bound-pre} and dividing by $1+|x|^\theta$ yields \eqref{eq:U0Bound}.

\smallskip

\noindent
{\it Step 2: Continuity in $C_{\theta}$.} We now prove item \ref{it:fc-heat-cont}, in three parts.
\begin{enumerate}[label=\textnormal{(\alph*)}]
\item \emph{Density.} The space $C_{0}$ from \eqref{eq: vanishinfty} is dense in $C_\theta$. Indeed, for $R>0$ let $\chi_R$ be a continuous cutoff function with $\mathbf 1_{\{|x|\le R\}}\le\chi_R\le\mathbf 1_{\{|x|\le2R\}}$. Then
\begin{equation}\label{eq:fc-density}
|\phi-\phi\chi_R|_{C_\theta}\le\sup_{|x|\ge R}\frac{|\phi(x)|}{1+|x|^\theta},
\end{equation}
and the right hand side of \eqref{eq:fc-density} vanishes as $R\to\infty$, by the definition \eqref{eq: weightSpace} of $C_{\theta}$. In the sequel we set $\phi_{n}:=\phi\chi_{n}$ for $n\geq 1$.
\item \emph{Compactly supported data.} Since $\phi_n$ is continuous with compact support, $H_t\phi_n\in C_0(\mathbb{R}^d)$ for every $t\ge0$. Moreover, standard estimates for the heat kernel show that $t\mapsto H_t\phi_n$ is continuous in the norm \eqref{eq:normBound}. 

\item \emph{Conclusion.} Let $t\in[0,T]$. By \eqref{eq:U0Bound} applied to $\phi-\phi_{n}$, and then by \eqref{eq:fc-density},
\begin{equation}\label{eq:fc-approx}
|H_t\phi-H_t\phi_n|_{C_\theta}\le c_\theta\big(1+m_\theta T^{\theta/2}\big)|\phi-\phi_n|_{C_\theta}\longrightarrow 0,\qquad\text{as }n\to\infty .
\end{equation}
By part (b), $H_t\phi_n\in C_0\subset C_\theta$.  Since the vanishing condition in \eqref{eq: weightSpace} passes to limits in the norm $|\cdot|_{C_{\theta}}$, relation \eqref{eq:fc-approx} shows that $H_{t}\phi\in C_{\theta}$. As for the continuity, we write, for $s,t\in[0,T]$,
\begin{equation*}
H_t\phi-H_s\phi=H_t(\phi-\phi_n)+(H_t\phi_n-H_s\phi_n)+H_s(\phi_n-\phi) .
\end{equation*}
Invoking \eqref{eq:fc-approx} for the first and third terms, we obtain
\begin{equation*}
|H_t\phi-H_s\phi|_{C_\theta}\le2c_\theta\big(1+m_\theta T^{\theta/2}\big)|\phi-\phi_n|_{C_\theta}+|H_t\phi_n-H_s\phi_n|_{C_{b}},\qquad s,t\in[0,T].
\end{equation*}
Choosing first $n$ large, and then $|t-s|$ small thanks to part (b), proves the continuity of $t\mapsto H_{t}\phi$ in $C_{\theta}$.
\end{enumerate}
This completes the proof of item \ref{it:fc-heat-cont}.
\end{proof}

\begin{proof}[Proof of Lemma~\ref{lem:weight-absorb}]
We write $z=\tau^{1/2}y$. By the explicit expression \eqref{eq: heatkernel} of $G$, we have
\begin{equation*}
G(\tau,z)^{a}=(2\pi\tau)^{-da/2}e^{-a|y|^{2}/2},
\qquad
G(2\tau,z)^{a}=(4\pi\tau)^{-da/2}e^{-a|y|^{2}/4},
\end{equation*}
and
\begin{equation*}
1+|z|^{\theta_{0}}\leq\bigl(1+\tau^{\theta_{0}/2}\bigr)\bigl(1+|y|^{\theta_{0}}\bigr).
\end{equation*}
Hence \eqref{eq:absorb-power} follows from the boundedness of $y\mapsto(1+|y|^{\theta_{0}})e^{-a|y|^{2}/4}$ on $\R^{d}$. As for \eqref{eq:absorb-laplace}, by \eqref{eq:laplaceG} and the same change of variables,
\begin{equation*}
\bigl(1+|z|^{\theta_{0}}\bigr)|\Delta_{z}G(\tau,z)|
\leq\frac{1}{\tau}\bigl(1+\tau^{\theta_{0}/2}\bigr)\bigl(1+|y|^{\theta_{0}}\bigr)\bigl(d+|y|^{2}\bigr)(2\pi\tau)^{-d/2}e^{-|y|^{2}/2},
\end{equation*}
and we conclude as before, with the bounded function $y\mapsto(1+|y|^{\theta_{0}})(d+|y|^{2})e^{-|y|^{2}/4}$.
\end{proof}

\begin{proof}[Proof of Lemma~\ref{lem:weights}]
A direct computation on the expression \eqref{eq: weight1} gives
\begin{equation*}
\nabla\hat{w}_{x_{0}}^{\theta_{0}}(x)=-\theta_{0}\,(x-x_{0})\bigl(1+|x-x_{0}|^{2}\bigr)^{-\frac{\theta_{0}}{2}-1},
\qquad x\in\R^{d},
\end{equation*}
and since $|x-x_{0}|\leq 1+|x-x_{0}|^{2}$, the first relation in \eqref{eq:what-grad} follows. The second relation is the triangle inequality in $\R^{d+1}$, which reads $(1+|y-x_{0}|^{2})^{1/2}\leq(1+|x-x_{0}|^{2})^{1/2}+|x-y|$, divided by $(1+|x-x_{0}|^{2})^{1/2}\geq 1$ and raised to the power $\theta_{0}$.
\end{proof}

\section{Moment estimates for stochastic convolutions}\label{app:conv}

This section contains the proofs of the spectral bound of Lemma~\ref{lem: spectral} and of the weighted moment bounds of Lemma~\ref{lem:weighted-moments}, on which the weighted estimate of Proposition~\ref{prop:weighted-conv} rests. We recall that the stochastic convolution $\mathcal{Z}^{h}$ is defined by \eqref{eq:StochConv} and that the quantity $\mathcal{N}^{x_{0}}_{p,\theta_{0},T}(h)$ is given by \eqref{eq:Nfunctional}.

\begin{proof}[Proof of Lemma~\ref{lem: spectral}]
By \eqref{eq: heatkernel}, $G(\tau,\cdot)^{\hat\alpha}$ is proportional to $G(\tau/\hat\alpha,\cdot)$, namely 
\[G(\tau,x)^{\hat\alpha}=(2\pi\tau)^{-d(\hat\alpha-1)/2}\hat\alpha^{-d/2}\,G(\tau/\hat\alpha,x),\]
so that formula~\eqref{eq: FourierPow} follows from the Fourier transform \eqref{eq:G-fourier} of the heat kernel. Squaring~\eqref{eq: FourierPow} and integrating against $\specm$:
\begin{equation}\label{eq: spectralBound-step}
    \int_{\R^{d}}\bigl|\mathcal{F}(G(\tau,\cdot)^{\hat\alpha})(\xi)\bigr|^{2}\specm(d\xi) = C_{\hat\alpha,d}^{2}\,\tau^{d(1-\hat\alpha)}\int_{\R^{d}}e^{-|\xi|^{2}\tau/\hat\alpha}\,\specm(d\xi).
\end{equation}
We split the remaining spectral integral into low and high frequencies:
\begin{equation}\label{eq: spectral-split}
    \int_{\R^{d}}e^{-|\xi|^{2}\tau/\hat\alpha}\specm(d\xi) \leq \specm(\{|\xi|\leq 1\})+\int_{\{|\xi|>1\}}\bigl((\tau/\hat\alpha)^{1/2}|\xi|\bigr)^{-2\vk}\bigl((\tau/\hat\alpha)^{1/2}|\xi|\bigr)^{2\vk}e^{-|\xi|^{2}\tau/\hat\alpha}\specm(d\xi).
\end{equation}
Using the elementary bound $\sup_{r\geq 0}r^{2\vk}e^{-r^{2}}\leq C_{\vk}<\infty$ with $r=(\tau/\hat\alpha)^{1/2}|\xi|$, and $\frac{1+|\xi|^{2\vk}}{|\xi|^{2\vk}}\leq 2$ for $|\xi|>1$, plus the fact that $\vk>1-\eta$, we end up with
\begin{equation}\label{eq: spectral-high}
    \int_{\R^{d}}e^{-|\xi|^{2}\tau/\hat\alpha}\specm(d\xi)\leq \specm(\{|\xi|\leq 1\})+C_{\vk}(\tau/\hat\alpha)^{-\vk}\int_{\R^{d}}\frac{\specm(d\xi)}{1+|\xi|^{2(1-\eta)}}\leq C_{\hat\alpha,\specm,\eta,\vk}\bigl(1+\tau^{-\vk}\bigr),
\end{equation}
where the second bound is ensured via the strong Dalang condition~\eqref{eq:strongDalang}. Inserting~\eqref{eq: spectral-high} into~\eqref{eq: spectralBound-step} yields~\eqref{eq: spectralBound}.
\end{proof}

\begin{proof}[Proof of Lemma~\ref{lem:weighted-moments}]
Throughout this proof, the value of the constant $C$ may change from line to line, but always depends only on the parameters explicitly indicated. We first derive a master inequality, valid for every deterministic kernel, and then apply it to the three kernels corresponding to items \ref{it:wm-point}, \ref{it:wm-space} and \ref{it:wm-time}.

\smallskip

\noindent
{\it Step 1: a master inequality.} Let $k\colon[0,T]\times\R^{d}\to\R$ be a deterministic measurable function, and consider the Walsh integral $M:=\int_{0}^{T}\int_{\R^{d}}k(s,y)h(s,y)\,W(ds,dy)$. Its quadratic variation is
\begin{equation}\label{eq:qv}
\langle M\rangle=\int_{0}^{T}\int_{\R^{d}}\int_{\R^{d}}k(s,y)h(s,y)\,\ell(y-z)\,h(s,z)k(s,z)\,dy\,dz\,ds,
\end{equation}
where we recall that $\ell$ is the covariance density of Assumption~\ref{assump:noise}. By the Burkholder-Davis-Gundy inequality, $\mathbf{E}[|M|^{p}]\leq c_{p}\,\mathbf{E}[\langle M\rangle^{p/2}]$. Since $\ell\geq 0$, we can bound $\langle M\rangle$ by the same integral with the absolute values $|k(s,y)h(s,y)|$ and $|h(s,z)k(s,z)|$ in place of the products. Then we apply Minkowski's inequality in $L^{p/2}(\Omega)$ (recall that $p\geq 2$), followed by the Cauchy-Schwarz inequality
\begin{equation*}
\|h(s,y)h(s,z)\|_{L^{p/2}(\Omega)}\leq\|h(s,y)\|_{L^{p}(\Omega)}\,\|h(s,z)\|_{L^{p}(\Omega)} .
\end{equation*}
This gives
\begin{equation}\label{eq:master-1}
\|M\|_{L^{p}(\Omega)}^{2}\leq c_{p}\int_{0}^{T}\int_{\R^{d}}\int_{\R^{d}}|k(s,y)|\,\|h(s,y)\|_{L^{p}(\Omega)}\,\ell(y-z)\,\|h(s,z)\|_{L^{p}(\Omega)}\,|k(s,z)|\,dy\,dz\,ds.
\end{equation}
Now by definition \eqref{eq:Nfunctional} of $N$ and by \eqref{eq:weight-triangle}, for every $x\in\R^{d}$ we have
\begin{equation}\label{eq:h-weighted}
\|h(s,y)\|_{L^{p}(\Omega)}\leq N\bigl(1+|y-x_{0}|^{\theta_{0}}\bigr)\leq c_{\theta_{0}}\,N\bigl(1+|x-x_{0}|^{\theta_{0}}\bigr)\bigl(1+|x-y|^{\theta_{0}}\bigr).
\end{equation}
Plugging \eqref{eq:h-weighted} into \eqref{eq:master-1}, we end up with our master inequality
\begin{equation}\label{eq:master}
\|M\|_{L^{p}(\Omega)}^{2}\leq C\,N^{2}\bigl(1+|x-x_{0}|^{\theta_{0}}\bigr)^{2}
\int_{0}^{T}\int_{\R^{d}}\int_{\R^{d}}\tilde{k}_{x}(s,y)\,\ell(y-z)\,\tilde{k}_{x}(s,z)\,dy\,dz\,ds,
\end{equation}
where $C=c_{p}c_{\theta_{0}}^{2}$ and where we have set
\begin{equation}\label{eq:ktilde}
\tilde{k}_{x}(s,y):=|k(s,y)|\bigl(1+|x-y|^{\theta_{0}}\bigr).
\end{equation}
Relation \eqref{eq:master} holds for every $x\in\R^{d}$, and the point $x$ will be chosen according to the kernel $k$. In each of the three steps below, the kernel $\tilde{k}_{x}$ is bounded thanks to Lemma~\ref{lem:weight-absorb} by heat kernels with a doubled variance, and the $dy\,dz$ integral in \eqref{eq:master} is then computed in Fourier modes thanks to relation \eqref{eq: relation}.

\smallskip

\noindent
{\it Step 2: the pointwise bound.} Fix $(t,x)\in[0,T]\times\R^{d}$ and take $k(s,y):=G(t-s,x-y)\mathbf{1}_{[0,t]}(s)$, so that $M=Z(t,x)$. By \eqref{eq:absorb-power} with $a=1$ and $\tau=t-s$, the kernel \eqref{eq:ktilde} satisfies, for any $s\in[0,t]$:
\begin{equation}\label{eq:ktilde-point}
\tilde{k}_{x}(s,y)\leq C\bigl(1+(t-s)^{\theta_{0}/2}\bigr)G(2(t-s),x-y)\leq C\bigl(1+T^{\theta_{0}/2}\bigr)G(2(t-s),x-y).
\end{equation}
By relation \eqref{eq: relation} and formula \eqref{eq: FourierPow} with $\hat{\alpha}=1$, the modulus of the Fourier transform of $G(2\tau,x-\cdot)$ being the one of $G(2\tau,\cdot)$, we have for every $\tau>0$
\begin{equation}\label{eq:fourier-double}
\int_{\R^{d}}\int_{\R^{d}}G(2\tau,x-y)\,\ell(y-z)\,G(2\tau,x-z)\,dy\,dz=\int_{\R^{d}}e^{-2\tau|\xi|^{2}}\,\specm(d\xi)\leq C\bigl(1+\tau^{-(1-\chik)}\bigr),
\end{equation}
where the last inequality is \eqref{eq:spectral-one}, applied with $2\tau$ in place of $\tau$ and with the exponent $1-\chik\in(1-\eta,1)$. Inserting \eqref{eq:ktilde-point} and \eqref{eq:fourier-double} into \eqref{eq:master}, we get
\begin{equation}\label{eq:point-1}
\|Z(t,x)\|_{L^{p}(\Omega)}^{2}\leq C\,N^{2}\bigl(1+|x-x_{0}|^{\theta_{0}}\bigr)^{2}\bigl(1+T^{\theta_{0}/2}\bigr)^{2}\int_{0}^{t}\bigl(1+(t-s)^{-(1-\chik)}\bigr)\,ds .
\end{equation}
Computing the time integral in \eqref{eq:point-1}, it is now readily checked that \eqref{eq:wm-point} is satisfied.

\smallskip

\noindent
{\it Step 3: the space increments.} Fix $t\in[0,T]$ and $x_{1},x_{2}\in\R^{d}$, and take $k(s,y):=[G(t-s,x_{1}-y)-G(t-s,x_{2}-y)]\mathbf{1}_{[0,t]}(s)$, so that $M=Z(t,x_{1})-Z(t,x_{2})$. We first bound the kernel. Namely, writing $\tau=t-s$, the sub-additivity of $r\mapsto r^{1-\lambda}$ gives
\begin{equation}\label{eq:space-holder}
|k(s,y)|\leq|k(s,y)|^{\lambda}\Bigl\{G(\tau,x_{1}-y)^{1-\lambda}+G(\tau,x_{2}-y)^{1-\lambda}\Bigr\},
\end{equation}
and the mean value theorem together with the gradient bound \eqref{eq: gradG} yields
\begin{equation}\label{eq:space-mvt}
|k(s,y)|^{\lambda}\leq C_{d}^{\lambda}\,\tau^{-\frac{\lambda(d+1)}{2}}\,|x_{1}-x_{2}|^{\lambda}.
\end{equation}
We now bring in the weight. Here the point $x$ of \eqref{eq:master} is taken to be $x_{1}$ for the first term in \eqref{eq:space-holder} and $x_{2}$ for the second one. More precisely, going back to \eqref{eq:h-weighted} and applying \eqref{eq:weight-triangle} with $x=x_{i}$ for the term involving $G(\tau,x_{i}-y)$, and then \eqref{eq:absorb-power} with $a=1-\lambda$, we get for $i=1,2$
\begin{equation}\label{eq:ktilde-space-i}
\|h(s,y)\|_{L^{p}(\Omega)}\,G(\tau,x_{i}-y)^{1-\lambda}\leq C\,N\bigl(1+|x_{i}-x_{0}|^{\theta_{0}}\bigr)\bigl(1+T^{\theta_{0}/2}\bigr)\,G(2\tau,x_{i}-y)^{1-\lambda}.
\end{equation}
Gathering \eqref{eq:space-holder}, \eqref{eq:space-mvt} and \eqref{eq:ktilde-space-i}, we end up with
\begin{align}\label{eq:ktilde-space}
|k(s,y)|\,\|h(s,y)\|_{L^{p}(\Omega)}\leq C\,N&\bigl(1+|x_{1}-x_{0}|^{\theta_{0}}+|x_{2}-x_{0}|^{\theta_{0}}\bigr)\bigl(1+T^{\theta_{0}/2}\bigr)\tau^{-\frac{\lambda(d+1)}{2}}|x_{1}-x_{2}|^{\lambda}\notag\\
&\times\Bigl\{G(2\tau,x_{1}-y)^{1-\lambda}+G(2\tau,x_{2}-y)^{1-\lambda}\Bigr\}.
\end{align}
We can now insert \eqref{eq:ktilde-space} into \eqref{eq:master-1}. The $dy\,dz$ integral is then handled in Fourier modes. Namely, by \eqref{eq: relation} and \eqref{eq: FourierPow} with $\hat{\alpha}=1-\lambda$, the Fourier transform of $G(2\tau,x_{i}-\cdot)^{1-\lambda}$ is $e^{i\xi\cdot x_{i}}$ times the one of $G(2\tau,\cdot)^{1-\lambda}$, and $|e^{i\xi\cdot x_{1}}+e^{i\xi\cdot x_{2}}|^{2}\leq 4$. Hence
\begin{align}\label{eq:fourier-space}
\int_{\R^{d}}\int_{\R^{d}}\sum_{i,j=1}^{2}G(2\tau,x_{i}-y)^{1-\lambda}\ell(y-z)\,G(2\tau,x_{j}-z)^{1-\lambda}\,dy\,dz
&\leq 4\!\int_{\R^{d}}\bigl|\mathcal{F}\bigl(G(2\tau,\cdot)^{1-\lambda}\bigr)(\xi)\bigr|^{2}\specm(d\xi)\notag\\
&\leq C\,\tau^{d\lambda}\bigl(1+\tau^{-(1-\chik)}\bigr),
\end{align}
where the last bound is \eqref{eq: spectralBound}, applied with $\hat{\alpha}=1-\lambda$, with $2\tau$ in place of $\tau$ and with the exponent $1-\chik$. Gathering \eqref{eq:ktilde-space} and \eqref{eq:fourier-space}, and observing that $\tau^{-\lambda(d+1)}\tau^{d\lambda}=\tau^{-\lambda}$, we obtain
\begin{align}\label{eq:space-1}
\|Z(t,x_{1})-Z(t,x_{2})\|_{L^{p}(\Omega)}^{2}\leq C\,N^{2}&\bigl(1+|x_{1}-x_{0}|^{\theta_{0}}+|x_{2}-x_{0}|^{\theta_{0}}\bigr)^{2}\bigl(1+T^{\theta_{0}}\bigr)\notag\\
&\times|x_{1}-x_{2}|^{2\lambda}\int_{0}^{t}\tau^{-\lambda}\bigl(1+\tau^{-(1-\chik)}\bigr)\,d\tau .
\end{align}
Since $\lambda<\chik$, the time integral in \eqref{eq:space-1} is bounded by $C(T^{1-\lambda}+T^{\chik-\lambda})\leq C(1+T)T^{\chik-\lambda}$, which gives \eqref{eq:wm-space}.

\smallskip

\noindent
{\it Step 4: the time increments.} Fix $x\in\R^{d}$ and $0\leq t_{1}<t_{2}\leq T$, and set $\delta:=t_{2}-t_{1}$. We decompose $Z(t_{2},x)-Z(t_{1},x)=I_{1}+I_{2}$, where $I_{1}$ and $I_{2}$ are the Walsh integrals of $h$ against the kernels
\begin{equation}\label{eq:time-kernels}
k_{1}(s,y):=G(t_{2}-s,x-y)\mathbf{1}_{[t_{1},t_{2}]}(s),
\qquad
k_{2}(s,y):=\bigl[G(t_{2}-s,x-y)-G(t_{1}-s,x-y)\bigr]\mathbf{1}_{[0,t_{1}]}(s).
\end{equation}
The term $I_{1}$ is handled exactly as in Step~2, the integration in time being now over $[t_{1},t_{2}]$. Namely, \eqref{eq:master}, \eqref{eq:ktilde-point} and \eqref{eq:fourier-double} give
\begin{equation}\label{eq:time-I1}
\|I_{1}\|_{L^{p}(\Omega)}^{2}\leq C\,N^{2}\bigl(1+|x-x_{0}|^{\theta_{0}}\bigr)^{2}\bigl(1+T^{\theta_{0}}\bigr)\int_{0}^{\delta}\bigl(1+r^{-(1-\chik)}\bigr)\,dr
\leq C\,N^{2}\bigl(1+|x-x_{0}|^{\theta_{0}}\bigr)^{2}c_{T}\,\delta^{\chik},
\end{equation}
where we let the patient reader check the bound on the time integral, which is the one of Step~2.

For the term $I_{2}$, observe that since $\partial_{\tau}G=\frac{1}{2}\Delta G$, for $s\in[0,t_{1})$ we have
\begin{equation}\label{eq:time-laplace}
k_{2}(s,y)=\frac{1}{2}\int_{t_{1}-s}^{t_{2}-s}\Delta_{z}G(\tau,x-y)\,d\tau,
\end{equation}
so that by \eqref{eq:absorb-laplace} the weighted kernel \eqref{eq:ktilde} satisfies
\begin{equation}\label{eq:ktilde-time}
\tilde{k}_{2,x}(s,y)\leq C\bigl(1+T^{\theta_{0}/2}\bigr)\int_{t_{1}-s}^{t_{2}-s}\frac{1}{\tau}\,G(2\tau,x-y)\,d\tau .
\end{equation}
We insert \eqref{eq:ktilde-time} into \eqref{eq:master} and apply Fubini's theorem, all the integrands being nonnegative. By \eqref{eq: relation} and \eqref{eq: FourierPow} with $\hat{\alpha}=1$, for $\tau,\tau'>0$ we have
\begin{align}\label{eq:fourier-two-times}
\int_{\R^{d}}\int_{\R^{d}}G(2\tau,x-y)\,\ell(y-z)\,G(2\tau',x-z)\,dy\,dz
&=\int_{\R^{d}}e^{-(\tau+\tau')|\xi|^{2}}\,\specm(d\xi)\notag\\
&\leq C\bigl(1+(\tau+\tau')^{-(1-\chik)}\bigr)\notag\\
&\leq C\bigl(1+(t_{1}-s)^{-(1-\chik)}\bigr),
\end{align}
where we have used \eqref{eq:spectral-one} and $\tau,\tau'\geq t_{1}-s$. Therefore
\begin{equation}\label{eq:time-I2-a}
\|I_{2}\|_{L^{p}(\Omega)}^{2}\leq C\,N^{2}\bigl(1+|x-x_{0}|^{\theta_{0}}\bigr)^{2}\bigl(1+T^{\theta_{0}}\bigr)
\int_{0}^{t_{1}}\Bigl(\int_{t_{1}-s}^{t_{2}-s}\frac{d\tau}{\tau}\Bigr)^{2}\bigl(1+(t_{1}-s)^{-(1-\chik)}\bigr)\,ds .
\end{equation}
The inner integral in \eqref{eq:time-I2-a} is $\log(1+\delta/r)$ with $r=t_{1}-s$, and we claim that
\begin{equation}\label{eq:log-integral}
\int_{0}^{t_{1}}\log^{2}\Bigl(1+\frac{\delta}{r}\Bigr)\bigl(1+r^{-(1-\chik)}\bigr)\,dr\leq C_{\chik}(1+T)\,\delta^{\chik}.
\end{equation}
Indeed, we split the integral at $r=\delta$ (if $\delta<t_{1}$, and otherwise only the first region is present). For $r\leq\delta$ we bound $\log(1+\delta/r)\leq\log(2\delta/r)$, and the change of variables $r=\delta v$ gives
\begin{equation}\label{eq:log-small}
\int_{0}^{\delta}\log^{2}\Bigl(\frac{2\delta}{r}\Bigr)\bigl(1+r^{-(1-\chik)}\bigr)\,dr=\delta\int_{0}^{1}\log^{2}\Bigl(\frac{2}{v}\Bigr)dv+\delta^{\chik}\int_{0}^{1}\log^{2}\Bigl(\frac{2}{v}\Bigr)v^{-(1-\chik)}\,dv\leq C_{\chik}\bigl(\delta+\delta^{\chik}\bigr),
\end{equation}
both integrals being finite since $\chik>0$. For $r\geq\delta$ we bound $\log(1+\delta/r)\leq\delta/r$, so that
\begin{equation}\label{eq:log-large}
\int_{\delta}^{t_{1}}\frac{\delta^{2}}{r^{2}}\bigl(1+r^{-(1-\chik)}\bigr)\,dr\leq\delta^{2}\Bigl(\frac{1}{\delta}+\frac{\delta^{-2+\chik}}{2-\chik}\Bigr)\leq\delta+\delta^{\chik}.
\end{equation}
Summing \eqref{eq:log-small} and \eqref{eq:log-large} and using $\delta\leq T^{1-\chik}\delta^{\chik}$, we get \eqref{eq:log-integral}. Plugging \eqref{eq:log-integral} into \eqref{eq:time-I2-a} and recalling that $(1+T^{\theta_{0}})(1+T)\leq 2c_{T}$, we obtain
\begin{equation}\label{eq:time-I2}
\|I_{2}\|_{L^{p}(\Omega)}^{2}\leq C\,N^{2}\bigl(1+|x-x_{0}|^{\theta_{0}}\bigr)^{2}c_{T}\,\delta^{\chik}.
\end{equation}
Finally, \eqref{eq:wm-time} follows from \eqref{eq:time-I1}, \eqref{eq:time-I2} and the triangle inequality in $L^{p}(\Omega)$.
\end{proof}

\section{Factorization and the penalized stochastic convolution}\label{app:factorization}

This section contains the proofs of the preliminary results of Section~\ref{sec:bounds}. We first prove the factorization formula of Lemma~\ref{lem: fact}, the integral bound of Lemma~\ref{lem:int-exp} and the moment bound on $Y_{\alpha}^{\beta}$ of Lemma~\ref{lem:Ybeta-moment}. We then verify, in Lemmas~\ref{lem: SbetaBound}, \ref{lem: SpaceIncrem} and~\ref{lem: TimeIncrem}, the three hypotheses of Proposition~\ref{prop: Salins} for the penalized stochastic convolution $\mathcal{S}^{\beta}u$ of \eqref{eq: Sbeta-recall}. Throughout, $\specm$ is the spectral measure of Assumption~\ref{assump:noise} and $M_{\sigma}$ the bound on $\sigma$ of Assumption~\ref{assump:diff}.

\begin{proof}[Proof of Lemma~\ref{lem: fact}]
Let $Z(t,x)$ be the double integral in the right hand side of $(\ref{eq:factorization})$, that is the right hand side without its prefactor $\sin(\alpha\pi)/\pi$. Replacing $Y_{\alpha}^{\beta}$ by its expression $(\ref{eq:Yprocess})$, we get
\begin{multline*}
    Z(t,x)=\int_{0}^{t}\int_{\R^d}\int_{0}^{s}\int_{\R^d}e^{-\beta(t-s)}G(t-s,x-y)(t-s)^{\alpha-1}\\\hspace{4cm}\times e^{-\beta(s-r)}G(s-r,y-z)\sigma(u(r,z))(s-r)^{-\alpha}W(dr,dz) dy\,ds.
    \end{multline*}
      Applying a stochastic version of Fubini's theorem, we integrate in the $y$ variable first and invoke the semigroup property \eqref{eq:G-semigroup} of the heat kernel $G$. This yields \begin{equation*}Z(t,x)=\int_{0}^{t}\int_{0}^{s}\int_{\R^d}e^{-\beta(t-r)}G(t-r,x-z)\sigma(u(r,z))W(dr,dz) (t-s)^{\alpha-1}(s-r)^{-\alpha}ds.
    \end{equation*}
    Note that we are not checking the integrability conditions for the application of the stochastic Fubini theorem. They rest on the moment bound of Lemma~\ref{lem:Ybeta-moment} on $Y_{\alpha}^{\beta}$, which is where the condition $\alpha<\frac{\eta}{2}$ enters, and we leave them to the reader. Next we resort again to Fubini's theorem, and we now integrate with respect to the $s$ variable to arrive at
    \begin{equation}\label{eq: fac1}
Z(t,x)=\int_{0}^{t}\int_{\R^d}\left(\int_{r}^{t}(t-s)^{\alpha-1}(s-r)^{-\alpha}ds\right)e^{-\beta(t-r)}G(t-r,x-z)\sigma(u(r,z))W(dr,dz),
\end{equation}
where the inner integral is a Beta integral. Namely, the change of variables $s=r+(t-r)v$ gives
\[\int_{r}^{t}(t-s)^{\alpha-1}(s-r)^{-\alpha}\,ds=\int_{0}^{1}(1-v)^{\alpha-1}v^{-\alpha}\,dv=\Gamma(\alpha)\Gamma(1-\alpha)=\frac{\pi}{\sin(\pi \alpha)},\]
by Euler's reflection formula.
Plugging this in $(\ref{eq: fac1})$, we end up with
\begin{equation*}
    Z(t,x)=\frac{\pi}{\sin(\alpha \pi)}\int_{0}^{t}\int_{\R^d}e^{-\beta(t-r)}G(t-r,x-z)\sigma(u(r,z))W(dr,dz)=\frac{\pi}{\sin(\alpha \pi)}\mathcal{S}^{\beta}u(t,x),
\end{equation*}
where the second identity is our definition $(\ref{eq:Sbeta})$. Multiplying by $\sin(\alpha\pi)/\pi$, which is the prefactor in \eqref{eq:factorization}, finishes the proof of our claim.
\end{proof}

\begin{proof}[Proof of Lemma~\ref{lem:int-exp}]
Since $u\le t$ for $u\in(0,t]$, we have $u^{\ze}\le t^{\ze}$, and hence
\begin{equation}\label{eq:int-exp-step1}
u^{-\aK} =  u^{\ze}   u^{-\aK-\ze}\;\le\; t^{\ze}\,u^{-\aK-\ze},
\qquad u\in(0,t].
\end{equation}
Applying~\eqref{eq:int-exp-step1} and then extending the domain of integration to $(0,\infty)$, we obtain
\begin{equation}\label{eq:int-exp-step2}
\int_{0}^{t}e^{-\beta u}\,u^{-\aK}\,du
\;\le\; t^{\ze}\int_{0}^{\infty}e^{-\beta u}\,u^{-\aK-\ze}\,du.
\end{equation}
Note that the condition $\aK+\ze<1$ ensures that the integral in the right hand side of~\eqref{eq:int-exp-step2} is finite. Next setting $v=\beta u$ in this integral gives
\begin{equation}\label{eq:int-exp-gamma}
\int_{0}^{\infty}e^{-\beta u}\,u^{-\aK-\ze}\,du
= \frac{1}{\beta^{1-\aK-\ze}}\int_{0}^{\infty}e^{-v}\,v^{-\aK-\ze}\,dv
= \frac{\Gamma(1-\aK-\ze)}{\beta^{1-\aK-\ze}}.
\end{equation}
Inserting~\eqref{eq:int-exp-gamma} into~\eqref{eq:int-exp-step2} yields~\eqref{eq:int-exp} and concludes our proof.
\end{proof}

\begin{proof}[Proof of Lemma~\ref{lem:Ybeta-moment}]
Throughout this proof, the value of the constant $C$ may change from line to line, but always depends only on the parameters explicitly indicated. Fix $s\in[0,T]$ and $y\in\R^d$. Applying the BDG inequality to the stochastic integral \eqref{eq:Yprocess} defining $Y_{\alpha}^{\beta}$, and computing its quadratic variation thanks to the covariance \eqref{eq:covWiener} of the noise, we get
\begin{align}\label{eq:Ybeta-BDG}
\mathbf{E}\big[|Y_{\alpha}^{\beta}(s,y)|^{p}\big]&\leq c_{p}\,\mathbf{E}\Bigg[\Bigg(\int_{0}^{s}\int_{\R^d}\int_{\R^d}(s-r)^{-2\alpha}e^{-2\beta(s-r)}G(s-r,y-z)\sigma(u(r,z))\nonumber \\
& \hspace{3.5cm}\times \ell(z-w)\,\sigma(u(r,w))\,G(s-r,y-w)\,dz\,dw\,dr\Bigg)^{p/2}\Bigg].
\end{align}
Since $|\sigma|\leq M_{\sigma}$ by Assumption~\ref{assump:diff}, the integral in \eqref{eq:Ybeta-BDG} is bounded by the deterministic quantity $M_{\sigma}^{2}A_{s,\beta}$, where
\begin{equation}\label{eq: Aest}
A_{s,\beta}:=\int_{0}^{s}\int_{\R^d}\int_{\R^d}(s-r)^{-2\alpha}e^{-2\beta(s-r)}G(s-r,y-z)\ell(z-w)G(s-r,y-w)\,dz\,dw\,dr.
\end{equation}
Hence \eqref{eq:Ybeta-BDG} reads
\begin{equation}\label{eq:Ybeta-BDG2}
\mathbf{E}\big[|Y_{\alpha}^{\beta}(s,y)|^{p}\big]\leq c_{p}\,M_{\sigma}^{p}\,A_{s,\beta}^{p/2}.
\end{equation}
We now implement the relation $(\ref{eq: relation})$ of Assumption~\ref{assump:noise}, which turns the space integral in \eqref{eq: Aest} into a spectral one. This yields
\begin{equation}\label{eq: A}
A_{s,\beta}=\int_{0}^{s}(s-r)^{-2\alpha}e^{-2\beta(s-r)}\int_{\R^d}e^{-(s-r)|\xi|^{2}}\specm(d\xi)\,dr .
\end{equation}
Next we apply the spectral bound \eqref{eq:spectral-one}, that is Lemma~\ref{lem: spectral} with $\hat{\alpha}=1$, $\tau=s-r$ and the exponent $\hat\vk\in(1-\eta,1)$. We obtain
\begin{equation}\label{eq: boundA}
A_{s,\beta}\leq C_{\specm,\eta,\hat\vk}\left(\int_{0}^{s}(s-r)^{-2\alpha}e^{-2\beta(s-r)}dr +\int_{0}^{s}(s-r)^{-2\alpha-\hat\vk}e^{-2\beta(s-r)}dr \right).
\end{equation}
Setting $u=s-r$ in both integrals of $(\ref{eq: boundA})$ and applying Lemma~\ref{lem:int-exp} with $\beta$ replaced by $2\beta$ and $\ze=0$ (with $\aK=2\alpha$ for the first integral, then with $\aK=2\alpha+\hat{\vk}<1$ for the second integral), we obtain
\begin{equation}\label{eq: integrals0}
\int_{0}^{s}(s-r)^{-2\alpha}e^{-2\beta(s-r)}\,dr \leq \frac{\Gamma(1-2\alpha)}{(2\beta)^{1-2\alpha}}, \quad
\int_{0}^{s}(s-r)^{-2\alpha-\hat\vk}e^{-2\beta(s-r)}\,dr \leq \frac{\Gamma(1-2\alpha-\hat\vk)}{(2\beta)^{1-2\alpha-\hat\vk}}.
\end{equation}
Inserting $(\ref{eq: integrals0})$ into $(\ref{eq: boundA})$ and absorbing $\Gamma(1-2\alpha)$ and $\Gamma(1-2\alpha-\hat{\vk})$ into the constant $C_{\specm,\eta,\hat\vk}$, we conclude
\begin{equation}\label{eq: boundAfinal}
A_{s,\beta} \leq C_{\specm,\eta,\hat\vk}\left(\frac{1}{\beta^{1-2\alpha}}+\frac{1}{\beta^{1-2\alpha-\hat\vk}}\right).
\end{equation}
Finally, inserting $(\ref{eq: boundAfinal})$ into $(\ref{eq:Ybeta-BDG2})$ we obtain our claim $(\ref{eq:Ybeta-moment})$. The constant does not depend on $T$, since neither \eqref{eq:Ybeta-BDG2} nor \eqref{eq: boundAfinal} does.
\end{proof}

\begin{proof}[Proof of Lemma~\ref{lem: SbetaBound}]
    Throughout this proof, the value of the constant $C$ may change from line to line, but always depends only on the parameters explicitly indicated. Fix $p\geq 1$, $t\geq 0$ and $x\in\R^d$. Applying the BDG inequality to the mild form $(\ref{eq: Sbeta-recall})$ of $\mathcal{S}^{\beta}u$, and computing the quadratic variation thanks to the covariance \eqref{eq:covWiener} of the noise, we get
\begin{align}\label{eq: BDG}
       & \mathbf{E}[|\mathcal{S}^{\beta}u(t,x)|^{p}] = \mathbf{E}\Bigg[\left|\int_{0}^{t}\int_{\R^d}e^{-\beta(t-s)}G(t-s,x-y)\sigma(u(s,y))W(ds,dy) \right|^{p} \Bigg] \nonumber \\
        &\leq c_{p}\,\mathbf{E}\Bigg[\Bigg|\int_{0}^{t}\int_{\R^d}\int_{\R^d}e^{-2\beta(t-s)}G(t-s,x-y)\sigma(u(s,y)) \nonumber \\
        &\hspace{4cm} \times \ell(y-z)\sigma(u(s,z))G(t-s,x-z)\,dy\,dz\,ds \Bigg|^{p/2}\Bigg].
\end{align}
Since $|\sigma|\leq M_{\sigma}$ by Assumption~\ref{assump:diff}, the integral in \eqref{eq: BDG} is bounded by the deterministic quantity $M_{\sigma}^{2}J_{t,\beta}$, where
\begin{equation}\label{eq: Jtbeta}J_{t,\beta}:=\int_{0}^{t}\int_{\R^d}\int_{\R^d}e^{-2\beta(t-s)}G(t-s,x-y)\ell(y-z)G(t-s,x-z)ds\,dy\,dz .\end{equation}
Hence \eqref{eq: BDG} reads
\begin{equation}\label{eq: boundSbeta}
        \mathbf{E}[|\mathcal{S}^{\beta}u(t,x)|^{p}] \leq c_{p}\,M_{\sigma}^{p}\,J_{t,\beta}^{p/2}.
\end{equation}
    With \eqref{eq: boundSbeta} in hand, in order to prove our claim~\eqref{eq: prelimLemmaBound1} it remains to bound $J_{t,\beta}$. Implementing the relationship $(\ref{eq: relation})$ from Assumption~\ref{assump:noise}, we have
 \begin{align}
        J_{t,\beta}
        &=\int_{0}^{t}\int_{\R^d}e^{-2\beta(t-s)}|\mathcal{F}(G(t-s,x-\cdot))(\xi)|^2 \specm(d\xi)\,ds \nonumber \\
        &=\int_{0}^{t}\int_{\R^d}e^{-2\beta(t-s)}e^{-(t-s)|\xi|^2}\specm(d\xi)\,ds \qquad\text{(by \eqref{eq:G-fourier})}
         \label{eq: J}=\int_{0}^{t}e^{-2\beta(t-s)}\widetilde{J}^{\specm}_{s,t}\,ds ,
    \end{align}
    where we have set
    \begin{equation}\label{eq: Jtilde}
        \widetilde{J}^{\specm}_{s,t}:=\int_{\R^{d}}\bigl|\mathcal{F}(G(t-s,\cdot))(\xi)\bigr|^{2}\specm(d\xi)=\int_{\R^{d}}e^{-(t-s)|\xi|^{2}}\specm(d\xi).
    \end{equation}
    Applying Lemma~\ref{lem: spectral} with $\hat\alpha=1$ and $\tau=t-s$, we obtain
    \begin{equation}\label{eq: boundJtilde}
        \widetilde{J}^{\specm}_{s,t}\leq C_{\specm,\eta,\vk}\bigl(1+(t-s)^{-\vk}\bigr),\qquad \vk\in(1-\eta,1).
    \end{equation}
    Plugging $(\ref{eq: boundJtilde})$ into~(\ref{eq: J}), we thus get
    \begin{align}\label{eq: boundJ}
        J_{t,\beta}&\leq {C_{\specm,\eta,\vk}}\left(\int_{0}^{t}e^{-2\beta(t-s)}\,ds +\int_{0}^{t}e^{-2\beta(t-s)}(t-s)^{-\vk}\,ds\right)
    \end{align}
    Setting $u=t-s$ in both integrals of $(\ref{eq: boundJ})$ and applying Lemma~\ref{lem:int-exp} with $\beta$ replaced by $2\beta$ and $\ze=0$ (with $\aK=0$ for the first integral, then with $\aK=\vk\in(1-\eta,1)$ for the second integral), we obtain
    \begin{equation}\label{eq: integrals}
        \int_{0}^{t}e^{-2\beta(t-s)}\,ds \leq \frac{1}{2\beta}, \qquad
        \int_{0}^{t}e^{-2\beta(t-s)}(t-s)^{-\vk}\,ds \leq \frac{\Gamma(1-\vk)}{(2\beta)^{1-\vk}}.
    \end{equation}
    Inserting $(\ref{eq: integrals})$ into $(\ref{eq: boundJ})$ and absorbing $\Gamma(1-\vk)$ into the constant $C_{\specm,\eta,\vk}$, we conclude
    \begin{equation}\label{eq: boundJfinal}
        J_{t,\beta} \leq C_{\specm,\eta,\vk}\left(\frac{1}{\beta}+\frac{1}{\beta^{1-\vk}}\right).
    \end{equation}
    Finally, inserting $(\ref{eq: boundJfinal})$ into $(\ref{eq: boundSbeta})$ we obtain our claim $(\ref{eq: prelimLemmaBound1})$.
\end{proof}

\begin{proof}[Proof of Lemma~\ref{lem: SpaceIncrem}]
Fix $p\geq 1$, $t\geq 0$ and $x,y\in \R^d$. The same argument as in $(\ref{eq: BDG})$--$(\ref{eq: boundSbeta})$, with the kernel $G(t-s,x-\cdot)$ replaced by the spatial difference $G(t-s,x-\cdot)-G(t-s,y-\cdot)$, gives
    \begin{equation}\label{eq: boundSbetaSpace}
       \mathbf{E}[|\mathcal{S}^{\beta}u(t,x)-\mathcal{S}^{\beta}u(t,y)|^{p}] \leq c_{p}\,M_{\sigma}^{p}\, K_{t,\beta}(x,y)^{p/2},
    \end{equation}
    where $K_{t,\beta}(x,y)$ is the spatial-increment analogue of $J_{t,\beta}$ defined in $(\ref{eq: Jtbeta})$, given explicitly by
    \begin{align}\label{eq: K}
        K_{t,\beta}(x,y)&:=\int_{0}^{t}\int_{\R^d}\int_{\R^d}e^{-2\beta(t-s)}|G(t-s,x-z)-G(t-s,y-z)| \nonumber\\ &\hspace{3cm}\times \ell(z-w)|G(t-s,x-w)-G(t-s,y-w)|\,dz\,dw\,ds.
    \end{align}
Next, we decompose the spatial difference of the heat kernels appearing in $K_{t,\beta}(x,y)$, following Step~3 of the proof of Lemma~\ref{lem:weighted-moments}. We denote this difference by
\begin{equation}
    H_{s,t}^{x,y,z}:= |G(t-s,x-z)-G(t-s,y-z)|.
\end{equation}
For an arbitrary parameter $\lambda \in (0,1)$, we have 
\begin{multline}\label{eq: decompH}
    H_{s,t}^{x,y,z}=|G(t-s,x-z)-G(t-s,y-z)|^{\lambda}|G(t-s,x-z)-G(t-s,y-z)|^{1-\lambda}\\
    \leq |G(t-s,x-z)-G(t-s,y-z)|^{\lambda}\left\{G(t-s,x-z)^{1-\lambda}+G(t-s,y-z)^{1-\lambda}\right\},
\end{multline}
where in the second inequality, we used the sub-additivity of concave functions. Next, the mean value theorem and the gradient bound $(\ref{eq: gradG})$ give $|G(t-s,x-z)-G(t-s,y-z)|\leq C_{d}(t-s)^{-(d+1)/2}|x-y|$, so that \eqref{eq: decompH} yields
\begin{equation}\label{eq: MVTforH}
   H_{s,t}^{x,y,z}
   \leq C_{\lambda}(t-s)^{-\frac{\lambda}{2}(d+1)}|x-y|^{\lambda} 
   \left\{G(t-s,x-z)^{1-\lambda}+G(t-s,y-z)^{1-\lambda}\right\}.
\end{equation}
Applying the estimate $(\ref{eq: MVTforH})$ inside $(\ref{eq: K})$, we thus obtain
\begin{align}\label{eq: K2}
K_{t,\beta}(x,y) &\leq C_{\lambda}\,|x-y|^{2\lambda}\int_{0}^{t}\int_{\R^d}\int_{\R^d}e^{-2\beta(t-s)}(t-s)^{-\lambda(d+1)} \nonumber\\
&\quad \times \{G(t-s,x-z)^{1-\lambda}+G(t-s,y-z)^{1-\lambda}\} \ell(w-z) \nonumber\\
&\quad \times \{G(t-s,x-w)^{1-\lambda}+G(t-s,y-w)^{1-\lambda}\}\,dz\,dw\,ds.
\end{align}
Invoking the relation $(\ref{eq: relation})$ to pass from the space integral to a spectral integral, as we did in \eqref{eq:fourier-space}, we end up with
\begin{align}\label{eq: K2b}
K_{t,\beta}(x,y) &\leq C_{\lambda}\,|x-y|^{2\lambda}\int_{0}^{t}\int_{\R^d}e^{-2\beta(t-s)}(t-s)^{-\lambda(d+1)} \nonumber\\
&\quad \times \bigl|\mathcal{F}\bigl(G(t-s,x-\cdot)^{1-\lambda}+G(t-s,y-\cdot)^{1-\lambda}\bigr)(\xi)\bigr|^{2}\,\specm(d\xi)\,ds.
\end{align}
Since $G(t,\cdot)^{1-\lambda}$ is again a Gaussian, its Fourier transform is given by  a small elaboration of formula~\eqref{eq: FourierPow} with $\hat{\alpha}=1-\lambda$.  That is for a general spatial shift $x$, we have $\mathcal{F}[G(t,x-\cdot)^{1-\lambda}](\xi)=C_{\lambda,d}\,t^{d\lambda/2}e^{i\xi\cdot x}\exp(-\frac{|\xi|^{2}t}{2(1-\lambda)})$. Hence we obtain
\begin{align}\label{eq: FourierComp1}
   \mathcal{F}\left(G(t-s,x-\cdot)^{1-\lambda}+G(t-s,y-\cdot)^{1-\lambda} \right)(\xi) &= C_{\lambda,d}\,(e^{i\xi\cdot x}+e^{i \xi \cdot y}) \nonumber\\
   &\quad \times (t-s)^{d\lambda/2}\exp\!\left(-\frac{|\xi|^{2}(t-s)}{2(1-\lambda)}\right).
\end{align}
Inserting $(\ref{eq: FourierComp1})$ in $(\ref{eq: K2b})$, and using $|e^{i\xi\cdot x}+e^{i\xi\cdot y}|^{2}\leq 4$, this yields
\begin{equation}\label{eq: K3}
    K_{t,\beta}(x,y)\leq C_{\lambda,d}|x-y|^{2\lambda}\int_{0}^{t}\int_{\R^d}e^{-2\beta(t-s)}(t-s)^{-\lambda}\exp\!\left(-\frac{|\xi|^{2}(t-s)}{1-\lambda}\right)\specm(d\xi)\,ds.
\end{equation}
By formula~\eqref{eq: FourierPow} with $\hat{\alpha}=1-\lambda$ and $\tau=t-s$, the Gaussian factor in $(\ref{eq: K3})$ satisfies
\begin{equation}\label{eq: K3b}
    \exp\!\left(-\frac{|\xi|^{2}(t-s)}{1-\lambda}\right) = C_{1-\lambda,d}^{-2}(t-s)^{-d\lambda}\bigl|\mathcal{F}\bigl(G(t-s,\cdot)^{1-\lambda}\bigr)(\xi)\bigr|^{2}.
\end{equation}
Substituting $(\ref{eq: K3b})$ into $(\ref{eq: K3})$ and absorbing the constant $C_{1-\lambda,d}^{-2}$ into $C_{\lambda,d}$, we obtain
\begin{equation}\label{eq: K3c}
    K_{t,\beta}(x,y)\leq C_{\lambda,d}|x-y|^{2\lambda}\int_{0}^{t}e^{-2\beta(t-s)}(t-s)^{-\lambda(d+1)}\int_{\R^d}\bigl|\mathcal{F}\bigl(G(t-s,\cdot)^{1-\lambda}\bigr)(\xi)\bigr|^{2}\specm(d\xi)\,ds.
\end{equation}
Applying Lemma~\ref{lem: spectral} with $\hat{\alpha}=1-\lambda$ and $\tau=t-s$ to bound the spectral integral in $(\ref{eq: K3c})$ by $C_{\lambda,\specm,\eta,\vk}(t-s)^{d\lambda}(1+(t-s)^{-\vk})$, we thus get
\begin{equation}\label{eq: K3d}
    K_{t,\beta}(x,y)\leq C_{\lambda,d,\specm,\eta,\vk}|x-y|^{2\lambda}\int_{0}^{t}e^{-2\beta(t-s)}\bigl[(t-s)^{-\lambda}+(t-s)^{-\lambda-\vk}\bigr]\,ds.
\end{equation}
Observe that the singularity exponents $\la$ and $\la+\vk$ are both in $[0,1)$, since $\vk<1-\la$. Hence one can set $u=t-s$ and apply Lemma~\ref{lem:int-exp} with $\beta$ replaced by $2\beta$ and $\ze=0$. We discover that
\begin{equation} \label{eq: K4}
    K_{t,\beta}(x,y)\leq {C_{\lambda,d,\specm,\eta,\vk}}\left(\frac{1}{\beta^{1-\lambda}}+\frac{1}{\beta^{1-\lambda-\vk}}\right)|x-y|^{2\lambda},
\end{equation}
Let us summarize now the conditions we have used on the parameter $\vk$: the lower bound $\vk>1-\eta$ is the hypothesis of Lemma~\ref{lem: spectral} that ensures the right hand side of~\eqref{eq: spectralBound} is finite, and the upper bound $\vk<1-\lambda$ ensures the exponents $\lambda$ and $\lambda+\vk$ are strictly less than $1$, which is required in order to apply Lemma~\ref{lem:int-exp} to the integrals in $(\ref{eq: K3d})$. Finally, inserting $(\ref{eq: K4})$ in $(\ref{eq: boundSbetaSpace})$ leads to our claim $(\ref{eq: prelimLemmaBound2})$.
\end{proof}

\begin{proof}[Proof of Lemma~\ref{lem: TimeIncrem}]
Throughout this proof, the value of the constant $C$ may change from line to line, but always depends only on the parameters explicitly indicated. Fix $p\geq 1$, $x\in \R^d$, and without loss of generality $t_{1}\leq t_{2}$. As in Step~4 of the proof of Lemma~\ref{lem:weighted-moments}, we decompose the time increment as $\mathcal{S}^{\beta}u(t_{2},x)-\mathcal{S}^{\beta}u(t_{1},x)=I_{1}+I_{2}$, where
\begin{align}\label{eq: I1I2}
I_{1}&:=\int_{t_{1}}^{t_{2}}\int_{\R^d}e^{-\beta(t_{2}-s)}G(t_{2}-s,x-y)\sigma(u(s,y))\,W(ds,dy),\\
I_{2}&:=\int_{0}^{t_{1}}\int_{\R^d}\bigl[e^{-\beta(t_{2}-s)}G(t_{2}-s,x-y)-e^{-\beta(t_{1}-s)}G(t_{1}-s,x-y)\bigr]\sigma(u(s,y))\,W(ds,dy). \nonumber
\end{align}

\smallskip\noindent
\textit{Step 1: Bounding $I_{1}$.}
The same argument as in $(\ref{eq: BDG})$--$(\ref{eq: boundSbeta})$, with the integration domain replaced by $[t_{1},t_{2}]\times\R^d$, gives
\begin{equation}\label{eq: I1bound}
    \mathbf{E}[|I_{1}|^{p}]\leq c_{p}\,M_{\sigma}^{p}\, L_{t_{1},t_{2},\beta}^{p/2},
\end{equation}
where, after applying $(\ref{eq: relation})$ to convert the spatial covariance to a spectral integral,
\begin{equation}\label{eq: L}
    L_{t_{1},t_{2},\beta}:=\int_{t_{1}}^{t_{2}}e^{-2\beta(t_{2}-s)}\widetilde{J}^{\specm}_{s,t_{2}}\,ds,
\end{equation}
with $\widetilde{J}^{\specm}_{s,t_{2}}$ as defined in $(\ref{eq: Jtilde})$. Substituting $u=t_{2}-s$ and applying the bound $(\ref{eq: boundJtilde})$, it is readily checked that
\begin{equation}\label{eq: Lbound1}
    L_{t_{1},t_{2},\beta}\leq C_{\specm,\eta,\vk}\int_{0}^{t_{2}-t_{1}}e^{-2\beta u}\bigl(1+u^{-\vk}\bigr)\,du.
\end{equation}
We can now apply Lemma~\ref{lem:int-exp} with $\beta$ replaced by $2\beta$ and $\ze=\la$ to each integral in~$(\ref{eq: Lbound1})$, the singularity exponent of the lemma taking the values $0$ and $\vk$ respectively (both in $[0,1-\lambda)$ since $\vk<1-\la$). This yields
\begin{equation}\label{eq: Lbound}
    L_{t_{1},t_{2},\beta}\leq C_{\specm,\eta,\vk}\frac{(t_{2}-t_{1})^{\lambda}}{\beta^{1-\lambda-\vk}}.
\end{equation}

\smallskip\noindent
\textit{Step 2: Bounding $I_{2}$.}
We introduce an additional notation for the weighted time increments of the heat kernel. Namely
set 
$
D_{s}^{t_{1},t_{2}}(z):=e^{-\beta(t_{2}-s)}G(t_{2}-s,z)-e^{-\beta(t_{1}-s)}G(t_{1}-s,z).
$ 
The same BDG+Jensen argument as in $(\ref{eq: BDG})$--$(\ref{eq: boundSbeta})$ (see also \eqref{eq: K}), applied to the kernel $D_{s}^{t_{1},t_{2}}(x-\cdot)$, gives
\begin{equation}\label{eq: I2bound}
    \mathbf{E}[|I_{2}|^{p}]\leq c_{p}\,M_{\sigma}^{p}\, P_{t_{1},t_{2},\beta}(x)^{p/2},
\end{equation}
where, after invoking $(\ref{eq:covWiener})$, the expression for $P_{t_{1},t_{2},\beta}(x)$ is spelled out as
\begin{equation}\label{eq: P}
    P_{t_{1},t_{2},\beta}(x):=\int_{0}^{t_{1}}\int_{\R^d}\int_{\R^d}|D_{s}^{t_{1},t_{2}}(x-z)|\,\ell(z-w)\,|D_{s}^{t_{1},t_{2}}(x-w)|\,dz\,dw\,ds.
\end{equation}

\smallskip\noindent
\textit{Step 3: Bounding $P_{t_{1},t_{2},\beta}$.}
We apply to $|D_{s}^{t_{1},t_{2}}|$ the same $\lambda/(1{-}\lambda)$ Hölder decomposition as in~\eqref{eq: decompH}, this time retaining the exponential factors. Specifically, by the sub-additivity of concave functions,
\begin{equation}\label{eq: decompDelta}
    |D_{s}^{t_{1},t_{2}}(z)|\leq |D_{s}^{t_{1},t_{2}}(z)|^{\lambda}\bigl\{e^{-\beta(1-\lambda)(t_{2}-s)}G(t_{2}-s,z)^{1-\lambda}+e^{-\beta(1-\lambda)(t_{1}-s)}G(t_{1}-s,z)^{1-\lambda}\bigr\}.
\end{equation}
Let us first bound the term $|D_{s}^{t_{1},t_{2}}(z)|^{\lambda}$ above. To this aim, we add and subtract $e^{-\beta(t_{1}-s)}G(t_{2}-s,z)$ to write
\begin{equation}\label{eq: MVTforDelta-split}
D_{s}^{t_{1},t_{2}}(z)
= e^{-\beta(t_{2}-s)}\bigl[G(t_{2}-s,z)-G(t_{1}-s,z)\bigr]
  +\bigl[e^{-\beta(t_{2}-s)}-e^{-\beta(t_{1}-s)}\bigr]G(t_{1}-s,z).
\end{equation}
We now bound the two term in the right hand side above separately. 
For the first term in~\eqref{eq: MVTforDelta-split}, we can argue as follows: since $t_{2}\geq t_{1}$ we have $e^{-\beta(t_{2}-s)}\leq e^{-\beta(t_{1}-s)}$, and the mean value theorem together with $(\ref{eq: dtG})$ give
\begin{equation}\label{eq: MVTforDelta-G}
\bigl|G(t_{2}-s,z)-G(t_{1}-s,z)\bigr|\leq C_{d}\,(t_{1}-s)^{-(d+2)/2}\,|t_{2}-t_{1}|.
\end{equation}
For the second term in~\eqref{eq: MVTforDelta-split}, we apply the mean value theorem to $\tau\mapsto e^{-\beta\tau}$ on $[t_{1}-s,t_{2}-s]$. This easily yields
\begin{equation*}
|e^{-\beta(t_{2}-s)}-e^{-\beta(t_{1}-s)}|\leq\beta e^{-\beta(t_{1}-s)}|t_{2}-t_{1}| .
\end{equation*}
Resorting to the inequality $G(t_{1}-s,z)\leq C_{d}(t_{1}-s)^{-d/2}$ from \eqref{eq: heatkernel} and writing $\beta(t_{1}-s)^{-d/2}e^{-\beta(t_{1}-s)}=(t_{1}-s)^{-(d+2)/2}\cdot\beta(t_{1}-s)e^{-\beta(t_{1}-s)}$, the elementary inequality $\sup_{r\geq 0}r\,e^{-r}=1/e$ gives
\begin{equation}\label{eq: MVTforDelta-exp}
\bigl|e^{-\beta(t_{2}-s)}-e^{-\beta(t_{1}-s)}\bigr|\,G(t_{1}-s,z)
\leq C_{d}\,(t_{1}-s)^{-(d+2)/2}|t_{2}-t_{1}|.
\end{equation}
Inserting $(\ref{eq: MVTforDelta-G})$--$(\ref{eq: MVTforDelta-exp})$ into $(\ref{eq: MVTforDelta-split})$ and using $e^{-\beta(t_{1}-s)}\leq 1$, we let the reader check that we end up with
\begin{equation}\label{eq: MVTforDelta}
    |D_{s}^{t_{1},t_{2}}(z)|^{\lambda}\leq C_{\lambda,d}\,
    \frac{|t_{2}-t_{1}|^{\lambda}}{(t_{1}-s)^{\frac{\lambda(d+2)}{2}}} .
\end{equation}

We now plug~\eqref{eq: MVTforDelta} into \eqref{eq: decompDelta}, and then \eqref{eq: decompDelta} into \eqref{eq: P}. This allows to
split $P_{t_{1},t_{2},\beta}$ according to the two Gaussian terms in~\eqref{eq: decompDelta}:
\begin{align}\label{eq: P2}
    P_{t_{1},t_{2},\beta}(x)&\leq C_{\lambda,d}\,|t_{2}-t_{1}|^{2\lambda}\int_{0}^{t_{1}}(t_{1}-s)^{-\lambda(d+2)}\nonumber\\
    &\quad\times\Bigl[e^{-2\beta(1-\lambda)(t_{1}-s)}\int_{\R^{d}}\bigl|\mathcal{F}(G(t_{1}-s,\cdot)^{1-\lambda})(\xi)\bigr|^{2}\specm(d\xi)\nonumber\\
    &\qquad+\,e^{-2\beta(1-\lambda)(t_{2}-s)}\int_{\R^{d}}\bigl|\mathcal{F}(G(t_{2}-s,\cdot)^{1-\lambda})(\xi)\bigr|^{2}\specm(d\xi)\Bigr]ds.
\end{align}
The two terms in the right hand side of~\eqref{eq: decompDelta} require slightly different treatments, that we detail below.

\noindent
\textit{Term involving $G(t_{1}-s,\cdot)$.} Applying Lemma~\ref{lem: spectral} with $\hat{\alpha}=1-\lambda$ and $\tau=t_{1}-s$ directly yields
\begin{equation}\label{eq: P3}
    \int_{\R^{d}}\bigl|\mathcal{F}(G(t_{1}-s,\cdot)^{1-\lambda})(\xi)\bigr|^{2}\specm(d\xi)
    \leq C_{\lambda,d,\specm,\eta,\vk}\,(t_{1}-s)^{d\lambda}\bigl(1+(t_{1}-s)^{-\vk}\bigr).
\end{equation}
Since $d\lambda-\lambda(d+2)=-2\lambda$, the contribution of this term to $P_{t_{1},t_{2},\beta}$ is bounded by
\begin{equation}\label{eq: P4a}
    C|t_{2}-t_{1}|^{2\lambda}\int_{0}^{t_{1}}e^{-2\beta(1-\lambda)(t_{1}-s)}\bigl[(t_{1}-s)^{-2\lambda}+(t_{1}-s)^{-2\lambda-\vk}\bigr]\,ds\leq \frac{C|t_{2}-t_{1}|^{2\lambda}}{\beta^{1-2\lambda-\vk}},
\end{equation}
where we applied Lemma~\ref{lem:int-exp} with $\beta$ replaced by $2\beta(1-\la)$ and $\ze=0$ to each power of $(t_{1}-s)$ separately, and where we recall that we assume $2\la+\vk<1$.

\noindent
\textit{Term involving $G(t_{2}-s,\cdot)$.} Applying Lemma~\ref{lem: spectral} with $\hat{\alpha}=1-\lambda$ and $\tau=t_{2}-s$ to the spectral integral in the $t_{2}$ term of~\eqref{eq: P2} gives
\begin{equation}\label{eq: P3t2}
    \int_{\R^{d}}\bigl|\mathcal{F}\bigl(G(t_{2}-s,\cdot)^{1-\lambda}\bigr)(\xi)\bigr|^{2}\specm(d\xi)
    \leq C\,(t_{2}-s)^{d\lambda}\bigl(1+(t_{2}-s)^{-\vk}\bigr).
\end{equation}
The factor $(t_{2}-s)^{d\lambda}$ in~\eqref{eq: P3t2} has to be handled slightly more carefully. Indeed, since $t_{2}-s\geq t_{1}-s$, it cannot be bounded by the integrable $(t_{1}-s)^{d\lambda}$. Instead, we absorb it into the exponential $e^{-2\beta(1-\lambda)(t_{2}-s)}$ from~\eqref{eq: P2} via the elementary inequality $\sup_{\tau>0}\tau^{a}e^{-c\tau}\leq(a/(ec))^{a}$ with $a=d\lambda$ and $c=\beta(1-\lambda)$:
\begin{equation}\label{eq: absorbPoly}
    (t_{2}-s)^{d\lambda}\,e^{-2\beta(1-\lambda)(t_{2}-s)}\leq C_{\beta}\,e^{-\beta(1-\lambda)(t_{2}-s)},
    \quad C_{\beta}:=\frac{c_{\lambda,d}}{\beta^{d\lambda}},\quad c_{\lambda,d}:=\Bigl(\frac{d\lambda}{e(1-\lambda)}\Bigr)^{d\lambda}.
\end{equation}
We now insert~\eqref{eq: P3t2}--\eqref{eq: absorbPoly} into the $t_{2}$ term of~\eqref{eq: P2}, using $e^{-\beta(1-\lambda)(t_{2}-s)}\leq e^{-\beta(1-\lambda)(t_{1}-s)}$ and $(t_{2}-s)^{-\vk}\leq(t_{1}-s)^{-\vk}$. Since the factor $(t_{2}-s)^{d\lambda}$ has been absorbed into $C_{\beta}$, the power of $(t_{1}-s)$ left from \eqref{eq: P2} is $-\lambda(d+2)$. We get that the contribution of the $t_{2}$ term to $P_{t_{1},t_{2},\beta}$ is bounded by
\begin{equation}\label{eq: P4b}
    C_{\beta}|t_{2}-t_{1}|^{2\lambda}\int_{0}^{t_{1}}e^{-\beta(1-\lambda)(t_{1}-s)}\bigl[(t_{1}-s)^{-\lambda(d+2)}+(t_{1}-s)^{-\lambda(d+2)-\vk}\bigr]\,ds
    \leq \frac{C_{\beta}|t_{2}-t_{1}|^{2\lambda}}{\beta^{1-\lambda(d+2)-\vk}},
\end{equation}
where we applied Lemma~\ref{lem:int-exp} with $\beta$ replaced by $\beta(1-\la)$ and $\ze=0$ to each power of $(t_{1}-s)$ separately. This is where the hypothesis $\vk<1-\lambda(d+2)$ of the lemma is used. Since $C_{\beta}=c_{\lambda,d}\beta^{-d\lambda}$ and $d\lambda+1-\lambda(d+2)-\vk=1-2\lambda-\vk$, the right hand side of~\eqref{eq: P4b} is bounded by $C|t_{2}-t_{1}|^{2\lambda}/\beta^{1-2\lambda-\vk}$.

Summarizing our considerations for this step, we combine $(\ref{eq: P4a})$ and $(\ref{eq: P4b})$ into \eqref{eq: P2}. We conclude
\begin{equation}\label{eq: P4final}
    P_{t_{1},t_{2},\beta}(x)\leq C\frac{|t_{2}-t_{1}|^{2\lambda}}{\beta^{1-2\lambda-\vk}},
\end{equation}
with a constant $C$ independent of $T$.

\smallskip\noindent
\textit{Step 4: Conclusion.}
Inserting $(\ref{eq: Lbound})$ into $(\ref{eq: I1bound})$ and $(\ref{eq: P4final})$ into $(\ref{eq: I2bound})$, we obtain
\begin{equation}\label{eq: conclusion}
    \mathbf{E}[|\mathcal{S}^{\beta}u(t_{2},x)-\mathcal{S}^{\beta}u(t_{1},x)|^{p}]
    \leq C\,M_{\sigma}^{p}\left(\frac{1}{\beta^{1-\lambda-\vk}}+\frac{1}{\beta^{1-2\lambda-\vk}}\right)^{p/2}|t_{2}-t_{1}|^{\lambda p/2},
\end{equation}
which is precisely $(\ref{eq: prelimLemmaBound3})$.
\end{proof}

\section{Measurability}\label{app:measurability}

This section contains the proofs of the two measurability lemmas of Section~\ref{subsec:feller-two}, namely Lemma~\ref{lem:two-spaces} on the Borel structure of the spaces $C_{\nrho}\subset C_{\nrhop}$, and Lemma~\ref{lem:plug-in} on the measure $\mathcal{P}_{t}^{*}\mu$ of Definition~\ref{def: transition}.

\begin{proof}[Proof of Lemma~\ref{lem:two-spaces}]
We prove items \ref{it:ts-inclusion} and \ref{it:ts-borel}, and we leave the easy proof of \ref{it:ts-solution} for the reader.

\smallskip

\noindent
{\it Item \ref{it:ts-inclusion}.} The two weights $w_{0}^{\nrhop}$ and $w_{0}^{\nrho}$ of \eqref{eq: weight1} satisfy $w_{0}^{\nrhop}\leq 2\,w_{0}^{\nrho}$, by item \ref{it:weight-two-exp} of Lemma~\ref{lem:weight-elem}. Owing to \eqref{eq:weight-two-exp}, if $w_{0}^{\nrho}\phi$ vanishes at infinity then so does $w_{0}^{\nrhop}\phi$, and $|\phi|_{C_{\nrhop}}\leq 2|\phi|_{C_{\nrho}}$. This proves item \ref{it:ts-inclusion}.

\smallskip

\noindent
{\it Item \ref{it:ts-borel}.} By item \ref{it:ts-inclusion}, the inclusion map from $C_{\nrho}$ into $C_{\nrhop}$ is continuous, and it is injective. Both spaces are separable Banach spaces (see Remark~\ref{rem: equiv-norms}), hence Polish spaces. By the Lusin-Souslin theorem (see \cite[Theorem~15.1 and Corollary~15.2]{kechris:95:classical}), the image $C_{\nrho}$ of this map is a Borel subset of $C_{\nrhop}$, and the map is a Borel isomorphism onto its image. The latter assertion means that a subset of $C_{\nrho}$ is Borel in $C_{\nrho}$ if and only if it is Borel in $C_{\nrhop}$, which is the second claim of item \ref{it:ts-borel}.
\end{proof}

\begin{proof}[Proof of Lemma~\ref{lem:plug-in}]
Item \ref{it:pi-kernel} is a direct consequence of item \ref{it:ft-meas} of Proposition~\ref{prop:feller-two}: $\phi\mapsto P_{t}(\phi,A)$ is Borel and bounded by $1$, so that \eqref{eq:Pt-star-kernel} is well defined, $\sigma$-additivity in $A$ follows from the monotone convergence theorem, and $\mathcal{P}_{t}^{*}\mu(C_{\nrho})=1$. The first identity in \eqref{eq: Pt} holds for $F=\mathbf{1}_{A}$ by \eqref{eq:Pt-star-kernel} and \eqref{eq:Pt-kernel}, hence for every $F\in B_{b}(C_{\nrho})$ by linearity and monotone convergence.

As for item \ref{it:pi-plug}, let $\Phi_{t}\colon C_{\nrho}\times\Omega\to C_{\nrho}$ be the solution map $(\phi,\omega)\mapsto u(t,\cdot\,;\phi)(\omega)$, which is well defined by item \ref{it:ts-solution} of Lemma~\ref{lem:two-spaces}. We admit that $\Phi_{t}$ is measurable with respect to the product of the Borel $\sigma$-algebra of $C_{\nrho}$ and $\mathcal{F}_{t}$, and that $u(t,\cdot\,;\xi)=\Phi_{t}(\xi,\cdot)$ almost surely (the \emph{plug-in property}). Both facts are standard in this context, see \cite{DPZ96}: the measurability in $(\phi,\omega)$ holds because $\Phi_{t}$ is obtained as the limit of the Picard iterations of Section~\ref{sec:extend-proof}, each of which is measurable in $(\phi,\omega)$, and the plug-in property holds because the same iterations started from $\xi$ converge to $u(\cdot\,;\xi)$, by the uniqueness in Theorem~\ref{thm:extend}. Granting these facts, since $\xi$ is $\mathcal{F}_{0}$-measurable and independent of the noise, hence of $\Phi_{t}(\phi,\cdot)$ for every $\phi$, the Fubini theorem gives, for $F\in B_{b}(C_{\nrho})$,
\begin{equation*}
\mathbf{E}\big[F(u(t,\cdot\,;\xi))\big]=\mathbf{E}\big[F(\Phi_{t}(\xi,\cdot))\big]=\int_{C_{\nrho}}\mathbf{E}\big[F(\Phi_{t}(\phi,\cdot))\big]\,\mu(d\phi)=\int_{C_{\nrho}}\mathcal{P}_{t}F(\phi)\,\mu(d\phi),
\end{equation*}
which is \eqref{eq:plug-in}, owing to item \ref{it:pi-kernel}.
\end{proof}

\bibliographystyle{plain}
\bibliography{ergodic-whole}

\end{document}